\documentclass[11pt,leqno]{article}
\usepackage{amssymb}
\usepackage[mathscr]{eucal}
\usepackage{amsmath,amssymb,latexsym,theorem,bbm}
\usepackage{color}
\usepackage{appendix}

\newcommand{\SC}{\scriptstyle}

\newcommand{\NN}{\mathbb{N}}

\newcommand{\RR}{\mathbb{R}}

\newcommand{\ZZ}{\mathbb{Z}}

\newcommand{\bA}{{\boldsymbol{A}}}

\newcommand{\bB}{{\boldsymbol{B}}}

\newcommand{\bC}{{\boldsymbol{C}}}

\newcommand{\bd}{{\boldsymbol{d}}}
\newcommand{\bD}{{\boldsymbol{D}}}

\newcommand{\bI}{{\boldsymbol{I}}}

\newcommand{\bM}{{\boldsymbol{M}}}

\newcommand{\bQ}{{\boldsymbol{Q}}}
\newcommand{\bS}{{\boldsymbol{S}}}

\newcommand{\bv}{{\boldsymbol{v}}}
\newcommand{\bV}{{\boldsymbol{V}}}

\newcommand{\bZ}{{\boldsymbol{Z}}}
\newcommand{\bU}{{\boldsymbol{U}}}

\newcommand{\bfeta}{{\boldsymbol{\eta}}}
\newcommand{\btheta}{{\boldsymbol{\theta}}}

\newcommand{\bSigma}{{\boldsymbol{\Sigma}}}
\newcommand{\bzero}{{\boldsymbol{0}}}

\newcommand{\cA}{{\mathcal A}}
\newcommand{\cB}{{\mathcal B}}

\newcommand{\cD}{{\mathcal D}}

\newcommand{\cG}{{\mathcal G}}
\newcommand{\cF}{{\mathcal F}}

\newcommand{\cN}{{\mathcal N}}

\newcommand{\cV}{{\mathcal V}}

\newcommand{\cX}{{\mathcal X}}
\newcommand{\cY}{{\mathcal Y}}

\newcommand{\cW}{{\mathcal W}}

\newcommand{\tcY}{\widetilde{\cY}}

\newcommand{\tcV}{\widetilde{\cV}}

\newcommand{\dd}{\mathrm{d}}
\newcommand{\ee}{\mathrm{e}}
\newcommand{\ii}{\mathrm{i}}

\newcommand{\slu}{{\SC\mathrm{lu}}}

\newcommand{\EE}{\operatorname{\mathbb{E}}}
\newcommand{\PP}{\operatorname{\mathbb{P}}}

\newcommand{\ha}{\widehat{a}}
\newcommand{\hb}{\widehat{b}}

\newcommand{\halpha}{\widehat{\alpha}}
\newcommand{\hbeta}{\widehat{\beta}}

\newcommand{\hbtheta}{\widehat{\btheta}}

\newcommand{\hbS}{\widehat{\bS}}

\newcommand{\tX}{\widetilde{X}}
\newcommand{\tY}{\widetilde{Y}}
\newcommand{\tx}{\widetilde{x}}
\newcommand{\ty}{\widetilde{y}}

\newcommand{\tW}{\widetilde{W}}

\renewcommand{\mid}{\,|\,}

\renewcommand{\leq}{\leqslant}
\renewcommand{\geq}{\geqslant}

\newcommand{\stoch}{\stackrel{\PP}{\longrightarrow}}
\newcommand{\distr}{\stackrel{\cD}{\longrightarrow}}
\newcommand{\distre}{\stackrel{\cD}{=}}

\newcommand{\as}{\stackrel{{\mathrm{a.s.}}}{\longrightarrow}}

\newcommand{\nT}{{\lfloor nT\rfloor}}

\newcommand{\proofend}{\hfill\mbox{$\Box$}}

\numberwithin{equation}{section}

\theoremstyle{change} \theorembodyfont{\em}
\newtheorem{Lem}{Lemma.}[section]
\newtheorem{Thm}[Lem]{Theorem.}
\newtheorem{Pro}[Lem]{Proposition.}
\newtheorem{Cor}[Lem]{Corollary.}
\newtheorem{Def}[Lem]{Definition.}

\theorembodyfont{\rm}
\newtheorem{Rem}[Lem]{Remark.}

\begin{document}

\begin{center}
 {\bfseries\Large
  Asymptotic properties of maximum likelihood estimators \\[2mm]
   for Heston models based on continuous time observations} \\[5mm]
 {\sc\large
  M\'aty\'as $\text{Barczy}^{*,\diamond}$ \ and \ Gyula $\text{Pap}^{**}$}
\end{center}

\vskip0.2cm

\noindent
 * Faculty of Informatics, University of Debrecen,
   Pf.~12, H--4010 Debrecen, Hungary.

\noindent
 ** Bolyai Institute, University of Szeged,
    Aradi v\'ertan\'uk tere 1, H--6720 Szeged, Hungary.

\noindent e--mails: barczy.matyas@inf.unideb.hu (M. Barczy),
                    papgy@math.u-szeged.hu (G. Pap).

\noindent $\diamond$ Corresponding author.

\renewcommand{\thefootnote}{}
\footnote{\textit{2010 Mathematics Subject Classifications\/}:
          60H10, 91G70, 60F05, 62F12.}
\footnote{\textit{Key words and phrases\/}:
 Heston model, maximum likelihood estimator.}
\vspace*{0.2cm}
\footnote{The research of M. Barczy and G. Pap was realized in the frames of
 T\'AMOP 4.2.4.\ A/2-11-1-2012-0001 ,,National Excellence Program --
 Elaborating and operating an inland student and researcher personal support
 system''.
The project was subsidized by the European Union and co-financed by the
 European Social Fund.}

\vspace*{-5mm}


\begin{abstract}
We study asymptotic properties of maximum likelihood estimators for Heston
 models based on continuous time observations of the log-price process.
We distinguish three cases: subcritical (also called ergodic), critical and
 supercritical.
In the subcritical case, asymptotic normality is proved for all the parameters,
 while in the critical and supercritical cases, non-standard asymptotic
 behavior is described.
\end{abstract}

\section{Introduction}

Affine processes and especially the Heston model have been frequently applied
 in financial mathematics since they can be well-fitted to financial time
 series, and also due to their computational tractability.
They are characterized by their characteristic function which is exponentially
 affine in the state variable.
A precise mathematical formulation and a complete characterization of regular
 affine processes are due to Duffie et al.\ \cite{DufFilSch}.
A very recent monograph of Baldeaux and Platen \cite{BalPla} gives a detailed
 survey on affine processes and their applications in financial mathematics.

Let us consider a Heston model
 \begin{align}\label{Heston_SDE}
  \begin{cases}
   \dd Y_t = (a - b Y_t) \, \dd t + \sigma_1 \sqrt{Y_t} \, \dd W_t , \\
   \dd X_t = (\alpha - \beta Y_t) \, \dd t
             + \sigma_2  \sqrt{Y_t}
               \bigl(\varrho \, \dd W_t
                     + \sqrt{1 - \varrho^2} \, \dd B_t\bigr) ,
  \end{cases} \qquad t \geq 0 ,
 \end{align}
 where \ $a > 0$, \ $b, \alpha, \beta \in \RR$, \ $\sigma_1 > 0$,
 \ $\sigma_2 > 0$, \ $\varrho \in (-1, 1)$ \ and \ $(W_t, B_t)_{t\geq0}$ \ is a
 2-dimensional standard Wiener process.
In this paper we study maximum likelihood estimator (MLE) of
 \ $(a, b, \alpha, \beta)$ \ based on continuous time observations
 \ $(X_t)_{t\in[0,T]}$ \ with \ $T > 0$,
 \ starting the process \ $(Y,X)$ \ from some known non-random initial value
 \ $(y_0,x_0)\in(0,\infty)\times\RR$.
\ We do not suppose the process \ $(Y_t)_{t\in[0,T]}$ \ being observed, since it
 can be determined using the observations \ $(X_t)_{t\in[0,T]}$, \ see Remark
 \ref{observation}.
We do not estimate the parameters \ $\sigma_1$, \ $\sigma_2$ \ and \ $\varrho$,
 \ since these parameters could ---in principle, at least--- be determined
 (rather than estimated) using the observations \ $(X_t)_{t\in[0,T]}$, \ see
 Remark \ref{Thm_MLE_cons_sigma_rho}.
Further, it will turn out that for the calculation of the MLE of
 \ $(a, b, \alpha, \beta)$, \ one does not need to know the values of the
 parameters \ $\sigma_1>0$, $\sigma_2 >0$, \ and \ $\varrho \in (-1, 1)$,
 \ see \eqref{MLE}.
Note also that \ $(Y_t, X_t)_{t\geq0}$ \ is a 2-dimensional affine diffusion process
 with state space \ $[0, \infty) \times \RR$, \ see Proposition
 \ref{Pro_Heston}.
In the language of financial mathematics, provided that
 \ $\beta = \sigma_2^2/2$, \ one can interpret
 \[
   S_t := \exp\left\{X_t - \alpha + \frac{\sigma_2^2}{2}t\right\}
 \]
 as the asset price, \ $X_t - \alpha + \frac{\sigma_2^2}{2}t$ \ as the
 log-price (log-spot) and \ $\sigma_2 \sqrt{Y_t}$ \ as the volatility of the
 asset price at time \ $t \geq 0$.
\ Indeed, using \eqref{Heston_SDE}, by an application of It\^{o}'s formula, if
 \ $\beta = \sigma_2^2/2$, \ then we have
 \begin{align*}
  \dd S_t = (\alpha + \sigma_2^2/2) S_t \, \dd t
            + \sigma_2  \sqrt{Y_t} S_t
              \bigl(\varrho \, \dd W_t
                    + \sqrt{1 - \varrho^2} \, \dd B_t\bigr) ,
  \qquad t \geq 0 ,
 \end{align*}
 which is Equation (19) in Heston \cite{Hes}.
The squared volatility process \ $(\sigma_2^2 Y_t)_{t\geq0}$ \ is a continuous
 time continuous state branching process with immigration, also called
 Cox--Ingersoll--Ross (CIR) process, first studied by Feller \cite{Fel}.

Parameter estimation for continuous time models has a long history, see, e.g.,
 the monographs of Liptser and Shiryaev \cite[Chapter 17]{LipShiII},
 Kutoyants \cite{Kut} and Bishwal \cite{Bis}.
For estimating continuous time models used in finance, Phillips and Yu
 \cite{PhiYu} gave an overview of maximum likelihood and Gaussian methods.
Since the exact likelihood can be constructed only in special cases
 (e.g., geometric Brownian motion, Ornstein--Uhlenbeck process, CIR process
 and inverse square-root process), much attention has been devoted to the
 development of methods designed to approximate the likelihood.

A\"{\i}t-Sahalia \cite{Ait} provides closed-form expansions for the
 log-likelihood function of multivariate diffusions based on discrete time
 observations.
He proved that, under some conditions, the approximate maximum likelihood
 exists almost surely, and the difference of the approximate and the true
 maximum likelihood converges in probability to 0 as the time interval
 separating observations tends to 0.
The above mentioned closed-form expansions for the Heston model can be found
 in A\"{\i}t-Sahalia and Kimmel \cite[Appendix A.1]{AitKim}.
We note that in S{\o}rensen \cite{Sor} one can find a brief and concise
 summary of the approach of A\"{\i}t-Sahalia.
In fact, S{\o}rensen \cite{Sor} gives a survey of estimation techniques for
 stationary and ergodic (one-dimensional) diffusion processes observed at
 discrete time points.
Besides the above mentioned approach of A\"{\i}t-Sahalia, she recalls
 estimating functions with special emphasis on martingale estimating functions
 and so-called simple estimating functions, together with Bayesian analysis of
 discretely observed diffusion processes.

Azencott and Gadhyan \cite{AzeGad} considered another parametrization of the
 Heston model \eqref{Heston_SDE}, and they investigated only the subcritical
 (also called ergodic) case, i.e., when \ $b > 0$ \ (see Definition
 \ref{Def_criticality}).
They developed an algorithm to estimate the parameters of the Heston model
 based on discrete time observations for the asset price and the volatility.
They supposed that \ $\sigma_2 = 1$ \ and \ $\beta = 1/2$, \ and estimated the
 parameters \ $\sigma_1$ \ and \ $\varrho$ \ as well.
They assumed the time interval separating two consecutive observations also to
 be unknown and used MLE based on Euler and Milstein discretization schemes.
They showed that parameter estimates derived from the Euler scheme using
 constrained optimization of the approximate MLE are strongly consistent.
Note that we obtain results also on the asymptotic behavior of the MLE, and
 not only in the subcritical case.

Hurn et al. \cite{HurLinMcC} developed  a quasi-maximum likelihood procedure
 for estimating the parameters of multi-dimensional diffusions based on
 discrete time obervations by replacing the original transition density by a
 multivariate Gaussian density with first and second moments approximating the
 true moments of the unknown density.
For affine drift and diffusion functions, these moments are exactly those of
 the true transitional density.
As an example, the Heston stochastic volatility model has been analyzed in the
 subcritical case.
However, they did not investigate consistency or asymptotic behavior of their
 estimators.

Recently, Varughese \cite{Var} has studied parameter estimation for time
 inhomogeneous multi-dimensional diffusion processes given by SDEs based on
 discrete time observations.
The likelihood of a diffusion process in question sampled at discrete time
 points has been estimated by a so-called saddlepoint approximation.
In general, the saddlepoint approximation is an algebraic expression based on
 a random variable's cumulant generation function.
In cases where the first few moments of a random variable are known but the
 corresponding probability density is difficult to obtain, the saddlepoint
 approximation to the density can be calculated.
The parameter estimates are taken to be the values that maximize this
 approximate likelihood, which may be estimated by a Markov Chain Monte Carlo
 (MCMC) procedure.
However, the asymptotic properties of the estimators have not been studied.
As an example, the saddlepoint MCMC is used to fit a subcritical Heston model
 to the S\&P 500 and the VIX indices over the period December 2009--November
 2010.

In case of the one-dimensional CIR process \ $Y$, \ the parameter estimation
 of \ $a$ \ and \ $b$ \ goes back to Overbeck and Ryd\'en \cite{OveRyd}
 (conditional least squares estimator (LSE)), Overbeck \cite{Ove} (MLE), and
 see also Bishwal \cite[Example 7.6]{Bis} and the very recent papers of Ben
 Alaya and Kebaier \cite{BenKeb1}, \cite{BenKeb2} (MLE).
We also note that Li and Ma \cite{LiMa} started to investigate the
 asymptotic behaviour of the (weighted) conditional LSE of the drift
 parameters for a CIR model driven by a stable noise (they call it a stable
 CIR model) from some discretely observed low frequency data set.

To the best knowledge of the authors the parameter estimation problem for
 multi-dimensional affine processes has not been tackled so far.
Since affine processes are frequently used in financial mathematics, the
 question of parameter estimation for them needs to be well-investigated.
In Barczy et al.~\cite{BarDorLiPap} we started the discussion with a simple
 non-trivial 2-dimensional affine diffusion process given by the SDE
 \begin{align}\label{2dim_affine}
  \begin{cases}
   \dd Y_t = (a - bY_t) \, \dd t + \sqrt{Y_t} \, \dd W_t , \\
   \dd X_t = (m - \theta X_t) \, \dd t + \sqrt{Y_t} \, \dd B_t ,
  \end{cases}
  \qquad t \geq 0 ,
 \end{align}
 where \ $a > 0$, \ $b, m, \theta \in \RR$, \ $(W_t, B_t)_{t\geq 0}$ \ is a
 2-dimensional standard Wiener process.
Chen and Joslin \cite{CheJos} have found several applications of the model
 \eqref{2dim_affine} in financial mathematics, see their equations (25) and
 (26).
In the special critical case \ $b = 0$, \ $\theta = 0$ \ we described the
 asymptotic behavior of the LSE of \ $(m, \theta)$ \ based on discrete time
 observations \ $X_0, X_1, \ldots, X_n$ \ as \ $n \to \infty$.
\ The description of the asymptotic behavior of the LSE of \ $(m, \theta)$
 \ in the other critical cases \ $b = 0$, \ $\theta > 0$ \ or \ $b > 0$,
 \ $\theta = 0$ \ remained opened.
In Barczy et al. \cite{BarDorLiPap3} we dealt with the same model
 \eqref{2dim_affine} but in the so-called subcritical (ergodic) case:
 \ $b > 0$, \ $\theta > 0$, \ and we considered the MLE of
 \ $(a, b, m, \theta)$ \ and the LSE of \ $(m, \theta)$ \ based on continuous
 time observations.
To carry out the analysis in the subcritical case, we needed to examine the
 question of existence of a unique stationary distribution and ergodicity for
 the model given by \eqref{2dim_affine}.
We solved this problem in a companion paper Barczy et
 al.~\cite{BarDorLiPap2}.

Next, we summarize our results comparing with those of Overbeck \cite{Ove}
 and Ben Alaya and Kebaier \cite{BenKeb1}, \cite{BenKeb2}, and give an
 overview of the structure of the paper.
Section \ref{Prel} is devoted to some preliminaries.
We recall that the SDE \eqref{Heston_SDE} has a pathwise unique strong
 solution and show that it is a regular affine process, see Proposition
 \ref{Pro_Heston}.
We describe the asymptotic behaviour of the first moment of
 \ $(Y_t, X_t)_{t\geq 0}$, \ and, based on it, we introduce a classification of
 Heston processes given by the SDE \eqref{Heston_SDE}, see Proposition
 \ref{Pro_moments} and Definition \ref{Def_criticality}.
Namely, we call \ $(Y_t, X_t)_{t\geq 0}$ \ subcritical, critical or supercritical
 if \ $b > 0$, \ $b = 0$, \ or \ $b < 0$, \ respectively.
We recall a result about existence of a unique stationary distribution and
 ergodicity for the process \ $(Y_t)_{t\geq 0}$ \ given by the first equation in
 \eqref{Heston_SDE} in the subcritical case, see Theorem \ref{Ergodicity}.
From Section \ref{section_EUMLE} we will consider the Heston model
 \eqref{Heston_SDE} with a non-random initial value.
In Section \ref{section_EUMLE} we study the existence and uniqueness of the
 MLE of \ $(a, b, \alpha, \beta)$ \ by giving an explicit formula for this MLE
 as well.
It turned out that the MLE of \ $(a, b)$ \ based on the observations
 \ $(Y_t)_{t\in[0,T]}$ \ for the CIR process \ $Y$ \ is the same as the MLE of
 \ $(a, b)$ \ based on the observations \ $(X_t)_{t\in[0,T]}$ \ for the
 Heston process \ $(Y, X)$ \ given by the SDE \eqref{Heston_SDE}, see formula
 \eqref{MLE} and Overbeck \cite[formula (2.2)]{Ove} or Ben Alaya and Kebaier
 \cite[Section 3.1]{BenKeb2}.

In Section \ref{section_CMLE} we investigate consistency of MLE.
For subcritical Heston models we prove that the MLE of
 \ $(a, b, \alpha, \beta)$ \ is strongly consistent whenever
 \ $a \in \bigl(\frac{\sigma_1^2}{2}, \infty\bigr)$ \ (which is an
 extension of strong consistency of the MLE of \ $(a, b)$ \ proved by Overbeck
 \cite[Theorem 2 (ii)]{Ove}, see Remark \ref{Rem_subcrit_comp}), \ and weakly
 consistent whenever \ $a = \frac{\sigma_1^2}{2}$ \ (which is an
 extension of weak consistency of the MLE of \ $(a, b)$ \ following from
 part 1 of Theorem 7 in Ben Alaya and Kebaier \cite{BenKeb2}, see Remark
 \ref{Rem_subcrit_comp}), see Theorem \ref{Thm_MLE_cons_sub}.
For critical Heston models with
 \ $a \in \bigl(\frac{\sigma_1^2}{2}, \infty\bigr)$, \ we obtain weak
 consistency of the MLE of \ $(a, b, \alpha, \beta)$ \ (as a consequence of
 Theorem \ref{Thm_MLE}), which is an extension of weak consistency of
 the MLE of \ $(a, b)$ \ following from Theorem 6 in Ben Alaya and Kebaier
 \cite{BenKeb2}, see Remark \ref{Rem_comparison5}.
For supercritical Heston models
 \ $a \in \bigl[\frac{\sigma_1^2}{2}, \infty\bigr)$, \ we get strong
 consistency of the MLE of \ $b$, \ see Theorem
 \ref{Thm_MLE_cons_super}, and weak consistency of the MLE of \ $\beta$,
 \ see Theorem \ref{Thm_MLE_super}, and it turns out that the MLE of
 \ $a$ \ and \ $\alpha$ \ is not even weakly consistent, see Corollary
 \ref{Cor_super}.
This is an extension of Overbeck
 \cite[Theorem 2, parts (i) and (v)]{Ove}, see Remark \ref{Rem_comparison1}.

Sections \ref{section_AMLE_subcritical}, \ref{section_AMLE_critical} and
 \ref{section_AMLE_supercritical} are devoted to study asymptotic
 behaviour of the MLE of \ $(a, b, \alpha, \beta)$ \ for subcritical, critical
 and supercritical Heston models, respectively.
In Section \ref{section_AMLE_subcritical} we show that the MLE of
 \ $(a, b, \alpha, \beta)$ \ is asymptotically normal in the subcritical case
 with \ $a \in \bigl(\frac{\sigma_1^2}{2}, \infty\bigr)$, \ which is a
 generalization of the asymptotic normality of the MLE of \ $(a, b)$ \ proved
 by Ben Alaya and Kebaier \cite[Theorem 5]{BenKeb2}, see Remark
 \ref{Rem_comparison2}.
We also show asymptotic normality with random scaling for the MLE
 of \ $(a, b, \alpha, \beta)$ \ generalizing the asymptotic normality with
 random scaling for the MLE of \ $(a, b)$ \ due to Overbeck
 \cite[Theorem 3 (iii)]{Ove}, see Remark \ref{Rem_comparison2}.
In Section \ref{section_AMLE_critical} we describe the asymptotic behaviour of
 the MLE in the critical case with
 \ $a \in \bigl(\frac{\sigma_1^2}{2}, \infty\bigr)$ \ generalizing the
 second part of Theorem 6 in Ben Alaya and Kebaier \cite{BenKeb2}, see Remark
 \ref{Rem_comparison6}.
It turns out that the MLE of \ $a$ \ and \ $\alpha$ \ is asymptotically
 normal, but we have a different limit behaviour for the MLE of \ $b$ \ and
 \ $\beta$, \ see Theorem \ref{Thm_MLE}.
In Theorem \ref{Thm_MLE_random} we incorporate random scaling for the
 MLE of \ $(a, b, \alpha, \beta)$ \ in case of critical Heston models
 generalizing part (ii) of Theorem 3 in Overbeck \cite{Ove}, see Remark
 \ref{Rem_comparison3}.
In Section \ref{section_AMLE_supercritical} for supercritical Heston models
 with \ $a \in \bigl[\frac{\sigma_1^2}{2}, \infty\bigr)$, \ we prove that the
 MLE of \ $a$ \ and \ $\alpha$ \ has a weak limit without any scaling
 (consequently, not weakly consistent, see Corollary \ref{Cor_super}), and the
 appropriately normalized MLE of \ $b$ \ and \ $\beta$ \ has a mixed normal
 limit distribution, which is a generalization of the second part of
 Theorem 3 (i) of Overbeck \cite{Ove}, see Remark \ref{Rem_comparison4}.
We also show asymptotic normality with random scaling for the MLE
 of \ $(b, \beta)$ \ generalizing the asymptotic normality with random
 scaling for the MLE of \ $b$ \ due to Overbeck
 \cite[first part of Theorem 3 (i)]{Ove}, see Remark \ref{Rem_comparison4}.
In the Appendix we recall some limit theorems for continuous local
 martingales for studying asymptotic behaviour of the MLE of
 \ $(a, b, \alpha, \beta)$.

In the proofs, mainly for the critical and supercritical cases, we extensively
 used the following results of Ben Alaya and Kebaier \cite[Propositions 3 and 4]{BenKeb1},
 \cite[Theorems 4 and 6]{BenKeb2}:
 for \ $b>0$ \ and \ $a=\frac{\sigma_1^2}{2}$, \ weak convergence of
 \ $\frac{1}{T^2}\int_0^T \frac{\dd s}{Y_s}$ \ as \ $T\to\infty$; \ for \ $b=0$ \ and
 \ $a>\frac{\sigma_1^2}{2}$, \ the explicit form of the moment generating function of the
 quadruplet \ $\bigl(\log Y_T, Y_T, \int_0^T Y_s\,\dd s, \int_0^T \frac{\dd s}{Y_s}\bigr)$, \ $T>0$;
 \ for \ $b<0$ \ and \ $a\geq \frac{\sigma_1^2}{2}$, \ a representation of the weak limit of
 \ $\bigl(\ee^{bT} Y_T, \int_0^T \frac{\dd s}{Y_s}\bigr)$ \ as \ $T\to\infty$.
\ However, our results are not simple consequences of those of Ben Alaya and Kebaier,
 we will have to find appropriate decompositions of the derived MLEs and then to investigate
 the joint weak convergence of the components via continuity theorem.

In Barczy et al. \cite{BarPapSza},
 we study conditional least squares estimation for the drift
 parameters \ $(a,b,\alpha,\beta)$ \ of the Heston model \eqref{Heston_SDE}
 starting from some known non-random initial value \ $(y_0,x_0)\in[0,\infty)\times\RR$ \
 based on discrete time observations \ $(Y_i,X_i)_{i\in\{1,\ldots,n\}}$, and in the subcritical case
 we describe its asymptotic properties.

Finally, note that Benke and Pap \cite{BenPap} study local asymptotic properties of likelihood ratios
 of the Heston model \eqref{Heston_SDE} under the assumption
 \ $a \in \bigl(\frac{\sigma_1^2}{2}, \infty\bigr)$.
\ Local asymptotic normality has been proved in the subcritical case and for the submodel when
 \ $b = 0$ \ and \ $\beta \in \RR$ \ are known in the critical case.
Moreover, local asymptotic mixed normality has been shown for the submodel when
 \ $a \in \bigl(\frac{\sigma_1^2}{2}, \infty\bigr)$ \ and \ $\alpha \in \RR$ \ are known in the
 supercritical case.
As a consequence, there exist asymptotic minimax bounds for arbitrary estimators in these models,
 the MLE (for the appropriate submodels in the critical and supercritical cases)
 attains this bound for bounded loss function, and the MLE is asymptotically efficient in
 H\'ajek's convolution theorem sense, see Benke and Pap \cite{BenPap}.

\section{Preliminaries}
\label{Prel}

Let \ $\NN$, \ $\ZZ_+$, \ $\RR$, \ $\RR_+$, \ $\RR_{++}$, \ $\RR_-$ \ and
 \ $\RR_{--}$ \ denote the sets of positive integers, non-negative integers,
 real numbers, non-negative real numbers, positive real numbers, non-positive
 real numbers and negative real numbers, respectively.
For \ $x , y \in \RR$, \ we will use the notations \ $x \land y := \min(x, y)$
 \ and \ $x \lor y := \max(x, y)$.
\ By \ $\|x\|$ \ and \ $\|A\|$, \ we denote the Euclidean norm of a vector
 \ $x \in \RR^d$ \ and the induced matrix norm of a matrix
 \ $A \in \RR^{d \times d}$, \ respectively.
By \ $\bI_d \in \RR^{d \times d}$, \ we denote the $d$-dimensional unit matrix.

Let \ $\bigl(\Omega, \cF, \PP\bigr)$ \ be a probability space.
By \ $C^2_c(\RR_+\times\RR, \RR)$ \ and \ $C^{\infty}_c(\RR_+\times\RR, \RR)$,
 \ we denote the set of twice continuously differentiable real-valued
 functions on \ $\RR_+\times\RR$ \ with compact support, and the set of
 infinitely differentiable real-valued functions on \ $\RR_+\times\RR$ \ with
 compact support, respectively.

The next proposition is about the existence and uniqueness of a strong
 solution of the SDE \eqref{Heston_SDE} stating also that \ $(Y,X)$ \ is a
 regular affine process.
Note that these statements for the first equation of \eqref{Heston_SDE}
 are well known.

\begin{Pro}\label{Pro_Heston}
Let \ $(\eta_0, \zeta_0)$ \ be a random vector independent of
 \ $(W_t, B_t)_{t\in\RR_+}$ \ satisfying \ $\PP(\eta_0 \in \RR_+) = 1$.
\ Then for all \ $a \in \RR_{++}$, \ $b, \alpha, \beta \in \RR$,
 \ $\sigma_1, \sigma_2 \in \RR_{++}$, \ $\varrho \in (-1, 1)$, \ there is a
 (pathwise) unique strong solution \ $(Y_t, X_t)_{t\in\RR_+}$ \ of the SDE
 \eqref{Heston_SDE} such that \ $\PP((Y_0, X_0) = (\eta_0, \zeta_0)) = 1$ \ and
 \ $\PP(\text{$Y_t \in \RR_+$ \ for all \ $t \in \RR_+$}) = 1$.
\ Further, for all \ $s, t \in \RR_+$ \ with \ $s \leq t$,
 \begin{align}\label{Solutions}
  \begin{cases}
   Y_t = \ee^{-b(t-s)}
         \left( Y_s
                + a \int_s^t \ee^{-b(s-u)} \, \dd u
                + \sigma_1
                  \int_s^t \ee^{-b(s-u)} \sqrt{Y_u} \, \dd W_u \right) , \\
   X_t = X_s + \int_s^t (\alpha - \beta Y_u) \, \dd u
         + \sigma_2
           \int_s^t
            \sqrt{Y_u} \,
            (\varrho\, \dd  W_u + \sqrt{1 - \varrho^2}\, \dd  B_u) .
  \end{cases}
 \end{align}
Moreover, \ $(Y_t, X_t)_{t\in\RR_+}$ \ is a regular affine process with
 infinitesimal generator
 \begin{align}\label{infgen}
  \begin{split}
   (\cA f)(y,x)
   &= (a - by) f_1'(y, x) + (\alpha - \beta y) f_2'(y, x) \\
   &\quad
      + \frac{1}{2} y
        \big( \sigma_1^2 f_{1,1}''(y, x)
              + 2 \varrho \sigma_1 \sigma_2 f_{1,2}''(y, x)
              + \sigma_2^2 f_{2,2}''(y, x) \big) ,
  \end{split}
 \end{align}
 where \ $(y, x) \in \RR_+\times\RR$, \ $f \in C^2_c(\RR_+\times\RR, \RR)$,
 \ and \ $f_i'$ \ and \ $f_{i,j}''$, \ $i, j \in \{1, 2\}$, \ denote the first
 and second order partial derivatives of \ $f$ \ with respect to its $i$-th,
 and $i$-th and $j$-th variables, respectively.
\end{Pro}

\noindent{\bf Proof.}
By a theorem due to Yamada and Watanabe (see, e.g., Karatzas and Shreve
 \cite[Proposition 5.2.13]{KarShr}), the strong uniqueness holds for the first
 equation in \eqref{Heston_SDE}.
By Ikeda and Watanabe \cite[Example 8.2, page 221]{IkeWat}, there is a
 (pathwise) unique non-negative strong solution \ $(Y_t)_{t\in\RR_+}$ \ of the
 first equation in \eqref{Heston_SDE} with any initial value \ $\eta_0$ \ such
 that \ $\PP(\eta_0 \in \RR_+) = 1$.
\ Clearly, the second equation in \eqref{Solutions} gives the (pathwise)
 unique strong solution \ $(X_t)_{t\in\RR_+}$ \ of the second equation in
 \eqref{Heston_SDE}.
Next, by an application of the It\^o's formula for the process
 \ $(Y_t)_{t\in\RR_+}$, \ we obtain
 \begin{align*}
  \dd(\ee^{bt} Y_t)
  = b \ee^{bt} Y_t \, \dd t + \ee^{bt} \dd Y_t
  = b \ee^{bt} Y_t \, \dd t
    + \ee^{bt} \bigl( (a - b Y_t) \, \dd t
                     + \sigma_1 \sqrt{Y_t} \, \dd W_t \bigr)
  = a \ee^{bt} \, \dd t + \sigma_1 \ee^{bt} \sqrt{Y_t} \, \dd W_t
 \end{align*}
 for all \ $t \in \RR_+$, \ which implies the first equation in
 \eqref{Solutions}.

Now we turn to check that \ $(Y_t, X_t)_{t\in\RR_+}$ \ is an affine
 process with the given infinitesimal generator.
We may and do suppose that the initial value is deterministic, say,
 \ $(Y_0, X_0) = (y_0, x_0) \in \RR_+ \times \RR$, \ since the infinitesimal
 generator of a time homogeneous Markov process does not depend on the initial
 value of the Markov process.
By It\^{o}'s formula, for all \ $f \in C^2_c(\RR_+\times\RR, \RR)$ \ we have
 \begin{align*}
  f(Y_t, X_t)
  &= f(y_0, x_0)
     + \sigma_1 \int_0^t f_1'(Y_s, X_s) \sqrt{Y_s} \, \dd W_s
     + \sigma_2
       \int_0^t
        f_2'(Y_s, X_s) \sqrt{Y_s}
        \, \big(\varrho \dd W_s + \sqrt{1 - \varrho^2} \, \dd B_s \big) \\
  &\quad
     + \int_0^t f_1'(Y_s, X_s) (a - bY_s) \, \dd s
     + \int_0^t f_2'(Y_s, X_s) (\alpha - \beta Y_s) \, \dd s \\
  &\quad
     + \frac{1}{2}
       \left( \sigma_1^2 \int_0^t f_{1,1}''(Y_s, X_s) Y_s \, \dd s
              + 2 \varrho \sigma_1 \sigma_2
                \int_0^t f_{1,2}''(Y_s, X_s) Y_s \, \dd s
              + \sigma_2^2 \int_0^t f_{2,2}''(Y_s, X_s) Y_s \, \dd s \right) \\
  &= f(y_0, x_0) + \int_0^t (\cA f)(Y_s, X_s) \, \dd s + M_t(f) ,
   \qquad t \in \RR_+ ,
 \end{align*}
 where
 \[
   M_t(f):= \sigma_1 \int_0^t f_1'(Y_s, X_s) \sqrt{Y_s} \, \dd W_s
            + \sigma_2
              \int_0^t
               f_2'(Y_s, X_s) \sqrt{Y_s} \,
               \big( \varrho \dd W_s + \sqrt{1 - \varrho^2} \, \dd B_s \big) ,
   \qquad t \in \RR_+ ,
 \]
 and \ $\cA f$ \ is given by \eqref{infgen}.
It is enough to show that \ $(M_t(f))_{t\in\RR_+}$ \ is a local martingale with
 respect to the augmented filtration corresponding to
 \ $(W_t, B_t)_{t\in\RR_+}$ \ and \ $(\eta_0, \zeta_0)$, \ constructed as in
 Karatzas and Shreve \cite[Section 5.2]{KarShr}.
However, it turns out that it is a square integrable martingale with respect
 to this filtration, since
 \begin{align*}
  &\int_0^t \EE((f_1'(Y_s, X_s))^2 Y_s) \, \dd s
   \leq C_1 \int_0^t \EE(Y_s) \, \dd s
   < \infty , \qquad t \in \RR_+ , \\
  &\int_0^t \EE((f_2'(Y_s, X_s))^2 Y_s) \, \dd s
   \leq C_2 \int_0^t \EE(Y_s) \, \dd s
   < \infty , \qquad t \in \RR_+ ,
 \end{align*}
 with some constants \ $C_1, C_2 \in \RR_{++}$, \ where the finiteness of the
 integrals follows by
 \begin{equation}\label{EY}
  \EE(Y_s) = \ee^{-bs} y_0 + a \int_0^s \ee^{-bu} \, \dd u, \qquad s \in \RR_+ ,
 \end{equation}
 see, e.g., Cox et al.\ \cite[Equation (19)]{CoxIngRos} or Jeanblanc et al.\
 \cite[Theorem 6.3.3.1]{JeaYorChe}.

Finally, we check that the transition semigroup \ $(P_t)_{t\in\RR_+}$ \ with
 state space \ $\RR_+\times\RR$ \ corresponding to \ $(Y_t, X_t)_{t\in\RR_+}$
 \ is a regular affine semigroup having infinitesimal generator given by
 \eqref{infgen}.
With the notations of Dawson and Li \cite{DawLi},
 \[
   \left(\begin{bmatrix}
          0 & 0 \\
          0 & 0
         \end{bmatrix} ,
         \frac{1}{2} \bS,
         \begin{bmatrix}
          a \\
          \alpha
         \end{bmatrix} ,
         \begin{bmatrix}
          -b & 0 \\
          -\beta & 0
         \end{bmatrix},
         0, 0\right)
 \]
 is a set of admissible parameters corresponding to the affine process
 \ $(Y_t, X_t)_{t\in\RR_+}$, \ where
 \begin{equation}\label{bSigma}
  \bS := \begin{bmatrix}
          \sigma_1^2 & \varrho \sigma_1 \sigma_2 \\
          \varrho \sigma_1 \sigma_2 &  \sigma_2^2
         \end{bmatrix} .
 \end{equation}

Hence Theorem 2.7 in Duffie et al.\ \cite{DufFilSch} (see also Theorem 6.1 in
 Dawson and Li \cite{DawLi}) yields that for this set of admissible
 parameters, there exists a regular affine semigroup \ $(Q_t)_{t\in\RR_+}$ \ with
 infinitesimal generator given by \eqref{infgen}.
By Theorem 2.7 in Duffie et al.\ \cite{DufFilSch},
 \ $C^\infty_c(\RR_+\times\RR, \RR)$ \ is a core of the infinitesimal generator
 corresponding to the affine semigroup \ $(Q_t)_{t\in\RR_+}$.
\ Since we have checked that the infinitesimal generators corresponding to the
 transition semigroups \ $(P_t)_{t\in\RR_+}$ \ and \ $(Q_t)_{t\in\RR_+}$ \ (defined
 on the Banach space of bounded real-valued functions on \ $\RR_+\times\RR$)
 \ coincide on \ $C^\infty_c(\RR_+\times\RR, \RR)$, \ by the definition of a
 core, we get they coincide on the Banach space of bounded real-valued
 functions on \ $\RR_+\times\RR$.
\ This yields that \ $(Y_t, X_t)_{t\in\RR_+}$ \ is a regular affine process with
 infinitesimal generator \eqref{infgen}.
We also note that we could have used Lemma 10.2 in Duffie et al.\
 \cite{DufFilSch} for concluding that \ $(Y_t, X_t)_{t\in\RR_+}$ \ is a regular
 affine process with infinitesimal generator \eqref{infgen}, since we have
 checked that \ $(M_t(f))_{t\in\RR_+}$ \ is a martingale with respect to the
 filtration \ $(\cF_t)_{t\in\RR_+}$ \ for any
 \ $f \in C^2_c(\RR_+\times\RR, \RR)$.
\proofend

Next we present a result about the first moment of \ $(Y_t, X_t)_{t\in\RR_+}$.
\ We note that Hurn et al.\ \cite[Equation (23)]{HurLinMcC} derived the same
 formula for the expectation of \ $(Y_t, X_t)$, \ $t \in \RR_+$, \ by a
 different method.
Note also that the formula for \ $\EE(Y_t)$, \ $t \in \RR_+$, \ is well known.

\begin{Pro}\label{Pro_moments}
Let \ $(Y_t, X_t)_{t\in\RR_+}$ \ be the unique strong solution of the SDE
 \eqref{Heston_SDE} satisfying \ $\PP(Y_0 \in \RR_+) = 1$ \ and
 \ $\EE(Y_0) < \infty$, \ $\EE(|X_0|) < \infty$.
\ Then
 \begin{align*}
  \begin{bmatrix}
   \EE(Y_t) \\
   \EE(X_t) \\
  \end{bmatrix}
  = \begin{bmatrix}
     \ee^{-bt} & 0 \\
     - \beta \int_0^t \ee^{-bu} \, \dd u & 1 \\
    \end{bmatrix}
    \begin{bmatrix}
     \EE(Y_0) \\
     \EE(X_0) \\
    \end{bmatrix}
    + \begin{bmatrix}
       \int_0^t \ee^{-bu} \, \dd u & 0 \\
       - \beta \int_0^t \left(\int_0^u \ee^{-bv} \, \dd v \right) \dd u & t \\
      \end{bmatrix}
      \begin{bmatrix}
       a \\
       \alpha \\
      \end{bmatrix} ,
  \qquad t \in \RR_+ .
 \end{align*}
Consequently, if \ $b \in \RR_{++}$, \ then
 \[
   \lim_{t\to\infty} \EE(Y_t) = \frac{a}{b} , \qquad
   \lim_{t\to\infty} t^{-1} \EE(X_t) = \alpha - \frac{\beta a}{b} ,
 \]
 if \ $b = 0$, \ then
 \[
   \lim_{t\to\infty} t^{-1} \EE(Y_t) = a , \qquad
   \lim_{t\to\infty} t^{-2} \EE(X_t) = - \frac{1}{2} \beta a ,
 \]
 if \ $b \in \RR_{--}$, \ then
 \[
   \lim_{t\to\infty} \ee^{bt} \EE(Y_t) = \EE(Y_0) - \frac{a}{b} , \qquad
   \lim_{t\to\infty} \ee^{bt} \EE(X_t)
   = \frac{\beta}{b} \EE(Y_0) - \frac{\beta a}{b^2} .
 \]
\end{Pro}

\noindent{\bf Proof.}
It is sufficient to prove the statement in the case when
 \ $(Y_0, X_0) = (y_0, x_0)$ \ with an arbitrary
 \ $(y_0, x_0) \in \RR_{++} \times \RR$, \ since then the statement of the
 proposition follows by the law of total expectation.

The formula for \ $\EE(Y_t)$, \ $t \in \RR_+$, \ can be found, e.g., in Cox et
 al.\ \cite[Equation (19)]{CoxIngRos} or Jeanblanc et al.\
 \cite[Theorem 6.3.3.1]{JeaYorChe}.
Next we observe that
 \begin{align}\label{help1}
  \left( \int_0^t
          \sqrt{Y_u} \,
          \dd(\varrho W_u + \sqrt{1 - \varrho^2} B_u) \right)_{t\in\RR_+}
 \end{align}
 is a square integrable martingale, since
 \begin{align*}
  \EE\left[\left(\int_0^t
                  \sqrt{Y_u} \,
                  \dd(\varrho W_u + \sqrt{1 - \varrho^2} B_u) \right)^2\right]
  = \int_0^t \EE(Y_u) \, \dd u
  < \infty ,
 \end{align*}
 where the finiteness of the integral follows from \eqref{EY}.

Taking expectations of both sides of the second equation in \eqref{Solutions}
 and using the martingale property of the process in \eqref{help1}, we have
 \begin{align*}
  \EE(X_t)
  &= x_0 + \int_0^t (\alpha - \beta \EE(Y_u)) \, \dd u \\
  &= x_0 + \alpha t
     - \beta
       \int_0^t
        \left( \ee^{-bu} y_0 + a \int_0^u \ee^{-bv} \, \dd v \right)
        \dd u \\
  &= x_0
     - \beta y_0 \int_0^t \ee^{-bu} \, \dd u
     + \alpha t
     - \beta a \int_0^t \left( \int_0^u \ee^{-bv} \, \dd v \right) \dd u
 \end{align*}
 for all \ $t \in \RR_+$.

Further, if \ $b \in \RR_{++}$, \ then
 \begin{align*}
  &\lim_{t\to\infty} \EE(Y_t)
   = \lim_{t\to\infty} \Big( \ee^{-bt} y_0 - \frac{a}{b} (\ee^{-bt} -1) \Big)
   = \frac{a}{b},\\
  &\lim_{t\to\infty} t^{-1}\EE(X_t)
   = \lim_{t\to\infty}
      \left( \frac{x_0}{t} + \frac{\beta}{b} y_0 \frac{\ee^{-bt}-1}{t}
             + \alpha
             + \frac{\beta a}{bt}\left( \frac{\ee^{-bt} - 1}{-b} - t\right)
      \right)
   = \alpha - \frac{\beta a}{b}.
 \end{align*}
If \ $b = 0$, \ then
 \begin{align*}
  &\lim_{t\to\infty} t^{-1} \EE(Y_t)
   = \lim_{t\to\infty} t^{-1} (y_0 + at) = a , \\
  &\lim_{t\to\infty} t^{-2}\EE(X_t)
   = \lim_{t\to\infty}
      \left( \frac{x_0}{t^2} - \frac{\beta y_0}{t}
             + \frac{\alpha}{t} - \frac{\beta a}{2} \right)
   = - \frac{\beta a}{2} .
 \end{align*}
If \ $b \in \RR_{--}$, \ then
 \begin{align*}
  \lim_{t\to\infty} \ee^{bt}\EE(Y_t)
  &= \lim_{t\to\infty} \left(y_0 + \frac{a}{b} (\ee^{bt} - 1) \right)
   = y_0 - \frac{a}{b} , \\
  \lim_{t\to\infty} \ee^{bt}\EE(X_t)
  &= x_0 \lim_{t\to\infty} \ee^{bt}
     + \frac{\beta}{b} y_0 \lim_{t\to\infty} (1 - \ee^{bt})
     + \alpha \lim_{t\to\infty} t\ee^{bt}
     + \frac{\beta a}{b}
       \lim_{t\to\infty} \left(\frac{1-\ee^{bt}}{-b} - t\ee^{bt}\right) \\
   &= \frac{\beta}{b} y_0 - \frac{\beta a}{b^2} . \\[-16mm]
 \end{align*}
\proofend

Based on the asymptotic behavior of the expectations \ $(\EE(Y_t), \EE(X_t))$
 \ as \ $t \to \infty$, \ we introduce a classification of Heston processes
 given by the SDE \eqref{Heston_SDE}.

\begin{Def}\label{Def_criticality}
Let \ $(Y_t, X_t)_{t\in\RR_+}$ \ be the unique strong solution of the SDE
 \eqref{Heston_SDE} satisfying \ $\PP(Y_0 \in \RR_+) = 1$.
\ We call \ $(Y_t, X_t)_{t\in\RR_+}$ \ subcritical, critical or supercritical if
 \ $b \in \RR_{++}$, \ $b = 0$ \ or \ $b \in \RR_{--}$, \ respectively.
\end{Def}

In the sequel \ $\stoch$, \ $\distr$ \ and \ $\as$ \ will denote convergence
 in probability, in distribution and almost surely, respectively.

The following result states the existence of a unique stationary distribution
 and the ergodicity for the process \ $(Y_t)_{t\in\RR_+}$ \ given by the first
 equation in \eqref{Heston_SDE} in the subcritical case, see, e.g., Feller
 \cite{Fel}, Cox et al.\ \cite[Equation (20)]{CoxIngRos}, Li and Ma
 \cite[Theorem 2.6]{LiMa} or Theorem 3.1 with \ $\alpha = 2$ \ and
 Theorem 4.1 in Barczy et al.\ \cite{BarDorLiPap2}.

\begin{Thm}\label{Ergodicity}
Let \ $a, b, \sigma_1 \in \RR_{++}$.
\ Let \ $(Y_t)_{t\in\RR_+}$ \ be the unique strong solution of the first equation
 of the SDE \eqref{Heston_SDE} satisfying \ $\PP(Y_0 \in \RR_+) = 1$.
 \renewcommand{\labelenumi}{{\rm(\roman{enumi})}}
 \begin{enumerate}
  \item
   Then \ $Y_t \distr Y_\infty$ \ as \ $t \to \infty$, \ and the distribution of
    \ $Y_\infty$ \ is given by
    \begin{align}\label{Laplace}
     \EE(\ee^{-\lambda Y_\infty})
     = \left(1 + \frac{\sigma_1^2}{2b} \lambda \right)^{-2a/\sigma_1^2} ,
     \qquad \lambda \in \RR_+ ,
    \end{align}
    i.e., \ $Y_\infty$ \ has Gamma distribution with parameters
    \ $2a / \sigma_1^2$ \ and \ $2b / \sigma_1^2$, \ hence
    \[
      \EE(Y_\infty^\kappa)
      = \frac{\Gamma\left(\frac{2a}{\sigma_1^2} + \kappa\right)}
             {\left(\frac{2b}{\sigma_1^2}\right)^\kappa
              \Gamma\left(\frac{2a}{\sigma_1^2}\right)} , \qquad
      \kappa \in \left( -\frac{2a}{\sigma_1^2}, \infty \right) .
    \]
   Especially, \ $\EE(Y_\infty) = \frac{a}{b}$.
   \ Further, if \ $a \in \left( \frac{\sigma_1^2}{2}, \infty \right)$,
    \ then \ $\EE\left(\frac{1}{Y_\infty}\right) = \frac{2b}{2a - \sigma_1^2}$.
 \item
  Supposing that the random initial value \ $Y_0$ \ has the same distribution
   as \ $Y_\infty$, \ the process \ $(Y_t)_{t\in\RR_+}$ \ is strictly stationary.
 \item
  For all Borel measurable functions \ $f : \RR \to \RR$ \ such that
   \ $\EE(|f(Y_\infty)|) < \infty$, \ we have
   \begin{equation}\label{ergodic}
    \frac{1}{T} \int_0^T f(Y_s) \, \dd s \as \EE(f(Y_\infty)) \qquad
    \text{as \ $T \to \infty$.}
   \end{equation}
\end{enumerate}
\end{Thm}

In the next remark we explain why we suppose only that the process \ $X$
 \ is observed.

\begin{Rem}\label{observation}
If \ $a \in \RR_{++}$, \ $b, \alpha, \beta \in \RR$,
 \ $\sigma_1, \sigma_2 \in \RR_{++}$, \ $\varrho \in (-1, 1)$, \ and
 \ $(Y_0, X_0) = (y_0, x_0) \in \RR_{++} \times \RR$, \ then, by the SDE
 \eqref{Heston_SDE},
 \[
   \langle X \rangle_t = \sigma_2^2 \int_0^t Y_s \, \dd s , \qquad t \in \RR_+ .
 \]
By Theorems I.4.47 a) and I.4.52 in Jacod and Shiryaev \cite{JSh},
 \[
   \sum_{i=1}^{\lfloor nt\rfloor} (X_{\frac{i}{n}} - X_{\frac{i-1}{n}})^2
   \stoch \langle X \rangle_t \qquad \text{as \ $n \to \infty$,} \quad
   t \in \RR_+ .
 \]
This convergence holds almost surely along a suitable subsequence, the members
 of this sequence are measurable functions of \ $(X_s)_{s\in[0,t]}$, \ hence,
 using Theorems 4.2.2 and 4.2.8 in Dudley \cite{Dud}, we obtain that
 \ $\langle X \rangle_t = \sigma_2^2 \int_0^t Y_s \, \dd s$ \ is a measurable
 function of \ $(X_s)_{s\in[0,t]}$.
\ Moreover,
 \begin{equation}\label{QV}
   \frac{\langle X \rangle_{t+h} - \langle X \rangle_t}{h}
   = \frac{\sigma_2^2}{h} \int_t^{t+h} Y_s \, \dd s
   \as \sigma_2^2 Y_t  \qquad \text{as \ $h \to 0$,} \quad t \in \RR_+ ,
 \end{equation}
 since \ $Y$ \ has almost surely continuous sample paths.
In particular,
 \[
   \frac{\langle X \rangle_h}{hy_0}
   = \frac{\sigma_2^2}{hy_0} \int_0^h Y_s \, \dd s
   \as \sigma_2^2 \frac{Y_0}{y_0} = \sigma_2^2  \qquad \text{as \ $h \to 0$,}
 \]
 hence, for any fixed \ $T > 0$, \ $\sigma_2^2$ \ is a measurable function of
 \ $(X_s)_{s\in[0,T]}$, \ i.e., it can be determined from a sample
 \ $(X_s)_{s\in[0,T]}$ \ (provided that \ $(Y,X)$ \ starts from some known non-random
  initial value \ $(y_0,x_0)\in(0,\infty)\times\RR$).
\ However, we also point out that this measurable function remains abstract.
Consequently, by \eqref{QV}, for all \ $t \in [0, T]$, \ $Y_t$ \ is a
 measurable function of \ $(X_s)_{s\in[0,T]}$, \ i.e., it can be determined from
 a sample \ $(X_s)_{s\in[0,T]}$ \
 (provided that \ $(Y,X)$ \ starts from some known non-random
  initial value \ $(y_0,x_0)\in(0,\infty)\times\RR$).
\ Finally, we note that the sample size \ $T$ \ is fixed above, and
 it is enough to know any short sample \ $(X_s)_{s\in[0,T]}$ \
 to carry out the above calculations.
\proofend
\end{Rem}

Next we give statistics for the parameters \ $\sigma_1$, $\sigma_2$ \ and
 \ $\varrho$ \ using continuous time observations \ $(X_t)_{t\in[0,T]}$ \ with
 some \ $T > 0$
 \ (provided that \ $(Y,X)$ \ starts from some known non-random
  initial value \ $(y_0,x_0)\in(0,\infty)\times\RR$).
\ Due to this result we do not consider the estimation of these parameters,
 they are supposed to be known.

\begin{Rem}\label{Thm_MLE_cons_sigma_rho}
If \ $a \in \RR_{++}$, \ $b, \alpha, \beta \in \RR$,
 \ $\sigma_1, \sigma_2 \in \RR_{++}$, \ $\varrho \in (-1, 1)$, \ and
 \ $(Y_0, X_0) = (y_0, x_0) \in \RR_{++} \times \RR$, \ then for all \ $T > 0$,
 \[
   \bS = \frac{1}{\int_0^T Y_s \, \dd s}
            \begin{bmatrix}
             \langle Y \rangle_T & \langle Y, X \rangle_T \\
             \langle Y, X \rangle_T & \langle X \rangle_T
            \end{bmatrix}
            =: \hbS_T \qquad \text{almost surely,}
 \]
 where \ $(\langle Y, X \rangle_t)_{t\in\RR_+}$ \ denotes the quadratic
 cross-variation process of \ $Y$ \ and \ $X$, \ since, by the SDE
 \eqref{Heston_SDE},
 \begin{gather*}
  \langle Y \rangle_T
   = \sigma_1^2 \int_0^T Y_s \, \dd s , \qquad
  \langle X \rangle_T
  = \sigma_2^2 \int_0^T Y_s \, \dd s , \qquad
  \langle Y, X \rangle_T
  = \varrho \sigma_1 \sigma_2 \int_0^T Y_s \, \dd s .
 \end{gather*}
Here \ $\hbS_T$ \ is a statistic, i.e., there exists a measurable
 function \ $\Xi : C([0,T], \RR) \to \RR^{2\times2}$ \ such that
 \ $\hbS_T = \Xi((X_s)_{s\in[0,T]})$, \ where \ $C([0,T], \RR)$ \ denotes the
 space of continuous real-valued functions defined on \ $[0,T]$, \ since
 \begin{equation}\label{sigma_rho}
   \frac{1}{\frac{1}{n} \sum_{i=1}^{\lfloor nT\rfloor} Y_{\frac{i-1}{n}}}
   \sum_{i=1}^{\lfloor nT\rfloor}
    \begin{bmatrix}
     Y_{\frac{i}{n}} - Y_{\frac{i-1}{n}} \\
     X_{\frac{i}{n}} - X_{\frac{i-1}{n}}
    \end{bmatrix}
    \begin{bmatrix}
     Y_{\frac{i}{n}} - Y_{\frac{i-1}{n}} \\
     X_{\frac{i}{n}} - X_{\frac{i-1}{n}}
    \end{bmatrix}^\top
   \stoch
   \hbS_T \qquad \text{as \ $n \to \infty$,}
 \end{equation}
 where \ $\lfloor x \rfloor$ \ denotes the integer part of a real number
 \ $x \in \RR$, \ the convergence in \eqref{sigma_rho} holds almost
 surely along a suitable subsequence, by Remark \ref{observation}, the members
 of the sequence in \eqref{sigma_rho} are measurable functions of
 \ $(X_s)_{s\in[0,T]}$, \ and one can use Theorems 4.2.2 and 4.2.8 in Dudley \cite{Dud}.
\ Next we prove \eqref{sigma_rho}.
By Theorems I.4.47 a) and I.4.52 in Jacod and Shiryaev \cite{JSh},
 \begin{gather*}
  \sum_{i=1}^{\lfloor nT\rfloor} (Y_{\frac{i}{n}} - Y_{\frac{i-1}{n}})^2
  \stoch \langle Y \rangle_T , \qquad
  \sum_{i=1}^{\lfloor nT\rfloor} (X_{\frac{i}{n}} - X_{\frac{i-1}{n}})^2
  \stoch \langle X \rangle_T , \\
  \sum_{i=1}^{\lfloor nT\rfloor}
   (Y_{\frac{i}{n}} - Y_{\frac{i-1}{n}}) (X_{\frac{i}{n}} - X_{\frac{i-1}{n}})
  \stoch \langle Y, X \rangle_T
 \end{gather*}
 as \ $n \to \infty$.
\ Consequently,
 \[
   \sum_{i=1}^{\lfloor nT\rfloor}
    \begin{bmatrix}
     Y_{\frac{i}{n}} - Y_{\frac{i-1}{n}} \\
     X_{\frac{i}{n}} - X_{\frac{i-1}{n}}
    \end{bmatrix}
    \begin{bmatrix}
     Y_{\frac{i}{n}} - Y_{\frac{i-1}{n}} \\
     X_{\frac{i}{n}} - X_{\frac{i-1}{n}}
    \end{bmatrix}^\top
   \stoch
   \left( \int_0^T Y_s \, \dd s \right)
   \hbS_T
 \]
 as \ $n \to \infty$, \ see, e.g., van der Vaart
 \cite[Theorem 2.7, part (vi)]{Vaart}.
Moreover,
 \[
   \frac{1}{n} \sum_{i=1}^{\lfloor nT\rfloor} Y_{\frac{i-1}{n}}
   \as \int_0^T Y_s \, \dd s  \qquad \text{as \ $n \to \infty$}
 \]
 since \ $Y$ \ has almost surely continuous sample paths.
Here \ $\PP\bigl(\int_0^T Y_s \, \dd s \in \RR_{++}\bigr) = 1$.
\ Indeed, if \ $\omega \in \Omega$ \ is such that
 \ $[0, T] \ni s \mapsto Y_s(\omega)$ \ is continuous and
 \ $Y_t(\omega) \in \RR_+$ \ for all \ $t \in \RR_+$, \ then we have
 \ $\int_0^T Y_s(\omega)\,\dd s = 0$ \ if and only if \ $Y_s(\omega) = 0$ \ for
 all \ $s \in [0, T]$.
\ Using the method of the proof of Theorem 3.1 in Barczy et.\ al
 \cite{BarDorLiPap}, we get \ $\PP(\int_0^T Y_s = 0) = 0$, \ as desired.
Hence \eqref{sigma_rho} follows by properties of convergence in probability.
\proofend
\end{Rem}

\section{Existence and uniqueness of MLE}
\label{section_EUMLE}

From this section, we will consider the Heston model \eqref{Heston_SDE} with a
 known non-random initial value \ $(y_0, x_0) \in \RR_{++} \times \RR$, \ and we
 equip \ $\bigl(\Omega, \cF, \PP\bigr)$ \ with the augmented filtration
 \ $(\cF_t)_{t\in\RR_+}$ \ corresponding to \ $(W_t, B_t)_{t\in\RR_+}$,
 \ constructed as in Karatzas and Shreve \cite[Section 5.2]{KarShr}.
Note that \ $(\cF_t)_{t\in\RR_+}$ \ satisfies the usual conditions, i.e.,
 the filtration \ $(\cF_t)_{t\in\RR_+}$ \ is right-continuous and
 \ $\cF_0$ \ contains all the $\PP$-null sets in \ $\cF$.

Let \ $\PP_{(Y,X)}$ \ denote the probability measure induced by
 \ $(Y_t, X_t)_{t\in\RR_+}$ \ on the measurable space
 \ $(C(\RR_+, \RR_+\times\RR), \cB(C(\RR_+, \RR_+\times\RR)))$ \ endowed with
 the natural filtration \ $(\cG_t)_{t\in\RR_+}$, \ given by
 \ $\cG_t := \varphi_t^{-1}(\cB(C(\RR_+, \RR_+\times\RR)))$, \ $t \in \RR_+$,
 \ where \ $\varphi_t : C(\RR_+, \RR_+\times\RR) \to C(\RR_+, \RR_+\times\RR)$
 \ is the mapping \ $\varphi_t(f)(s) := f(t \land s)$, \ $s, t \in \RR_+$,
 \ $f\in C(\RR_+, \RR_+\times\RR)$.
\ Here \ $C(\RR_+, \RR_+\times\RR)$ \ denotes the set of
 \ $\RR_+\times\RR$-valued continuous functions defined on \ $\RR_+$, \ and
 \ $\cB(C(\RR_+, \RR_+\times\RR))$ \ is the Borel \ $\sigma$-algebra on it.
Further, for all \ $T \in \RR_{++}$, \ let
 \ $\PP_{(Y,X),T} := \PP_{(Y,X)}|_{\cG_T}$ \ be the restriction of \ $\PP_{(Y,X)}$
 \ to \ $\cG_T$.

\begin{Lem}\label{RN}
Let \ $a \in \bigl[ \frac{\sigma_1^2}{2}, \infty \bigr)$,
 \ $b, \alpha, \beta \in \RR$, \ $\sigma_1, \sigma_2 \in \RR_{++}$, \ and
 \ $\varrho \in (-1, 1)$.
\ Let \ $(Y_t, X_t)_{t\in\RR_+}$ \ and \ $(\tY_t, \tX_t)_{t\in\RR_+}$ \ be the
 unique strong solutions of the SDE \eqref{Heston_SDE} with initial values
 \ $(y_0, x_0) \in \RR_{++} \times \RR$,
 \ $(\ty_0, \tx_0) \in \RR_{++} \times \RR$ \ such that
 \ $(y_0, x_0) = (\ty_0, \tx_0)$, \ corresponding to the parameters
 \ $(a, b, \alpha, \beta, \sigma_1, \sigma_2, \varrho)$ \ and
 \ $(\sigma_1^2, 0, 0, 0, \sigma_1, \sigma_2, \varrho)$, \ respectively.
Then for all \ $T \in \RR_{++}$, \ the measures \ $\PP_{(Y,X),T}$ \ and
 \ $\PP_{(\tY,\tX),T}$ \ are absolutely continuous with respect to each other,
 and the Radon--Nikodym derivative of \ $\PP_{(Y,X),T}$ \ with respect to
 \ $\PP_{(\tY,\tX),T}$ \ (the so called likelihood ratio) takes the form
 \begin{align*}
  L^{(Y,X),(\tY,\tX)}_T \bigl((Y_s, X_s)_{s\in[0,T]}\bigr)
  =\exp\Bigg\{&\int_0^T
                \frac{1}{Y_s}
                \begin{bmatrix}
                 a - b Y_s - \sigma_1^2 \\
                 \alpha - \beta Y_s
                \end{bmatrix}^\top
                \bS^{-1}
                \begin{bmatrix}
                 \dd Y_s \\
                 \dd X_s
                \end{bmatrix} \\
              &-\frac{1}{2}
                \int_0^T
                 \frac{1}{Y_s}
                \begin{bmatrix}
                 a - b Y_s - \sigma_1^2 \\
                 \alpha - \beta Y_s
                \end{bmatrix}^\top
                \bS^{-1}
                \begin{bmatrix}
                 a - b Y_s + \sigma_1^2 \\
                 \alpha - \beta Y_s
                \end{bmatrix}
                \dd s \Bigg\} ,
 \end{align*}
 where \ $\bS$ \ is defined in \eqref{bSigma}.
\end{Lem}

\noindent{\bf Proof.}
First note that the SDE \eqref{Heston_SDE} can be written in the matrix form
 \begin{align}\label{Heston_SDE_matrix}
  \begin{bmatrix} \dd Y_t \\ \dd X_t \end{bmatrix}
  = \left( \begin{bmatrix} - b & 0 \\ - \beta & 0 \end{bmatrix}
           \begin{bmatrix} Y_t \\ X_t \end{bmatrix}
           + \begin{bmatrix} a \\ \alpha \end{bmatrix} \right) \dd t
    + \sqrt{Y_t}
      \begin{bmatrix}
       \sigma_1 & 0 \\
       \sigma_2 \varrho & \sigma_2 \sqrt{1 - \varrho^2}
      \end{bmatrix}
      \begin{bmatrix}
       \dd W_t \\
       \dd B_t
      \end{bmatrix} ,
  \qquad t \in \RR_+ .
 \end{align}
Note also that under the condition
 \ $a \in \bigl[ \frac{\sigma_1^2}{2}, \infty \bigr)$, \ we have
 \ $\PP(\text{$Y_t \in \RR_{++}$ \ for all \ $t \in \RR_+$}) = 1$,
 \ see, e.g., page 442 in Revuz and Yor \cite{RevYor}.

We intend to use formula (7.139) in Section 7.6.4 of Liptser and Shiryaev
 \cite{LipShiI}.
We have to check their condition (7.137) which takes the form
 \begin{equation}\label{LS_cond}
  \PP\Bigg( \int_0^T
             \frac{1}{Y_s}
             \begin{bmatrix}
              a - b Y_s \\
              \alpha - \beta Y_s
             \end{bmatrix}^\top
              \bS^{-1}
             \begin{bmatrix}
              a - b Y_s \\
              \alpha - \beta Y_s
             \end{bmatrix}
           + \frac{1}{Y_s}
             \begin{bmatrix}
               \sigma_1^2 \\
                0
             \end{bmatrix}^\top
             \bS^{-1}
             \begin{bmatrix}
              \sigma_1^2 \\
              0
             \end{bmatrix}
             \dd s
            < \infty \Bigg)
  = 1, \qquad \forall \; T \in \RR_+.
 \end{equation}
Here note that the matrix \ $\bS$ \ is invertible, since
 \ $\sigma_1, \sigma_2 \in \RR_{++}$ \ and \ $\varrho \in (-1, 1)$.
\ Since \ $Y$ \ has continuous sample paths almost surely, condition
 \eqref{LS_cond} holds if
 \begin{align}\label{int_Y}
  \PP\left(\int_0^T \frac{1}{Y_s} \, \dd s < \infty\right) = 1 \qquad
  \text{for all \ $T \in \RR_+$.}
 \end{align}
Since \ $Y$ \ has continuous sample paths almost surely and
 \ $\PP(Y_t \in \RR_{++}, \, \forall \, t \in \RR_+) = 1$,
 \ we have \ $\PP(\inf_{t\in[0,T]} Y_t \in \RR_{++}) = 1$ \ for all
 \ $T \in \RR_+$, \ which yields \eqref{int_Y}.
Note that under the condition
 \ $a \in \bigl[ \frac{\sigma_1^2}{2}, \infty \bigr)$,
 \ Theorems 1 and 3 in Ben Alaya and Kebaier \cite{BenKeb1} also imply
 \eqref{int_Y}.
Applying formula (7.139) in Section 7.6.4 of Liptser and Shiryaev
 \cite{LipShiI} we obtain the statement.

We call the attention that conditions (4.110) and (4.111) are also required
 for Section 7.6.4 in Liptser and Shiryaev \cite{LipShiI}, but the Lipschitz
 condition (4.110) in Liptser and Shiryaev \cite{LipShiI} does not hold for
 the SDE \eqref{Heston_SDE}.
However, we can use formula (7.139) in Liptser and Shiryaev \cite{LipShiI},
 since they use their conditions (4.110) and (4.111) only in order to ensure
 that the SDE they consider in Section 7.6.4 has a unique strong solution
 (see, the proof of Theorem 7.19 in Liptser and Shiryaev \cite{LipShiI}).
By Proposition \ref{Pro_Heston}, under the conditions of the present lemma,
 there is a (pathwise) unique strong solution of the SDE \eqref{Heston_SDE}.
\proofend

By Lemma \ref{RN}, under its conditions the log-likelihood function satisfies
 \begin{align*}
  &(1 - \varrho^2)
    \log L^{(Y,X),(\tY,\tX)}_T \bigl((Y_s, X_s)_{s\in[0,T]}\bigr) \\
   &=\int_0^T
      \frac{1}{Y_s}
      \left[ \left( \frac{a - b Y_s - \sigma_1^2}{\sigma_1^2}
                    - \frac{\varrho (\alpha - \beta Y_s)}
                           {\sigma_1 \sigma_2} \right) \dd Y_s
             + \left( - \frac{\varrho (a - b Y_s - \sigma_1^2)}
                             {\sigma_1 \sigma_2}
                      + \frac{\alpha - \beta Y_s}{\sigma_2^2} \right)
               \dd X_s \right] \\
   &\quad
     - \frac{1}{2}
       \int_0^T
        \frac{1}{Y_s}
        \left[ \frac{(a - b Y_s)^2 - \sigma_1^4}{\sigma_1^2}
               - \frac{2 \varrho (a - b Y_s) (\alpha - \beta Y_s)}
                      {\sigma_1 \sigma_2}
               + \frac{(\alpha - \beta Y_s)^2}{\sigma_2^2} \right] \dd s \\
   &= a \int_0^T
         \left( \frac{\dd Y_s}{\sigma_1^2 Y_s}
                - \frac{\varrho \, \dd X_s}{\sigma_1 \sigma_2 Y_s} \right)
      + b \int_0^T
           \left( - \frac{\dd Y_s}{\sigma_1^2}
                  + \frac{\varrho \, \dd X_s}{\sigma_1 \sigma_2} \right) \\
   &\quad
      + \alpha \int_0^T
           \left( - \frac{\varrho \, \dd Y_s}{\sigma_1 \sigma_2 Y_s}
                  + \frac{\dd X_s}{\sigma_2^2 Y_s} \right)
      + \beta
        \int_0^T
         \left( \frac{\varrho \, \dd Y_s}{\sigma_1 \sigma_2}
                - \frac{\dd X_s}{\sigma_2^2} \right) \\
   &\quad
      - \frac{1}{2} a^2 \int_0^T \frac{\dd s}{\sigma_1^2 Y_s}
      + a b \int_0^T \frac{\dd s}{\sigma_1^2}
      - \frac{1}{2} b^2 \int_0^T \frac{Y_s \, \dd s}{\sigma_1^2}
      - \frac{1}{2} \alpha^2 \int_0^T \frac{\dd s}{\sigma_2^2 Y_s}
      + \alpha \beta \int_0^T \frac{\dd s}{\sigma_2^2}
      - \frac{1}{2} \beta^2 \int_0^T \frac{Y_s \, \dd s}{\sigma_2^2} \\
   &\quad
      + a \alpha \int_0^T \frac{\varrho \, \dd s}{\sigma_1 \sigma_2 Y_s}
      - (b \alpha + a \beta)
        \int_0^T \frac{\varrho \, \dd s}{\sigma_1 \sigma_2}
      + b \beta \int_0^T \frac{\varrho Y_s \, \dd s}{\sigma_1 \sigma_2}
      - \int_0^T \frac{\dd Y_s}{Y_s}
      + \int_0^T \frac{\varrho \sigma_1 \, \dd X_s}{\sigma_2 Y_s}
      + \frac{1}{2} \int_0^T \frac{\sigma_1^2 \, \dd s}{Y_s} \\
   &= \btheta^\top \! \bd_T - \frac{1}{2} \btheta^\top \! \bA_T \, \btheta
      - \int_0^T \frac{\dd Y_s}{Y_s}
      + \int_0^T \frac{\varrho \sigma_1 \, \dd X_s}{\sigma_2 Y_s}
      + \frac{1}{2} \int_0^T \frac{\sigma_1^2 \, \dd s}{Y_s} ,
 \end{align*}
 where\vspace*{-3mm}
 \[
   \btheta
   := \begin{bmatrix}
       a \\
       b \\
       \alpha \\
       \beta
      \end{bmatrix} , \qquad
   \bd_T
   := \bd_T^{(\sigma_1, \sigma_2, \varrho)}\bigl((Y_s, X_s)_{s\in[0,T]}\bigr)
   := \begin{bmatrix}
       \int_0^T
        \left( \frac{\dd Y_s}{\sigma_1^2 Y_s}
               - \frac{\varrho \, \dd X_s}
                      {\sigma_1 \sigma_2 Y_s} \right) \\[2mm]
       \int_0^T
        \left( - \frac{\dd Y_s}{\sigma_1^2}
               + \frac{\varrho \, \dd X_s}{\sigma_1 \sigma_2} \right) \\[2mm]
       \int_0^T
        \left( - \frac{\varrho \, \dd Y_s}{\sigma_1 \sigma_2 Y_s}
               + \frac{\dd X_s}{\sigma_2^2 Y_s} \right) \\[2mm]
       \int_0^T
        \left( \frac{\varrho \, \dd Y_s}{\sigma_1 \sigma_2}
               - \frac{\dd X_s}{\sigma_2^2} \right)
      \end{bmatrix} ,
 \]
 \[
   \bA_T
   := \bA_T^{(\sigma_1, \sigma_2, \varrho)}\bigl((Y_s, X_s)_{s\in[0,T]}\bigr)
   := \begin{bmatrix}
       \int_0^T \frac{\dd s}{\sigma_1^2 Y_s}
        & -\int_0^T \frac{\dd s}{\sigma_1^2}
        & -\int_0^T \frac{\varrho \, \dd s}{\sigma_1 \sigma_2 Y_s}
        & \int_0^T \frac{\varrho \, \dd s}{\sigma_1 \sigma_2} \\[2mm]
       -\int_0^T \frac{\dd s}{\sigma_1^2}
        & \int_0^T \frac{Y_s \, \dd s}{\sigma_1^2}
        & \int_0^T \frac{\varrho \, \dd s}{\sigma_1 \sigma_2}
        & -\int_0^T \frac{\varrho Y_s \, \dd s}{\sigma_1 \sigma_2} \\[2mm]
       -\int_0^T \frac{\varrho \, \dd s}{\sigma_1 \sigma_2 Y_s}
        & \int_0^T \frac{\varrho \, \dd s}{\sigma_1 \sigma_2}
        & \int_0^T \frac{\dd s}{\sigma_2^2 Y_s}
        & -\int_0^T \frac{\dd s}{\sigma_2^2} \\[2mm]
       \int_0^T \frac{\varrho \, \dd s}{\sigma_1 \sigma_2}
        & -\int_0^T \frac{\varrho Y_s \, \dd s}{\sigma_1 \sigma_2}
        & -\int_0^T \frac{\dd s}{\sigma_2^2}
        & \int_0^T \frac{Y_s \, \dd s}{\sigma_2^2}
      \end{bmatrix} .
 \]
If we fix \ $\sigma_1, \sigma_2 \in \RR_{++}$, \ $\varrho \in (-1, 1)$, \ the
 initial value \ $(y_0, x_0)\in\RR_{++}\times\RR$, \ and
 \ $T \in \RR_{++}$, \ then the probability measures \ $\PP_{(Y,X),T}$ \ induced
 by \ $(Y_t, X_t)_{t\in\RR_+}$ \ corresponding to the parameters
 \ $(a, b, \alpha, \beta, \sigma_1, \sigma_2, \varrho)$, \ where
 \ $a \in \bigl[ \frac{\sigma_1^2}{2}, \infty \bigr)$,
 \ $b, \alpha, \beta \in \RR$, \  are absolutely continuous with respect to
 each other.
Hence it does not matter which measure is taken as a reference measure for
 defining the MLE (we have chosen the measure corresponding to the parameters
 \ $(\sigma_1^2, 0, 0, 0, \sigma_1, \sigma_2, \varrho)$).
\ For more details, see, e.g., Liptser and Shiryaev \cite[page 35]{LipShiI}.

The random symmetric matrix \ $\bA_T$ \ can be written as a Kronecker product
 of a deterministic symmetric matrix and a random symmetric matrix, namely,
 \[
   \bA_T = \begin{bmatrix}
            \frac{1}{\sigma_1^2} & - \frac{\varrho}{\sigma_1 \sigma_2} \\
            - \frac{\varrho}{\sigma_1 \sigma_2} & \frac{1}{\sigma_2^2}
           \end{bmatrix}
           \otimes
           \begin{bmatrix}
            \int_0^T \frac{\dd s}{Y_s} & -\int_0^T 1 \, \dd s \\[2mm]
            -\int_0^T 1 \, \dd s & \int_0^T Y_s \, \dd s
           \end{bmatrix} .
 \]
The first matrix is strictly positive definite.
The second matrix is strictly positive definite if and only if
 \ $\int_0^T Y_s \, \dd s \int_0^T \frac{\dd s}{Y_s} > T^2$.
\ The eigenvalues of \ $\bA_T$ \ coincides with the products of the
 eigenvalues of the two matrices in question (taking into account their
 multiplicities), hence the matrix \ $\bA_T$ \ is strictly positive definite
 if and only if \ $\int_0^T Y_s \, \dd s \int_0^T \frac{\dd s}{Y_s} > T^2$,
 \ and in this case the inverse \ $\bA_T^{-1}$ \ has the form (applying the
 identity \ $(\bA\otimes\bB)^{-1} = \bA^{-1}\otimes \bB^{-1}$)
 \[
   \bA_T^{-1}
   =\begin{bmatrix}
     \frac{1}{\sigma_1^2} & - \frac{\varrho}{\sigma_1 \sigma_2} \\
     - \frac{\varrho}{\sigma_1 \sigma_2} & \frac{1}{\sigma_2^2}
    \end{bmatrix}^{-1}
    \otimes
    \begin{bmatrix}
     \int_0^T \frac{\dd s}{Y_s} & -T \\[2mm]
     -T & \int_0^T Y_s \, \dd s
    \end{bmatrix}^{-1}
    = \frac{\bS \otimes
            \begin{bmatrix}
             \int_0^T Y_s \, \dd s  & T \\[2mm]
             T & \int_0^T \frac{\dd s}{Y_s}
            \end{bmatrix}}
           {(1 - \varrho^2)
            \left( \int_0^T Y_s \, \dd s \int_0^T \frac{\dd s}{Y_s}
                   - T^2 \right)} .
 \]
Hence we have
 \begin{multline*}
  2 (1 - \varrho^2) \log L^{(Y,X),(\tY,\tX)}_T \bigl((Y_s, X_s)_{s\in[0,T]}\bigr) \\
  = - (\btheta - \bA_T^{-1} \bd_T)^\top \bA_T (\btheta - \bA_T^{-1} \bd_T)
    + \bd_T^\top \bA_T^{-1} \bd_T
    - 2\int_0^T \frac{\dd Y_s}{Y_s}
    + 2\int_0^T \frac{\varrho \sigma_1 \, \dd X_s}{\sigma_2 Y_s}
    +  \int_0^T \frac{\sigma_1^2 \, \dd s}{Y_s} ,
 \end{multline*}
 provided that \ $\int_0^T Y_s \, \dd s \int_0^T \frac{\dd s}{Y_s} > T^2$.
\ Recall that \ $\sigma_1, \sigma_2 \in \RR_{++}$ \ and \ $\varrho \in (-1, 1)$
 \ are supposed to be known.
Then maximizing
 \ $(1 - \varrho^2) \log L^{(Y,X),(\tY,\tX)}_T \bigl((Y_s, X_s)_{s\in[0,T]}\bigr)$
 \ in \ $(a, b, \alpha, \beta) \in \RR^4$ \ gives the MLE of
 \ $(a, b, \alpha, \beta)$ \ based on the observations \ $(X_t)_{t\in[0,T]}$
 \ having the form
 \begin{align*}
  \hbtheta_T
  = \begin{bmatrix}
     \ha_T \\
     \hb_T \\
     \halpha_T \\
     \hbeta_T
    \end{bmatrix}
  = \bA_T^{-1} \bd_T ,
 \end{align*}
 provided that \ $\int_0^T Y_s \, \dd s \int_0^T \frac{\dd s}{Y_s} > T^2$.
\ The random vector \ $\bd_T$ \ can be expressed as
 \[
   \bd_T = \begin{bmatrix}
            \frac{1}{\sigma_1^2} \\
            - \frac{\varrho}{\sigma_1 \sigma_2}
           \end{bmatrix}
           \otimes
           \begin{bmatrix}
            \int_0^T \frac{\dd Y_s}{Y_s} \\[1mm]
            -\int_0^T \dd Y_s
           \end{bmatrix}
           +
           \begin{bmatrix}
            - \frac{\varrho}{\sigma_1 \sigma_2} \\
            \frac{1}{\sigma_2^2}
           \end{bmatrix}
           \otimes
           \begin{bmatrix}
            \int_0^T \frac{\dd X_s}{Y_s} \\[1mm]
            -\int_0^T \dd X_s
           \end{bmatrix} .
 \]
Applying the identity
 \ $(\bA \otimes \bB) (\bC \otimes \bD) = (\bA \bC) \otimes (\bB \bD)$, \ we
 can calculate
 \begin{align*}
  &\left( \bS \otimes
          \begin{bmatrix}
           \int_0^T Y_s \, \dd s  & T \\[1mm]
           T & \int_0^T \frac{\dd s}{Y_s}
          \end{bmatrix} \right)
   \bd_T \\
  &\quad
   = \left( \bS
            \begin{bmatrix}
             \frac{1}{\sigma_1^2} \\
             - \frac{\varrho}{\sigma_1 \sigma_2}
            \end{bmatrix} \right)
     \otimes
     \left( \begin{bmatrix}
             \int_0^T Y_s \, \dd s  & T \\[1mm]
             T & \int_0^T \frac{\dd s}{Y_s}
            \end{bmatrix}
            \begin{bmatrix}
             \int_0^T \frac{\dd Y_s}{Y_s} \\[1mm]
             -\int_0^T \dd Y_s
            \end{bmatrix} \right) \\
  &\quad\quad
     + \left( \bS
              \begin{bmatrix}
               - \frac{\varrho}{\sigma_1 \sigma_2} \\
               \frac{1}{\sigma_2^2}
              \end{bmatrix} \right)
       \otimes
       \left( \begin{bmatrix}
               \int_0^T Y_s \, \dd s  & T \\[1mm]
               T & \int_0^T \frac{\dd s}{Y_s}
              \end{bmatrix}
              \begin{bmatrix}
               \int_0^T \frac{\dd X_s}{Y_s} \\[1mm]
               -\int_0^T \dd X_s
              \end{bmatrix} \right) \\
  &\quad
   = \begin{bmatrix}
      1 - \varrho^2 \\
      0
     \end{bmatrix}
     \otimes
     \begin{bmatrix}
      \int_0^T Y_s \, \dd s \int_0^T \frac{\dd Y_s}{Y_s}
       - T (Y_T - y_0) \\[2mm]
      T \int_0^T \frac{\dd Y_s}{Y_s} - (Y_T - y_0) \int_0^T \frac{\dd s}{Y_s}
     \end{bmatrix}
     + \begin{bmatrix}
        0 \\
        1 - \varrho^2
       \end{bmatrix}
       \otimes
       \begin{bmatrix}
        \int_0^T Y_s \, \dd s \int_0^T \frac{\dd X_s}{Y_s}
         - T (X_T - x_0) \\[2mm]
        T \int_0^T \frac{\dd X_s}{Y_s} - (X_T - x_0) \int_0^T \frac{\dd s}{Y_s}
       \end{bmatrix} .
 \end{align*}
Consequently, we obtain
 \begin{align}\label{MLE}
  \begin{bmatrix}
   \ha_T \\
   \hb_T \\
   \halpha_T \\
   \hbeta_T
  \end{bmatrix}
  &= \frac{1}{\int_0^T Y_s \, \dd s \int_0^T \frac{\dd s}{Y_s} - T^2}
     \begin{bmatrix}
      \int_0^T Y_s \, \dd s \int_0^T \frac{\dd Y_s}{Y_s} - T (Y_T - y_0) \\[1mm]
      T \int_0^T \frac{\dd Y_s}{Y_s}
       - (Y_T - y_0) \int_0^T \frac{\dd s}{Y_s} \\[1mm]
      \int_0^T Y_s \, \dd s \int_0^T \frac{\dd X_s}{Y_s} - T (X_T - x_0) \\[1mm]
      T \int_0^T \frac{\dd X_s}{Y_s} - (X_T - x_0) \int_0^T \frac{\dd s}{Y_s}
     \end{bmatrix} ,
 \end{align}
 provided that \ $\int_0^T Y_s \, \dd s \int_0^T \frac{\dd s}{Y_s} > T^2$.
\ In fact, it turned out that for the calculation of the MLE of
 \ $(a, b, \alpha, \beta)$, \ one does not need to know the values of the
 parameters \ $\sigma_1, \sigma_2 \in \RR_{++}$ \ and \ $\varrho \in (-1, 1)$.
\ Note that the MLE of \ $(a, b)$ \ based on the observations
 \ $(X_t)_{t\in[0,T]}$ \ for the Heston model \ $(Y,X)$ \ is the same
 as the MLE of \ $(a, b)$ \ based
 on the observations \ $(Y_t)_{t\in[0,T]}$ \ for the CIR process \ $Y$, \
 see, e.g., Overbeck \cite[formula (2.2)]{Ove} or Ben Alaya and Kebaier
 \cite[Section 3.1]{BenKeb2}.

In the next remark we point out that the MLE \eqref{MLE} of
 \ $(a, b, \alpha, \beta)$ \ can be approximated using discrete time
 observations for \ $X$, \ which can be reassuring for practical applications,
 where data in continuous record is not available.

\begin{Rem}
For the stochastic integrals \ $\int_0^T \frac{\dd X_s}{Y_s}$
 \ and \ $\int_0^T \frac{\dd Y_s}{Y_s}$ \ in \eqref{MLE}, we have
 \begin{equation}\label{measurability}
   \sum_{i=1}^\nT
    \frac{X_{\frac{i}{n}} - X_{\frac{i-1}{n}}}
         {Y_{\frac{i-1}{n}}}
   \stoch
   \int_0^T \frac{\dd X_s}{Y_s} \qquad \text{and} \qquad
   \sum_{i=1}^\nT
    \frac{Y_{\frac{i}{n}} - Y_{\frac{i-1}{n}}}
         {Y_{\frac{i-1}{n}}}
   \stoch
   \int_0^T \frac{\dd Y_s}{Y_s} \qquad
   \text{as \ $n \to \infty$,}
 \end{equation}
 following from Proposition I.4.44 in Jacod and Shiryaev \cite{JSh} with the
 Riemann sequence of deterministic subdivisions
 \ $\left(\frac{i}{n} \land T\right)_{i\in\NN}$, \ $n \in \NN$.
\ Thus, there exist measurable functions
 \ $\Phi, \Psi : C([0,T], \RR) \to \RR$
 \ such that \ $\int_0^T \frac{\dd X_s}{Y_s} = \Phi((X_s)_{s\in[0,T]})$
 \ and \ $\int_0^T \frac{\dd Y_s}{Y_s} = \Psi((X_s)_{s\in[0,T]})$, \ since
 the convergences in \eqref{measurability} hold almost surely along suitable
 subsequences, by Remark \ref{observation}, the members of both sequences in
 \eqref{measurability} are measurable functions of \ $(X_s)_{s\in[0,T]}$,
 \ and one can use Theorems 4.2.2 and 4.2.8 in Dudley \cite{Dud}.
\ Moreover, since \ $Y$ \ has continuous sample paths almost surely,
 \begin{align*}
   \frac{1}{n} \sum_{i=1}^\nT Y_{\frac{i-1}{n}}
     \as \int_0^T Y_s \, \dd s \qquad \text{as \ $n \to \infty$,}\qquad\text{and}\qquad
        \frac{1}{n} \sum_{i=1}^\nT \frac{1}{Y_{\frac{i-1}{n}}}
           \as \int_0^T \frac{\dd s}{Y_s} \qquad \text{as \ $n \to \infty$,}
 \end{align*}
 hence the right hand side of \eqref{MLE} is a measurable function of
 \ $(X_s)_{s\in[0,T]}$, \ i.e., it is a statistic.
Further, one can define a sequence \ $(\hbtheta_{T,n})_{n\in\NN}$ \ of estimators
 of \ $\btheta = (a,b,\alpha,\beta)^\top$ \ based only on the discrete time
 observations \ $(Y_{\frac{i}{n}}, X_{\frac{i}{n}})_{i\in\{1,\ldots,\nT\}}$ \ such that
 \ $\hbtheta_{T,n} \stoch \hbtheta_T$ \ as \ $n \to \infty$.
\ This is also called infill asymptotics.
This phenomenon is similar to the approximate MLE, used by A\"it-Sahalia
 \cite{Ait}, as discussed in the Introduction.
\proofend
\end{Rem}

Using the SDE \eqref{Heston_SDE} one can check that
 \begin{align}\nonumber
  \begin{bmatrix}
   \ha_T - a \\
   \hb_T - b \\
   \halpha_T - \alpha \\
   \hbeta_T - \beta
  \end{bmatrix}
  &= \frac{1}{\int_0^T Y_s \, \dd s \int_0^T \frac{\dd s}{Y_s} - T^2}
     \begin{bmatrix}
      \int_0^T Y_s \, \dd s \int_0^T \frac{\dd Y_s}{Y_s} - T (Y_T - y_0)
       - a \int_0^T Y_s \, \dd s \int_0^T \frac{\dd s}{Y_s} + a T^2 \\[1mm]
      T \int_0^T \frac{\dd Y_s}{Y_s} - (Y_T - y_0) \int_0^T \frac{\dd s}{Y_s}
       - b \int_0^T Y_s \, \dd s \int_0^T \frac{\dd s}{Y_s} + b T^2 \\[1mm]
      \int_0^T Y_s \, \dd s \int_0^T \frac{\dd X_s}{Y_s} - T (X_T - x_0)
       - \alpha \int_0^T Y_s \, \dd s \int_0^T \frac{\dd s}{Y_s}
       + \alpha T^2 \\[1mm]
      T \int_0^T \frac{\dd X_s}{Y_s} - (X_T - x_0) \int_0^T \frac{\dd s}{Y_s}
       - \beta \int_0^T Y_s \, \dd s \int_0^T \frac{\dd s}{Y_s} + \beta T^2
     \end{bmatrix} \\
  &= \frac{1}{\int_0^T Y_s \, \dd s \int_0^T \frac{\dd s}{Y_s} - T^2}
     \begin{bmatrix}
      \sigma_1 \int_0^T Y_s \, \dd s \int_0^T \frac{\dd W_s}{\sqrt{Y_s}}
       - \sigma_1 T \int_0^T \sqrt{Y_s} \, \dd W_s \\[1mm]
      \sigma_1 T \int_0^T \frac{\dd W_s}{\sqrt{Y_s}}
       - \sigma_1 \int_0^T \frac{\dd s}{Y_s}
         \int_0^T \sqrt{Y_s} \, \dd W_s \\[1mm]
      \sigma_2 \int_0^T Y_s \, \dd s \int_0^T \frac{\dd \tW_s}{\sqrt{Y_s}}
      - \sigma_2 T \int_0^T \sqrt{Y_s} \, \dd \tW_s \\[1mm]
      \sigma_2 T \int_0^T \frac{\dd \tW_s}{\sqrt{Y_s}}
      - \sigma_2 \int_0^T \frac{\dd s}{Y_s} \int_0^T \sqrt{Y_s} \, \dd \tW_s
     \end{bmatrix} , \label{MLE-}
 \end{align}
 provided that \ $\int_0^T Y_s \, \dd s \int_0^T \frac{\dd s}{Y_s} > T^2$,
 \ where the process
 \[
   \tW_s := \varrho W_s + \sqrt{1 - \varrho^2} B_s , \qquad s \in \RR_+ ,
 \]
 is a standard Wiener process.

The next lemma is about the existence of
 \ $\bigl( \ha_T, \hb_T, \halpha_T, \hbeta_T \bigr)$.

\begin{Lem}\label{LEMMA_MLE_exist}
If \ $a \in \bigl[ \frac{\sigma_1^2}{2}, \infty \bigr)$, \ $b \in \RR$,
 \ $\sigma_1 \in \RR_{++}$, \ and \ $Y_0 = y_0 \in \RR_{++}$, \ then
 \begin{align}\label{positive}
  \PP\left( \int_0^T Y_s \, \dd s \int_0^T \frac{1}{Y_s} \, \dd s > T^2 \right)
  = 1 \qquad \text{for all \ $T \in \RR_{++}$,}
 \end{align}
 and hence, supposing also that \ $\alpha, \beta \in \RR$,
 \ $\sigma_2 \in\RR_{++}$, \ $\varrho \in(-1, 1)$, \ and \ $X_0=x_0\in\RR$,
 \ there exists a unique MLE
 \ $\bigl( \ha_T, \hb_T, \halpha_T, \hbeta_T \bigr)$ \ for all
 \ $T \in \RR_{++}$.
\end{Lem}

\noindent{\bf Proof.}
First note that
 \ $\PP(\text{$Y_t \in \RR_{++}$ \ for all \ $t \in \RR_+$}) = 1$ \ as it was
 detailed in the proof of Lemma \ref{RN}.
We have \ $\PP(\int_0^T Y_s \, \dd s < \infty) = 1$ \ for all \ $T \in \RR_+$,
 \ since \ $Y$ \ has continuous trajectories almost surely, and further,
 \ $\PP(\int_0^T \frac{1}{Y_s} \, \dd s < \infty) = 1$ \ by \eqref{int_Y}.
For each \ $T \in \RR_{++}$, \ put
 \[
   A_T := \left\{ \omega \in \Omega
                  : \text{$t \mapsto Y_t(\omega)$ \ is continuous and positive
                           on \ $[0, T]$} \right\} .
 \]
Then \ $A_T \in \cF$, \ $\PP(A_T) = 1$, \ and for all \ $\omega \in A_T$, \ by
 the Cauchy--Schwarz's inequality, we have
 \[
   \int_0^T Y_s(\omega) \, \dd s \int_0^T \frac{1}{Y_s(\omega)} \, \dd s
   \in [T^2, \infty) ,
 \]
 and
 \ $\int_0^T Y_s(\omega) \, \dd s \int_0^T \frac{1}{Y_s(\omega)} \, \dd s
    = T^2$
 \ if and only if
 \ $K_T(\omega) Y_s(\omega) = \frac{L_T(\omega)}{Y_s(\omega)}$ \ for almost
 every \ $s \in [0, T]$ \ with some \ $K_T(\omega), L_T(\omega)\in \RR_+$
 \ satisfying \ $K_T(\omega)^2 + L_T(\omega)^2 \in \RR_{++}$.
\ Clearly, \ $K_T(\omega) = 0$ \ would imply \ $L_T(\omega) = 0$, \ thus
 \ $K_T(\omega) \ne 0$ \ and
 \ $Y_s(\omega) = \left(\frac{L_T(\omega)}{K_T(\omega)}\right)^{1/2}$ \ for
 almost every \ $s \in [0, T]$.
\ Hence \ $Y_s(\omega) = y_0$ \ for all \ $s \in [0, T]$ \ if
 \ $\omega \in A_T$ \ and
 \ $\int_0^T Y_s(\omega) \, \dd s \int_0^T \frac{1}{Y_s(\omega)} \, \dd s = T^2$.
\ Since the quadratic variation of a deterministic process is the identically zero process, the quadratic variation process \ $(\langle Y \rangle_t)_{t\in[0,T]}$
 \ of \ $(Y_t)_{t\in[0,T]}$ \ should be identically zero on the event
  \[
    A_T\cap \left\{ \omega\in\Omega :
                      \int_0^T Y_s(\omega) \, \dd s \int_0^T \frac{1}{Y_s(\omega)} \, \dd s = T^2
             \right\}.
  \]
Since \ $\langle Y \rangle_t = \sigma_1^2\int_0^t Y_s\,\dd s$, $t\in\RR_+$, \ we have
 \ $\int_0^t Y_s(\omega)\,\dd s=0$ \ for all \ $t\in[0,T]$ \ on the event
   \[
    A_T\cap \left\{ \omega\in\Omega :
                      \int_0^T Y_s(\omega) \, \dd s \int_0^T \frac{1}{Y_s(\omega)} \, \dd s = T^2
           \right\}.
  \]
However, \ $\left\{ \omega \in \Omega
            : \int_0^T Y_s(\omega) \, \dd s = 0 \right\} \bigcap A_T = \emptyset$,
 \ since \ $t \mapsto Y_t(\omega)$ \ is continuous and positive on \ $[0, T]$
 \ for all \ $\omega \in A_T$.
\ Consequently, since \ $\PP(A_T)=1$, \ we have
 \ $\PP\left( \int_0^T Y_s \, \dd s \int_0^T \frac{1}{Y_s} \, \dd s = T^2 \right)
    = 0$.
\proofend

\section{Consistency of MLE}
\label{section_CMLE}

First we consider the case of subcritical Heston models, i.e., when
 \ $b \in \RR_{++}$.

\begin{Thm}\label{Thm_MLE_cons_sub}
If \ $b \in \RR_{++}$, \ $\alpha, \beta \in \RR$,
 \ $\sigma_1, \sigma_2 \in \RR_{++}$, \ $\varrho \in (-1, 1)$, \ and
 \ $(Y_0, X_0) = (y_0, x_0) \in \RR_{++} \times \RR$, \ then the MLE of
 \ $(a, b, \alpha, \beta)$ \ is strongly consistent, i.e.,
 \ $\bigl(\ha_T, \hb_T, \halpha_T, \hbeta_T\bigr)
    \as (a, b, \alpha, \beta)$
 \ as \ $T \to \infty$, \ whenever
 \ $a \in \left( \frac{\sigma_1^2}{2}, \infty \right)$, \ and it is weakly
 consistent, i.e.,
 \ $\bigl(\ha_T, \hb_T, \halpha_T, \hbeta_T\bigr)
    \stoch (a, b, \alpha, \beta)$
 \ as \ $T \to \infty$, \ whenever \ $a = \frac{\sigma_1^2}{2}$.
\end{Thm}

\noindent{\bf Proof.}
In both cases we have to show coordinate-wise convergences.
Indeed, for the almost sure convergence, one can use that the intersection of
 four events with probability one is an event with probability one, and for
 the convergence in probability one can apply, e.g., van der Vaart
 \cite[Theorem 2.7, part (vi)]{Vaart}.

By Lemma \ref{LEMMA_MLE_exist}, there exists a unique MLE
 \ $\bigl(\ha_T, \hb_T, \halpha_T, \hbeta_T\bigr)$ \ of
 \ $(a, b, \alpha, \beta)$ \ for all \ $T \in \RR_{++}$, \ which has the form
 given in \eqref{MLE}.
By \eqref{MLE-}, we have
 \begin{align}\label{deco_a}
  \halpha_T - \alpha
  = \frac{\sigma_2
          \cdot
          \frac{\int_0^T \frac{\dd \tW_s}{\sqrt{Y_s}}}
               {\int_0^T \frac{\dd s}{Y_s}}
          - \frac{\sigma_2}{\frac{1}{T} \int_0^T \frac{\dd s}{Y_s}}
            \cdot\frac{\int_0^T \sqrt{Y_s} \, \dd \tW_s}
                      {\int_0^T Y_s \, \dd s}}
         {1 - \frac{1}{\frac{1}{T} \int_0^T Y_s \, \dd s \,
                       \cdot \frac{1}{T} \int_0^T \frac{\dd s}{Y_s}}}
 \end{align}
 provided that \ $\int_0^T Y_s \, \dd s \int_0^T \frac{\dd s}{Y_s} > T^2$
 \ (implying \ $\int_0^T Y_s \, \dd s \int_0^T \frac{\dd s}{Y_s} \in \RR_{++}$)
 \ which holds a.s.

First we consider the case of
 \ $a \in \left( \frac{\sigma_1^2}{2}, \infty \right)$.
\ The strong consistency of the MLE of \ $(a, b)$ \ has been proved by
 Overbeck \cite[Theorem 2, part (ii)]{Ove}.
By part (i) of Theorem \ref{Ergodicity}, \ $\EE(Y_\infty) = \frac{a}{b}$
 \ and \ $\EE\left(\frac{1}{Y_\infty}\right) = \frac{2b}{2a-\sigma_1^2}$,
 \ and hence, part (iii) of Theorem \ref{Ergodicity} implies
 \begin{align}\label{SLLN}
  \frac{1}{T} \int_0^T Y_s \, \dd s \as \EE(Y_\infty) \qquad \text{and} \qquad
  \frac{1}{T} \int_0^T \frac{\dd s}{Y_s}
  \as \EE\left(\frac{1}{Y_\infty}\right) \qquad
  \text{as \ $T \to \infty$.}
 \end{align}
Further, since \ $\EE(Y_\infty), \EE\left(\frac{1}{Y_\infty}\right) \in \RR_{++}$,
 \ \eqref{SLLN} yields
 \begin{align*}
  \int_0^T Y_s \, \dd s \as \infty \qquad \text{and} \qquad
  \int_0^T \frac{\dd s}{Y_s} \as \infty \qquad
  \text{as \ $T \to \infty$.}
 \end{align*}
Applying a strong law of large numbers for continuous local martingales
 (see, e.g., Theorem \ref{DDS_stoch_int}), we obtain
 \begin{align*}
  \halpha_T - \alpha
  \as \frac{\sigma_2 \cdot 0
            - \frac{\sigma_2}
                   {\frac{2b}{2a-\sigma_1^2}} \cdot 0}
           {1 - \frac{1}
                     {\frac{a}{b} \cdot \frac{2b}{2a-\sigma_1^2}}}
  = 0 \qquad \text{as \ $T \to \infty$,}
 \end{align*}
 where we also used that the denominator above is not zero due to
 \ $\sigma_1\in\RR_{++}$.

Next we consider the case of \ $a = \frac{\sigma_1^2}{2}$.
\ Weak consistency of the MLE of \ $(a, b)$ \ follows from part 1 of Theorem 7
 in Ben Alaya and Kebaier \cite{BenKeb2}.
We have again \ $\EE(Y_\infty) = \frac{a}{b} \in \RR_{++}$, \ implying
 \begin{align}\label{SLLN=}
  \frac{1}{T} \int_0^T Y_s \, \dd s \as \EE(Y_\infty) \qquad \text{and} \qquad
  \int_0^T Y_s \, \dd s \as \infty \qquad
  \text{as \ $T \to \infty$.}
 \end{align}
Due to Ben Alaya and Kebaier \cite[Proposition 4]{BenKeb1}, we have
 \begin{align}\label{AK}
  \frac{1}{T^2} \int_0^T \frac{\dd s}{Y_s} \distr \tau \qquad
  \text{as \ $T \to \infty$,}
 \end{align}
 where \ $\tau := \inf\{ t \in \RR_{++} : \cW_t = \frac{b}{\sigma_1}\}$ \ with
 a standard Wiener process \ $(\cW_t)_{t\in\RR_+}$.
\ Since \ $\PP(\tau \in \RR_{++}) = 1$, \ we conclude
 \[
   \frac{1}{\frac{1}{T}\int_0^T \frac{\dd s}{Y_s}}
    = \frac{1}{T} \, \frac{1}{\frac{1}{T^2} \int_0^T \frac{\dd s}{Y_s}}
    \distr 0 \cdot \frac{1}{\tau} = 0 \qquad
   \text{as \ $T \to \infty$,}
 \]
 and hence,
 \begin{equation}\label{stoch_1/Y}
   \frac{1}{\frac{1}{T}\int_0^T \frac{\dd s}{Y_s}} \stoch 0 \qquad
   \text{as \ $T \to \infty$,}
 \end{equation}
 implying also
 \[
   \frac{1}{\int_0^T \frac{\dd s}{Y_s}}
   = \frac{1}{T} \, \frac{1}{\frac{1}{T} \int_0^T \frac{\dd s}{Y_s}}
   \stoch 0 \qquad
   \text{as \ $T \to \infty$.}
 \]
Since the function
 \ $\RR_{++} \ni T \mapsto \frac{1}{\int_0^T \frac{\dd s}{Y_s}}$ \ is
 monotone decreasing, we obtain
 \[
   \frac{1}{\int_0^T \frac{\dd s}{Y_s}} \as 0 \qquad \text{and} \qquad
   \int_0^T \frac{\dd s}{Y_s} \as \infty \qquad
   \text{as \ $T \to \infty$.}
 \]
Using \eqref{deco_a} and a strong law of large numbers for continuous local
 martingales (see, e.g., Theorem \ref{DDS_stoch_int}), we obtain
 \begin{align*}
  \halpha_T - \alpha
  \stoch \frac{\sigma_2 \cdot 0 - 0 \cdot 0}
              {1 - \frac{b}{a} \cdot 0}
  = 0 \qquad \text{as \ $T \to \infty$.}
 \end{align*}
Here we have convergence only in probability because of \eqref{stoch_1/Y}.

By \eqref{MLE-}, we have
 \begin{align}\label{deco_b}
  \hbeta_T - \beta
  = \frac{\frac{\sigma_2}{\frac{1}{T} \int_0^T Y_s \, \dd s}
          \cdot
          \frac{\int_0^T \frac{\dd \tW_s}{\sqrt{Y_s}}}
               {\int_0^T \frac{\dd s}{Y_s}}
          - \sigma_2
            \cdot\frac{\int_0^T \sqrt{Y_s} \, \dd \tW_s}
                      {\int_0^T Y_s \, \dd s}}
         {1 - \frac{1}{\frac{1}{T} \int_0^T Y_s \, \dd s \,
                       \cdot \frac{1}{T} \int_0^T \frac{\dd s}{Y_s}}}
 \end{align}
 provided that \ $\int_0^T Y_s \, \dd s \int_0^T \frac{\dd s}{Y_s} > T^2$
 \ (implying \ $\int_0^T Y_s \, \dd s \int_0^T \frac{\dd s}{Y_s} \in \RR_{++}$)
 \ which holds a.s.

First we consider the case of
 \ $a \in \left( \frac{\sigma_1^2}{2}, \infty \right)$.
\ Applying again a strong law of large numbers for continuous local
 martingales (see, e.g., Theorem \ref{DDS_stoch_int}), we obtain
 \begin{align*}
  \hbeta_T - \beta
  \as \frac{\frac{\sigma_2}{\EE(Y_\infty)} \cdot 0
            - \sigma_2 \cdot 0}
           {1 - \frac{1}
                     {\EE(Y_\infty) \EE\left(\frac{1}{Y_\infty}\right)}}
  = 0 \qquad \text{as \ $T \to \infty$,}
 \end{align*}
 where we also used that the denominator above is not zero due to
 \ $\sigma_1\in\RR_{++}$.

Next we consider the case of \ $a = \frac{\sigma_1^2}{2}$.
\ Using \eqref{stoch_1/Y} and \eqref{deco_b}, we obtain
 \begin{align*}
  \hbeta_T - \beta
  \stoch \frac{\frac{\sigma_2}{\EE(Y_\infty)} \cdot 0 - \sigma_2 \cdot 0}
              {1 - \frac{1}{\EE(Y_\infty)} \cdot 0}
  = 0 \qquad \text{as \ $T \to \infty$.}\\[-16mm]
 \end{align*}
\proofend

In order to handle supercritical Heston models, i.e., when \ $b \in \RR_{--}$,
 \ we need the following integral version of the Toeplitz Lemma, due to Dietz
 and Kutoyants \cite{DieKut}.

\begin{Lem}\label{int_Toeplitz}
Let \ $\{\varphi_T : T \in \RR_+\}$ \ be a family of probability measures on
 \ $\RR_+$ \ such that \ $\varphi_T([0,T]) = 1$ \ for all \ $T \in \RR_+$,
 \ and \ $\lim_{T\to\infty} \varphi_T([0,K]) = 0$ \ for all \ $K \in \RR_{++}$.
\ Then for every bounded and measurable function \ $f : \RR_+ \to \RR$ \ for
 which the limit \ $f(\infty) := \lim_{t\to\infty} f(t)$ \ exists, we have
 \[
   \lim_{T\to\infty} \int_0^\infty f(t) \, \varphi_T(\dd t) = f(\infty) .
 \]
\end{Lem}

As a special case, we have the following integral version of the Kronecker
 Lemma, see K\"uchler and S{\o}rensen \cite[Lemma B.3.2]{KS}.

\begin{Lem}\label{int_Kronecker}
Let \ $a : \RR_+ \to \RR_+$ \ be a measurable function.
Put \ $b(T) := \int_0^T a(t) \, \dd t$, \ $T \in \RR_+$.
\ Suppose that \ $\lim_{T\to\infty} b(T) = \infty$.
\ Then for every bounded and measurable function \ $f : \RR_+ \to \RR$ \ for
 which the limit \ $f(\infty) := \lim_{t\to\infty} f(t)$ \ exists, we have
 \[
   \lim_{T\to\infty} \frac{1}{b(T)} \int_0^T a(t) f(t) \, \dd t = f(\infty) .
 \]
\end{Lem}

The next theorem states strong consistency of the MLE of \ $b$ \ in the
 supercritical case.
Overbeck \cite[Theorem 2, part (i)]{Ove} contains this result for CIR
 processes with a slightly incomplete proof.

\begin{Thm}\label{Thm_MLE_cons_super}
If \ $a \in \left[ \frac{\sigma_1^2}{2}, \infty \right)$, \ $b \in \RR_{--}$,
 \ $\alpha, \beta \in \RR$, \ $\sigma_1, \sigma_2 \in \RR_{++}$,
 \ $\varrho \in (-1, 1)$, \ and
 \ $(Y_0, X_0) = (y_0, x_0) \in \RR_{++} \times \RR$, \ then the MLE
 of \ $b$ \ is strongly consistent, i.e., \ $\hb_T \as b$ \ as
 \ $T \to \infty$.
\end{Thm}

\noindent{\bf Proof.}
By Lemma \ref{LEMMA_MLE_exist}, there exists a unique MLE \ $\hb_T$ \ of
 \ $b$ \ for all \ $T\in\RR_{++}$ \ which has the form given in \eqref{MLE}.
First we check that
 \begin{align*}
  \EE( Y_t \mid \cF^Y_s )
  = \EE( Y_t \mid Y_s )
  = \ee^{-b(t-s)} Y_s + a \int_s^t \ee^{-b(t-u)} \, \dd u
 \end{align*}
 for all \ $s, t \in \RR_+$ \ with \ $0 \leq s \leq t$, \ where \ $\cF^{Y}_s$
 \ denotes the \ $\sigma$-algebra \ $\sigma(\{Y_u, \, u\in[0,s]\})$.
\ The first equality follows from the Markov property of the process
 \ $(Y_t)_{t\in\RR_+}$.
\ The second equality is a consequence of the time-homogeneity of the Markov
 process \ $Y$ \ and
 \begin{align*}
  \EE\bigl( Y_t \mid (Y_0, X_0) = (y_0, x_0) \bigr)
  = \ee^{-bt} y_0 + a \int_0^t \ee^{-b(t-u)} \, \dd u , \qquad t \in \RR_+ ,
 \end{align*}
 valid for all \ $(y_0, x_0) \in \RR_+ \times \RR$, \ following from
 Proposition \ref{Pro_moments}.
Thus
 \[
   \EE( \ee^{bt} Y_t \mid \cF^Y_s )
   = \ee^{bs} Y_s + a \int_s^t \ee^{bu} \, \dd u
   \geq \ee^{bs} Y_s
 \]
 for all \ $s, t \in \RR_+$ \ with \ $0 \leq s \leq t$, \ consequently, the
 process \ $(\ee^{bt} Y_t)_{t\in\RR_+}$ \ is a non-negative submartingale with
 respect to the filtration \ $(\cF^Y_t)_{t\in\RR_+}$.
\ Moreover,
 \[
   \EE(\ee^{bt} Y_t)
   = y_0 + a \int_0^t \ee^{bu} \, \dd u
   \leq y_0 + a \int_0^\infty \ee^{bu} \, \dd u
   = y_0 - \frac{a}{b}
   < \infty , \qquad t \in \RR_+ ,
 \]
 hence, by the submartingale convergence theorem, there exists a non-negative
 random variable \ $V$ \ such that
 \begin{align}\label{lim_Y}
  \ee^{bt} Y_t \as V \qquad \text{as \ $t \to \infty$.}
 \end{align}
Note that the distribution of \ $V$ \ coincides with the distribution of
 \ $\tcY_{-1/b}$, \ where \ $(\tcY_t)_{t\in\RR_+}$ \ is a CIR process given by
 the SDE
 \begin{align*}
  \dd \tcY_t = a \dd t + \sigma_1 \sqrt{\tcY_t} \, \dd \cW_t ,
  \qquad t \in \RR_+ ,
 \end{align*}
 with initial value \ $\tcY_0 = y_0$, \ where \ $(\cW_t)_{t\in\RR_+}$ \ is a
 standard Wiener process, see Ben Alaya and Kebaier
 \cite[Proposition 3]{BenKeb1}.
Consequently, \ $\PP(V \in \RR_{++}) = 1$ \ due to
 \ $\PP(\tcY_t\in\RR_{++}, \, \forall\,t\in\RR_+)=1$.
\ If \ $\omega \in \Omega$ \ such that \ $\RR_+ \ni t \mapsto Y_t(\omega)$
 \ is continuous and \ $\ee^{bt} Y_t(\omega) \to V(\omega)$ \ as
 \ $t \to \infty$, \ then, by the integral Kronecker Lemma \ref{int_Kronecker}
 with \ $f(t) = \ee^{bt} Y_t(\omega)$ \ and \ $a(t) = \ee^{-bt}$,
 \ $t \in \RR_+$, \ we have
 \[
   \frac{1}{\int_0^t \ee^{-bu} \, \dd u}
   \int_0^t \ee^{-bu} (\ee^{bu} Y_u(\omega)) \, \dd u
   \to V(\omega)  \qquad \text{as \ $t\to\infty$.}
 \]
Here \ $\int_0^t \ee^{-bu} \, \dd u = \frac{\ee^{-bt} - 1}{-b}$, \ $t \in \RR_+$,
 \ thus we conclude
 \begin{align}\label{lim_intY}
  \ee^{bt} \int_0^t Y_u \, \dd u \as -\frac{V}{b}  \qquad
  \text{as \ $t \to \infty$.}
 \end{align}
Further,
 \begin{align}\label{lim_1/Y}
  \int_0^t \frac{\dd u}{Y_u} \as \int_0^\infty \frac{\dd u}{Y_u}  \qquad
  \text{as \ $t \to \infty$,}
 \end{align}
 where
 \ $\int_0^\infty \frac{\dd u}{Y_u} \distre \int_0^{-1/b} \tcY_u \, \dd u$,
 \ see Ben-Alaya and Kebaier \cite[Proposition 4]{BenKeb1}.
Consequently, \ $\PP(\int_0^\infty \frac{\dd u}{Y_u} \in \RR_{++}) = 1$ \ due to
 \ $\PP(\tcY_t\in\RR_{++}, \, \forall\,t\in\RR_+)=1$.

Since \ $\PP(\text{$Y_t \in \RR_{++}$ for all \ $t \in\RR_+$}) = 1$,
 \ one can apply It\^o's rule to the function \ $f(x) = \log x$,
 \ $x \in \RR_{++}$, \ for which \ $f'(x) = 1 / x$,
 \ $f''(x) = - 1 / x^2$, \ $x \in \RR_{++}$, \ and we obtain
 \begin{align}\label{logY}
   \log Y_t
   = \log y_0 + \int_0^t \frac{\dd Y_s}{Y_s}
     - \frac{\sigma_1^2}{2} \int_0^t \frac{\dd s}{Y_s} , \qquad t \in \RR_+ ,
 \end{align}
 for all \ $b \in \RR$, \ see von Weizs\"acker and Winkler
 \cite[Theorem 8.1.1]{WeiWin}.

Using \eqref{MLE} and \eqref{logY}, we have
 \begin{align*}
  \hb_T
  = \frac{\frac{T \int_0^T \frac{\dd Y_s}{Y_s}}
               {\int_0^T Y_s \, \dd s \int_0^T \frac{\dd s}{Y_s}}
          - \frac{Y_T - y_0}
                 {\int_0^T Y_s \, \dd s}}
         {1 - \frac{1}{\frac{1}{T^2} \int_0^T Y_s \, \dd s
                       \int_0^T \frac{\dd s}{Y_s}}}
  = \frac{\frac{T (\log Y_T - \log y_0)}
               {\int_0^T Y_s \, \dd s \int_0^T \frac{\dd s}{Y_s}}
          - \frac{Y_T - y_0 - \frac{\sigma_1^2}{2} T}
                 {\int_0^T Y_s \, \dd s}}
         {1 - \frac{1}{\frac{1}{T^2} \int_0^T Y_s \, \dd s
                       \int_0^T \frac{\dd s}{Y_s}}}
 \end{align*}
 provided that \ $\int_0^T Y_s \, \dd s \int_0^T \frac{\dd s}{Y_s} > T^2$
 \ (implying \ $\int_0^T Y_s \, \dd s \int_0^T \frac{\dd s}{Y_s} \in \RR_{++}$)
 \ which holds a.s.
\ Applying \eqref{lim_Y}, \eqref{lim_intY} and \eqref{lim_1/Y}, we conclude
 \begin{align*}
  \hb_T
  = \frac{\frac{T\ee^{bT}\log(\ee^{bT}Y_T) - bT^2\ee^{bT} - T\ee^{bT}\log y_0}
               {\left(\ee^{bT} \int_0^T Y_s \, \dd s\right)
                \int_0^T \frac{\dd s}{Y_s}}
         - \frac{\ee^{bT} Y_T - \ee^{bT} y_0 - \frac{\sigma_1^2}{2}T\ee^{bT}}
                {\ee^{bT} \int_0^T Y_s \, \dd s}}
         {1 - T^2\ee^{bT}\frac{1}{\left(\ee^{bT} \int_0^T Y_s \, \dd s\right)
                       \int_0^T \frac{\dd s}{Y_s}}}
  \as \frac{ \frac{0 \cdot \log V - 0 }
                         {-\frac{V}{b} \int_0^\infty \frac{\dd s}{Y_s}}
            - \frac{V - 0}{-\frac{V}{b}}}
           {1 - 0 \cdot \frac{1}{-\frac{V}{b}
                                 \int_0^\infty \frac{\dd s}{Y_s}}}
  = b
 \end{align*}
 as \ $T \to \infty$.
\proofend

\begin{Rem}\label{Rem_subcrit_comp}
For subcritical (i.e., \ $b\in\RR_{++}$) \ CIR models with
 \ $a \in \bigl(\frac{\sigma_1^2}{2}, \infty\bigr)$, \ Overbeck
 \cite[Theorem 2, part (ii)]{Ove} proved strong consistency of the MLE of
 \ $(a,b)$.
For subcritical (i.e., \ $b\in\RR_{++}$) \ CIR models with
 \ $a = \frac{\sigma_1^2}{2}$, \ weak consistency of the MLE of \ $(a, b)$
 \ follows from part 1 of Theorem 7 in Ben Alaya and Kebaier \cite{BenKeb2}.
\proofend
\end{Rem}

\begin{Rem}\label{Rem_comparison5}
For critical (i.e., \ $b=0$) \ CIR models with
 \ $a \in \bigl[\frac{\sigma_1^2}{2}, \infty\bigr)$,
 \ weak consistency of the MLE of \ $(a,b)$ \ follows from Theorem 2 (iii) in
 Overbeck \cite{Ove} or Theorem 6 in Ben Alaya and Kebaier \cite{BenKeb2}.
For critical Heston models with \ $a\in(\frac{\sigma_1^2}{2},\infty)$,
 \ weak consistency of the MLE of \ $(a, b, \alpha, \beta)$ \ is a consequence
 of Theorem \ref{Thm_MLE}.
\proofend
\end{Rem}

\begin{Rem}\label{Rem_comparison1}
For supercritical (i.e., \ $b\in\RR_{--}$) \ CIR models with
 \ $a \in \bigl[\frac{\sigma_1^2}{2}, \infty\bigr)$,
 \ Overbeck \cite[Theorem 2, parts (i) and (v)]{Ove} proved that the MLE of
 \ $b$ \ is strongly consistent, however, there is no strongly
 consistent estimator of \ $a$.
\ See also Ben Alaya and Kebaier \cite[Theorem 7, part 2]{BenKeb2}.
For supercritical Heston models with
 \ $a \in \bigl[\frac{\sigma_1^2}{2}, \infty\bigr)$, \ it will turn out that
 the MLE of \ $a$ \ and
 \ $\alpha$ \ is not even weakly consistent, but the MLE of \ $\beta$
 \ is weakly consistent, see Theorem \ref{Thm_MLE_super}.
\proofend
\end{Rem}

\section{Asymptotic behaviour of MLE: subcritical case}
\label{section_AMLE_subcritical}

We consider subcritical Heston models, i.e., when \ $b \in \RR_{++}$.

\begin{Thm}\label{Thm_MLE_sub}
If  \ $a \in \left( \frac{\sigma_1^2}{2}, \infty \right)$, \ $b \in \RR_{++}$,
 \ $\alpha, \beta \in \RR$, \ $\sigma_1, \sigma_2 \in \RR_{++}$,
 \ $\varrho \in (-1, 1)$, \ and
 \ $(Y_0, X_0) = (y_0, x_0) \in \RR_{++} \times \RR$, \ then the MLE of
 \ $(a, b, \alpha, \beta)$ \ is asymptotically normal, i.e.,
 \begin{align}\label{MLE_sub}
  \sqrt{T}
  \begin{bmatrix}
   \ha_T - a \\
   \hb_T - b \\
   \halpha_T - \alpha \\
   \hbeta_T - \beta
  \end{bmatrix}
  \distr \cN_4\left(\bzero, \bS
     \otimes
     \begin{bmatrix}
      \frac{2b}{2a - \sigma_1^2} & - 1 \\
      - 1 & \frac{a}{b}
     \end{bmatrix}^{-1} \right) \qquad
  \text{as \ $T \to \infty$,}
 \end{align}
 where \ $\bS$ \ is defined in \eqref{bSigma}.

With a random scaling, we have
 \begin{align}\label{MLE_subr}
  \frac{1}{\bigl(\int_0^T \frac{\dd s}{Y_s} \bigr)^{1/2}}
  \left( \bI_2
         \otimes
         \begin{bmatrix}
          \int_0^T \frac{\dd s}{Y_s} & - T \\
          0 & \bigl(\int_0^T Y_s \, \dd s
                      \int_0^T \frac{\dd s}{Y_s} - T^2\bigr)^{1/2}
         \end{bmatrix} \right)
  \begin{bmatrix}
   \ha_T - a \\
   \hb_T - b \\
   \halpha_T - \alpha \\
   \hbeta_T - \beta
  \end{bmatrix}
  \distr \cN_4\left(\bzero, \bS \otimes \bI_2 \right)
 \end{align}
 as \ $T \to \infty$.
\end{Thm}

\noindent{\bf Proof.}
By Lemma \ref{LEMMA_MLE_exist}, there exists a unique MLE
 \ $\bigl(\ha_T, \hb_T, \halpha_T, \hbeta_T\bigr)$
 \ of \ $(a, b, \alpha, \beta)$ \ for all \ $T\in\RR_{++}$, \ which has
 the form given in \eqref{MLE}.
By \eqref{MLE-}, we have
 \begin{align*}
  \sqrt{T} (\ha_T - a)
  &= \frac{\frac{1}{T} \int_0^T Y_s \, \dd s \,
           \cdot
           \frac{\sigma_1}{\sqrt{T}} \int_0^T \frac{\dd W_s}{\sqrt{Y_s}}
           - \frac{\sigma_1}{\sqrt{T}} \int_0^T \sqrt{Y_s} \, \dd W_s}
          {\frac{1}{T} \int_0^T Y_s \, \dd s \,
           \cdot \frac{1}{T} \int_0^T \frac{\dd s}{Y_s}
           - 1} , \\
  \sqrt{T} (\hb_T - b)
  &= \frac{\frac{\sigma_1}{\sqrt{T}} \int_0^T \frac{\dd W_s}{\sqrt{Y_s}}
           - \frac{1}{T} \int_0^T \frac{\dd s}{Y_s} \,
             \cdot
             \frac{\sigma_1}{\sqrt{T}} \int_0^T \sqrt{Y_s} \, \dd W_s}
          {\frac{1}{T} \int_0^T Y_s \, \dd s \,
           \cdot \frac{1}{T} \int_0^T \frac{\dd s}{Y_s}
           - 1} , \\
  \sqrt{T} (\halpha_T - \alpha)
  &= \frac{\frac{1}{T} \int_0^T Y_s \, \dd s \,
           \cdot
           \frac{\sigma_2}{\sqrt{T}} \int_0^T \frac{\dd \tW_s}{\sqrt{Y_s}}
           - \frac{\sigma_2}{\sqrt{T}} \int_0^T \sqrt{Y_s} \, \dd \tW_s}
          {\frac{1}{T} \int_0^T Y_s \, \dd s \,
           \cdot \frac{1}{T} \int_0^T \frac{\dd s}{Y_s}
           - 1} , \\
  \sqrt{T} (\hbeta_T - \beta)
  &= \frac{\frac{\sigma_2}{\sqrt{T}} \int_0^T \frac{\dd \tW_s}{\sqrt{Y_s}}
           - \frac{1}{T} \int_0^T \frac{\dd s}{Y_s} \,
             \cdot
             \frac{\sigma_2}{\sqrt{T}} \int_0^T \sqrt{Y_s} \, \dd \tW_s}
          {\frac{1}{T} \int_0^T Y_s \, \dd s \,
           \cdot \frac{1}{T} \int_0^T \frac{\dd s}{Y_s} - 1}
 \end{align*}
 provided that \ $\int_0^T Y_s \, \dd s \int_0^T \frac{1}{Y_s} \, \dd s > T^2$
 \ which holds a.s.
\ Consequently,
 \begin{align*}
   \sqrt{T}
   \begin{bmatrix}
    \ha_T - a \\
    \hb_T - b \\
    \halpha_T - \alpha \\
    \hbeta_T - \beta
   \end{bmatrix}
   &= \frac{1}{\frac{1}{T} \int_0^T Y_s \, \dd s \,
               \cdot \frac{1}{T} \int_0^T \frac{\dd s}{Y_s} - 1}
      \left( \bI_2
             \otimes
             \begin{bmatrix}
              \frac{1}{T} \int_0^T Y_s \, \dd s & 1 \\
              1 & \frac{1}{T} \int_0^T \frac{\dd s}{Y_s}
             \end{bmatrix} \right)
      \frac{1}{\sqrt{T}} \bM_T \\
   &= \left( \bI_2
             \otimes
             \begin{bmatrix}
              \frac{1}{T} \int_0^T \frac{\dd s}{Y_s} & - 1 \\
              - 1 & \frac{1}{T} \int_0^T Y_s \, \dd s
             \end{bmatrix}^{-1} \right)
      \frac{1}{\sqrt{T}} \bM_T ,
 \end{align*}
 where
 \[
   \bM_t := \begin{bmatrix}
             \sigma_1 \int_0^t \frac{\dd W_s}{\sqrt{Y_s}} \\[1mm]
             - \sigma_1 \int_0^t \sqrt{Y_s} \, \dd W_s \\[1mm]
             \sigma_2 \int_0^t \frac{\dd \tW_s}{\sqrt{Y_s}} \\[1mm]
             - \sigma_2 \int_0^t \sqrt{Y_s} \, \dd \tW_s
            \end{bmatrix} , \qquad
   t \in \RR_+ .
 \]
Next, we show that
 \begin{align}\label{MLE_CLT}
  \frac{1}{\sqrt{T}} \bM_T \distr \bfeta \bZ \qquad
  \text{as \ $T \to \infty$,}
 \end{align}
 where \ $\bZ$ \ is a 4-dimensional standard normally distributed random
 vector and \ $\bfeta \in \RR^{4 \times 4}$ \ such that
 \[
   \bfeta \bfeta^\top
   = \bS
     \otimes
     \begin{bmatrix}
      \EE\left(\frac{1}{Y_\infty}\right) & - 1 \\
      - 1 & \EE(Y_\infty)
     \end{bmatrix} .
 \]
Here the two symmetric matrices on the right hand side are positive definite,
 since
 \[
    \sigma_1^2 \sigma_2^2  (1-\varrho^2) \in \RR_{++} \qquad \text{and} \qquad
    \EE\left(\frac{1}{Y_\infty}\right) \EE(Y_\infty) - 1
     = \frac{\sigma_1^2}{2a-\sigma_1^2} \in \RR_{++},
 \]
 so \ $\bfeta$ \ can be chosen, for instance, as the uniquely defined
 symmetric positive definite square root of the Kronecker product of the two
 matrices in question.
The process \ $(\bM_t)_{t\in\RR_+}$ \ is a 4-dimensional continuous local
 martingale with quadratic variation process
 \[
   \langle \bM \rangle_t
   = \bS
     \otimes
     \begin{bmatrix}
      \int_0^t \frac{1}{Y_s} \, \dd s & - t \\
      - t & \int_0^t Y_s \, \dd s
     \end{bmatrix} , \qquad t \in \RR_+ .
 \]
By Theorem \ref{Ergodicity}, we have
 \[
   \bQ(t) \langle \bM \rangle_t \, \bQ(t)^\top
   \as \bS
       \otimes
       \begin{bmatrix}
        \EE\left(\frac{1}{Y_\infty}\right) & - 1 \\
        - 1 & \EE(Y_\infty)
       \end{bmatrix} \qquad
   \text{as \ $t \to \infty$}
 \]
 with \ $\bQ(t) := t^{-1/2} \bI_4$, \ $t \in \RR_{++}$.
\ Hence, Theorem \ref{THM_Zanten} yields \eqref{MLE_CLT}.
Then Slutsky's lemma yields
 \begin{align*}
  \sqrt{T}
  \begin{bmatrix}
   \ha_T - a \\
   \hb_T - b \\
   \halpha_T - \alpha \\
   \hbeta_T - \beta
  \end{bmatrix}
  \distr
  \left( \bI_2
         \otimes
         \begin{bmatrix}
          \EE\left(\frac{1}{Y_\infty}\right) & -1 \\
          - 1 & \EE(Y_\infty)
         \end{bmatrix}^{-1} \right) \,
  \bfeta \bZ
  \distre
  \cN_4(\bzero, \bSigma_1)
  \qquad \text{as \ $T \to \infty$,}
 \end{align*}
 where (applying the identities
 \ $(\bA\otimes\bB)^\top = \bA^\top \otimes \bB^\top$ \ and
  \ $(\bA\otimes\bB)(\bC\otimes\bD) = (\bA\bC)\otimes(\bB\bD)$)
 \begin{align*}
  \bSigma_1
  &:= \left( \bI_2
             \otimes
             \begin{bmatrix}
              \EE\left(\frac{1}{Y_\infty}\right) & -1 \\
              - 1 & \EE(Y_\infty)
             \end{bmatrix}^{-1} \right)
      \left( \bS
             \otimes
             \begin{bmatrix}
              \EE\left(\frac{1}{Y_\infty}\right) & - 1 \\
              - 1 & \EE(Y_\infty)
             \end{bmatrix} \right)
      \left( \bI_2
             \otimes
             \begin{bmatrix}
              \EE\left(\frac{1}{Y_\infty}\right) & -1 \\
              - 1 & \EE(Y_\infty)
             \end{bmatrix}^{-1} \right)^\top \\
  &= ( \bI_2 \bS \bI_2 )
       \otimes
        \left( \begin{bmatrix}
              \EE\left(\frac{1}{Y_\infty}\right) & -1 \\
              - 1 & \EE(Y_\infty)
             \end{bmatrix}^{-1}
             \begin{bmatrix}
              \EE\left(\frac{1}{Y_\infty}\right) & - 1 \\
              - 1 & \EE(Y_\infty)
             \end{bmatrix}
            \left( \begin{bmatrix}
              \EE\left(\frac{1}{Y_\infty}\right) & -1 \\
              - 1 & \EE(Y_\infty)
             \end{bmatrix}^{-1}\right)^\top \right)
 \end{align*}
 \begin{align*}
  &= \bS
     \otimes
     \begin{bmatrix}
      \EE\left(\frac{1}{Y_\infty}\right) & - 1 \\
      - 1 & \EE(Y_\infty)
     \end{bmatrix}^{-1} ,
 \end{align*}
 which yields \eqref{MLE_sub} recalling \ $\EE(Y_\infty) = \frac{a}{b}$
 \ and \ $\EE\left(\frac{1}{Y_\infty}\right) = \frac{2b}{2a-\sigma_1^2}$.

Slutsky's lemma and \eqref{MLE_sub} yield
 \begin{align*}
  &\frac{1}{\bigl(\int_0^T \frac{\dd s}{Y_s} \bigr)^{1/2}}
   \left( \bI_2
          \otimes
          \begin{bmatrix}
           \int_0^T \frac{\dd s}{Y_s} & - T \\
           0 & \bigl(\int_0^T Y_s \, \dd s
                     \int_0^T \frac{\dd s}{Y_s} - T^2\bigr)^{1/2}
          \end{bmatrix} \right)
   \begin{bmatrix}
    \ha_T - a \\
    \hb_T - b \\
    \halpha_T - \alpha \\
    \hbeta_T - \beta
   \end{bmatrix} \\
  &= \frac{1}{\bigl(\frac{1}{T} \int_0^T \frac{\dd s}{Y_s} \bigr)^{1/2}}
     \left( \bI_2
            \otimes
            \begin{bmatrix}
             \frac{1}{T} \int_0^T \frac{\dd s}{Y_s} & - 1 \\
             0 & \bigl(\frac{1}{T} \int_0^T Y_s \, \dd s \cdot
                       \frac{1}{T} \int_0^T \frac{\dd s}{Y_s} - 1\bigr)^{1/2}
         \end{bmatrix} \right)
  \sqrt{T}
  \begin{bmatrix}
   \ha_T - a \\
   \hb_T - b \\
   \halpha_T - \alpha \\
   \hbeta_T - \beta
  \end{bmatrix}
 \end{align*}
 \begin{align*}
  &\distr
   \frac{1}{\bigl(\EE\bigl(\frac{1}{Y_\infty}\bigr)\bigr)^{1/2}}
   \left( \bI_2
          \otimes
          \begin{bmatrix}
           \EE\bigl(\frac{1}{Y_\infty}\bigr) & - 1 \\
           0 & \bigl(\EE(Y_\infty)
                     \EE\bigl(\frac{1}{Y_\infty}\bigr)-1\bigr)^{1/2}
         \end{bmatrix} \right)
   \cN_4\left(\bzero, \bS
                      \otimes
                      \begin{bmatrix}
                       \EE\bigl(\frac{1}{Y_\infty}\bigr) & - 1 \\
                       - 1 & \EE(Y_\infty)
                      \end{bmatrix}^{-1} \right) \\
  &\distre \cN_4(\bzero, \bSigma_2) \qquad \text{as \ $T \to \infty$,}
 \end{align*}
 where (applying the identities
 \ $(\bA\otimes\bB)^\top = \bA^\top \otimes \bB^\top$ \ and
 \ $(\bA\otimes\bB)(\bC\otimes\bD) = (\bA\bC)\otimes(\bB\bD)$)
 \begin{align*}
  \bSigma_2 &:= \frac{1}{\EE\bigl(\frac{1}{Y_\infty}\bigr)}
                \left( \bI_2
                     \otimes
                     \begin{bmatrix}
                      \EE\bigl(\frac{1}{Y_\infty}\bigr) & - 1 \\
                      0 & \bigl(\EE(Y_\infty)
                                \EE\bigl(\frac{1}{Y_\infty}\bigr)-1\bigr)^{1/2}
                     \end{bmatrix} \right)
              \left( \bS
                     \otimes
                     \begin{bmatrix}
                      \EE\bigl(\frac{1}{Y_\infty}\bigr) & - 1 \\
                      - 1 & \EE(Y_\infty)
                     \end{bmatrix}^{-1} \right) \\
          &\quad\:
              \times
              \left( \bI_2
                     \otimes
                     \begin{bmatrix}
                      \EE\bigl(\frac{1}{Y_\infty}\bigr) & - 1 \\
                      0 & \bigl(\EE(Y_\infty)
                                \EE\bigl(\frac{1}{Y_\infty}\bigr)-1\bigr)^{1/2}
                     \end{bmatrix} \right)^\top\\
    &= \frac{1}{\EE\bigl(\frac{1}{Y_\infty}\bigr)} (\bI_2\bS\bI_2) \\
    &\quad \otimes \!
              \left( \begin{bmatrix}
                      \EE\bigl(\frac{1}{Y_\infty}\bigr) & - 1 \\
                      0 & \bigl(\EE(Y_\infty) \EE\bigl(\frac{1}{Y_\infty}\bigr)-1\bigr)^{1/2}
                     \end{bmatrix}
                     \begin{bmatrix}
                      \EE\bigl(\frac{1}{Y_\infty}\bigr) & - 1 \\
                      - 1 & \EE(Y_\infty)
                     \end{bmatrix}^{-1}
                     \begin{bmatrix}
                      \EE\bigl(\frac{1}{Y_\infty}\bigr) & 0 \\
                      -1 & \bigl(\EE(Y_\infty) \EE\bigl(\frac{1}{Y_\infty}\bigr)-1\bigr)^{1/2}
                     \end{bmatrix}
                     \right) \\
    &= \bS \otimes \bI_2 ,
 \end{align*}
 since
 \[
   \begin{bmatrix}
    \EE\bigl(\frac{1}{Y_\infty}\bigr) & - 1 \\
    - 1 & \EE(Y_\infty)
   \end{bmatrix}
   = \frac{1}{\EE\bigl(\frac{1}{Y_\infty}\bigr)}
     \begin{bmatrix}
      \EE\bigl(\frac{1}{Y_\infty}\bigr) & 0 \\
      -1 & \bigl(\EE(Y_\infty) \EE\bigl(\frac{1}{Y_\infty}\bigr)-1\bigr)^{1/2}
     \end{bmatrix}
     \begin{bmatrix}
      \EE\bigl(\frac{1}{Y_\infty}\bigr) & - 1 \\
      0 & \bigl(\EE(Y_\infty) \EE\bigl(\frac{1}{Y_\infty}\bigr)-1\bigr)^{1/2}
     \end{bmatrix} .
 \]
Thus we obtain \eqref{MLE_subr}.
\proofend

\begin{Rem}\label{Rem_comparison2}
For subcritical (i.e., \ $b\in\RR_{++}$) \ CIR models, for the MLE of
 \ $(a, b)$, \ Ben Alaya and Kebaier \cite[Theorems 5 and 7]{BenKeb2} proved
 asymptotic normality whenever
 \ $a \in \bigl(\frac{\sigma_1^2}{2}, \infty\bigr)$, \ and derived a limit
 theorem with a non-normal limit distribution whenever
 \ $a = \frac{\sigma_1^2}{2}$.
\ For subcritical (i.e., \ $b\in\RR_{++}$) \ CIR models, for the MLE of
 \ $(a, b)$, \ with random scaling, Overbeck \cite[Theorem 3 (iii)]{Ove}
 showed asymptotic normality.
\proofend
\end{Rem}

\section{Asymptotic behaviour of MLE: critical case}
\label{section_AMLE_critical}

We consider critical Heston models, i.e., when \ $b = 0$.
\ First we present an auxiliary lemma.

\begin{Lem}\label{Lem_intf_cont}
The mapping
 \ $C(\RR_+, \RR) \ni f \mapsto \bigl(\int_0^t f(u) \, \dd u\bigr)_{t\in\RR_+}
    \in C(\RR_+, \RR)$
 \ is continuous, hence measurable, where \ $C(\RR_+, \RR)$ \ denotes the set
 of real-valued continuous functions defined on \ $\RR_+$.
\end{Lem}

\noindent{\bf Proof.}
The space \ $C(\RR_+, \RR)$ \ is topologized by the locally uniform metric
 \[
   \delta_\slu(f, g)
   := \sum_{N=1}^\infty
       2^{-N} \min\biggl\{ 1, \sup_{t\in[0,N]} |f(t) - g(t)| \biggr\} ,
   \qquad f, g \in C(\RR_+, \RR),
 \]
 see, e.g., Jacod and Shiryaev \cite[Chapter VI, Section 1]{JSh}.
Let \ $f \in C(\RR_+, \RR)$ \ and \ $f_n \in C(\RR_+, \RR)$, \ $n \in \NN$,
 \ such that \ $\delta_\slu(f, f_n) \to 0$ \ as \ $n \to \infty$.
\ Put \ $F(t) := \int_0^t f(s) \, \dd s$, \ $t \in \RR_+$, \ and
 \ $F_n(t) := \int_0^t f_n(s) \, \dd s$, \ $t \in \RR_+$, \ $n \in \NN$.
\ Then \ $\sup_{t\in[0,N]} |F(t) - F_n(t)| \leq N \sup_{t\in[0,N]} |f(t) - f_n(t)|$
 \ for all \ $N \in \NN$, \ hence for each \ $K \in \NN$, \ we have
 \begin{align*}
  \delta_\slu(F, F_n)
  &= \sum_{N=K+1}^\infty
      2^{-N} \min\biggl\{ 1, \sup_{t\in[0,N]} |F(t) - F_n(t)| \biggr\}
     + \sum_{N=1}^K
        2^{-N} \min\biggl\{ 1, \sup_{t\in[0,N]} |F(t) - F_n(t)| \biggr\} \\
  &\leq \sum_{N=K+1}^\infty 2^{-N}
        + \sum_{N=1}^K
           2^{-N} \min\biggl\{ 1, N \sup_{t\in[0,N]} |f(t) - f_n(t)| \biggr\} \\
  &\leq 2^{-K}
        + \sum_{N=1}^K
           N 2^{-N} \min\biggl\{ 1, \sup_{t\in[0,N]} |f(t) - f_n(t)| \biggr\} \\
  &\leq 2^{-K} + \delta_\slu(f, f_n) \sum_{N=1}^K N
   = 2^{-K} + \frac{K(K+1)}{2} \delta_\slu(f, f_n)
   \to 2^{-K} \qquad \text{as \ $n \to \infty$.}
 \end{align*}
Consequently,
 \[
   \limsup_{n\to\infty} \delta_\slu(F, F_n) \leq 2^{-K} \qquad
   \text{for all \ $K \in \NN$,}
 \]
 thus we obtain the statement.

We present another short proof.
Applying Problem 3.11.26 in Ethier and Kurtz \cite{EthKur} and Proposition
 VI.1.17  in Jacod and Shiryaev \cite{JSh}, the mapping
 \ $C(\RR_+,\RR) \ni f \mapsto (\int_0^t f(u) \, \dd u)_{t\in\RR_+}
    \in C(\RR_+,\RR)$
 \ is continuous, hence measurable.
\proofend

The next result can be considered as a generalization of part 2 of Theorem 6
 in Ben Alaya and Kebaier \cite{BenKeb2} for critical Heston models.

\begin{Thm}\label{Thm_MLE}
If \ $a \in \bigl( \frac{\sigma_1^2}{2}, \infty \bigr)$, \ $b = 0$,
 \ $\alpha, \beta \in \RR$, \ $\sigma_1, \sigma_2 \in \RR_{++}$,
 \ $\varrho \in (-1, 1)$ \ and
 \ $(Y_0, X_0) = (y_0, x_0) \in \RR_{++} \times \RR$, \ then
 \begin{align}\label{aabb}
  \begin{bmatrix}
   \sqrt{\log T} (\ha_T - a) \\
   \sqrt{\log T} (\halpha_T - \alpha) \\
   T \hb_T \\
   T (\hbeta_T - \beta)
  \end{bmatrix}
  \distr
  \begin{bmatrix}
   \left(a - \frac{\sigma_1^2}{2}\right)^{1/2} \bS^{1/2} \bZ_2 \\
   \frac{a - \cY_1}{\int_0^1 \cY_s \, \dd s} \\[1mm]
   \frac{\alpha - \cX_1}{\int_0^1 \cY_s \, \dd s}
  \end{bmatrix}
  \qquad \text{as \ $T \to \infty$,}
 \end{align}
 where \ $(\cY_t, \cX_t)_{t\in\RR_+}$ \ is the unique strong solution of the SDE
 \begin{align}\label{help_limit_YX}
  \begin{cases}
   \dd \cY_t = a \, \dd t + \sigma_1 \sqrt{\cY_t} \, \dd \cW_t , \\
   \dd \cX_t = \alpha \, \dd t
               + \sigma_2  \sqrt{\cY_t}
                 \bigl(\varrho \, \dd \cW_t
                       + \sqrt{1 - \varrho^2} \, \dd \cB_t\bigr) ,
  \end{cases} \qquad t \in \RR_+ ,
 \end{align}
 with initial value \ $(\cY_0, \cX_0) = (0, 0)$, \ where
 \ $(\cW_t, \cB_t)_{t\in\RR_+}$ \ is a $2$-dimensional standard Wiener process,
 \ $\bZ_2$ \ is a $2$-dimensional standard normally distributed random vector
 independent of \ $\bigl(\cY_1, \int_0^1 \cY_t \, \dd t, \cX_1\bigr)$,
 \ $\bS$ \ is defined in \eqref{bSigma}, and \ $\bS^{1/2}$ \ denotes its
 uniquely determined symmetric, positive definite square root.
\end{Thm}

\noindent{\bf Proof.}
By Lemma \ref{LEMMA_MLE_exist}, there exists a unique MLE
 \ $\bigl(\ha_T, \hb_T, \halpha_T, \hbeta_T\bigr)$
 \ of \ $(a, b, \alpha, \beta)$ \ for all \ $T\in\RR_{++}$, \ which has the
 form given in \eqref{MLE}.
By \eqref{MLE-}, we have
 \begin{align*}
  \sqrt{\log T} (\ha_T - a)
  = \frac{\frac{1}
               {\bigl(\frac{1}{\log T} \int_0^T \frac{\dd s}{Y_s}\bigr)^{1/2}}
          \frac{\sigma_1 \int_0^T \frac{\dd W_s}{\sqrt{Y_s}}}
               {\bigl(\int_0^T \frac{\dd s}{Y_s}\bigr)^{1/2}}
          - \frac{\frac{1}{\sqrt{\log T}}}
                 {\frac{1}{\log T} \int_0^T \frac{\dd s}{Y_s}}
            \frac{T \sigma_1 \int_0^T \sqrt{Y_s} \, \dd W_s}
                 {\int_0^T Y_s \, \dd s}}
         {1 - \frac{1}{\frac{1}{T^2} \int_0^T Y_s \, \dd s}
              \frac{1}{\int_0^T \frac{\dd s}{Y_s}}} ,
 \end{align*}
 \begin{align*}
  \sqrt{\log T} (\halpha_T - \alpha)
  = \frac{\frac{1}
               {\bigl(\frac{1}{\log T} \int_0^T \frac{\dd s}{Y_s}\bigr)^{1/2}}
          \frac{\sigma_2 \int_0^T \frac{\dd \tW_s}{\sqrt{Y_s}}}
               {\bigl(\int_0^T \frac{\dd s}{Y_s}\bigr)^{1/2}}
          - \frac{\frac{1}{\sqrt{\log T}}}
                 {\frac{1}{\log T} \int_0^T \frac{\dd s}{Y_s}}
            \frac{T \sigma_2 \int_0^T \sqrt{Y_s} \, \dd \tW_s}
                 {\int_0^T Y_s \, \dd s}}
         {1 - \frac{1}{\frac{1}{T^2} \int_0^T Y_s \, \dd s}
              \frac{1}{\int_0^T \frac{\dd s}{Y_s}}} ,
 \end{align*}
 \begin{align*}
  T \hb_T
  = \frac{\frac{1}{\frac{1}{T^2} \int_0^T Y_s \, \dd s}
          \frac{1}
               {\bigl(\int_0^T \frac{\dd s}{Y_s}\bigr)^{1/2}}
          \frac{\sigma_1 \int_0^T \frac{\dd W_s}{\sqrt{Y_s}}}
               {\bigl(\int_0^T \frac{\dd s}{Y_s}\bigr)^{1/2}}
          - \frac{T \sigma_1 \int_0^T \sqrt{Y_s} \, \dd W_s}
                 {\int_0^T Y_s \, \dd s}}
         {1 - \frac{1}{\frac{1}{T^2} \int_0^T Y_s \, \dd s}
              \frac{1}{\int_0^T \frac{\dd s}{Y_s}}} ,
 \end{align*}
 and
 \begin{align*}
  T (\hbeta_T - \beta)
  = \frac{\frac{1}{\frac{1}{T^2} \int_0^T Y_s \, \dd s}
          \frac{1}
               {\bigl(\int_0^T \frac{\dd s}{Y_s}\bigr)^{1/2}}
          \frac{\sigma_2 \int_0^T \frac{\dd \tW_s}{\sqrt{Y_s}}}
               {\bigl(\int_0^T \frac{\dd s}{Y_s}\bigr)^{1/2}}
          - \frac{T \sigma_2 \int_0^T \sqrt{Y_s} \, \dd \tW_s}
                 {\int_0^T Y_s \, \dd s}}
         {1 - \frac{1}{\frac{1}{T^2} \int_0^T Y_s \, \dd s}
              \frac{1}{\int_0^T \frac{\dd s}{Y_s}}} ,
 \end{align*}
 provided that \ $\int_0^T Y_s \, \dd s \int_0^T \frac{1}{Y_s} \, \dd s > T^2$
 \ which holds a.s.
It is known that
 \begin{align}\label{1/Y}
  \frac{1}{\log T} \int_0^T \frac{\dd s}{Y_s}
  \stoch \left(a - \frac{\sigma_1^2}{2}\right)^{-1} \qquad
  \text{as \ $T \to \infty$,}
 \end{align}
 see, e.g., Overbeck \cite[Lemma 5]{Ove} or Ben Alaya and Kebaier
 \cite[Proposition 2]{BenKeb1}.
Consequently,
 \begin{align}\label{lim_int_1/Y}
 \frac{1}{\int_0^T \frac{\dd s}{Y_s}} \as 0
 \qquad \text{and} \qquad
  \int_0^T \frac{\dd s}{Y_s} \as \infty  \qquad
  \text{as \ $T \to \infty$,}
 \end{align}
 where we used that \ $\bigl(\int_0^t \frac{\dd s}{Y_s}\bigr)_{t\in\RR_+}$ \ is
 monotone increasing and convergence in probability implies the existence of a
 subsequence which converges almost surely.
Note that
 \begin{align}\label{help10}
  &\frac{T \sigma_1 \int_0^T \sqrt{Y_s} \, \dd W_s}{\int_0^T Y_s \, \dd s}
   =\frac{\frac{1}{T}(Y_T - y_0) - a}
         {\frac{1}{T^2} \int_0^T Y_s \, \dd s} , \qquad T \in \RR_{++} , \\
  &\frac{T \sigma_2 \int_0^T \sqrt{Y_s} \, \dd \tW_s}{\int_0^T Y_s \, \dd s}
   = \frac{\sigma_2 \varrho}{\sigma_1}
     \frac{T \sigma_1 \int_0^T \sqrt{Y_s} \, \dd W_s}
          {\int_0^T Y_s \, \dd s}
     + \frac{\sigma_2 \sqrt{1 - \varrho^2}}
            {\bigl(\frac{1}{T^2} \int_0^T Y_s \, \dd s\bigr)^{1/2}}
       \frac{\int_0^T \sqrt{Y_s} \, \dd B_s}
            {\bigl(\int_0^T Y_s \, \dd s\bigr)^{1/2}} ,
   \qquad T \in \RR_{++} . \label{tW2}
 \end{align}
Consequently, \eqref{aabb} will follow from
 \begin{align}\label{aabbt}
   \Biggl( \frac{\sigma_1 \int_0^T \frac{\dd W_s}{\sqrt{Y_s}}}
                {\bigl(\int_0^T \frac{\dd s}{Y_s}\bigr)^{1/2}} ,
           \frac{\sigma_2 \int_0^T \frac{\dd \tW_s}{\sqrt{Y_s}}}
                {\bigl(\int_0^T \frac{\dd s}{Y_s}\bigr)^{1/2}} ,
           \frac{\int_0^T \sqrt{Y_s} \, \dd B_s}
                {\bigl(\int_0^T Y_s \, \dd s\bigr)^{1/2}} ,
           \frac{1}{T} Y_T,
           \frac{1}{T^2} \int_0^T Y_s \, \dd s \Biggr)
    \distr \biggl( \bS^{1/2} \bZ_2, Z_3, \cY_1, \int_0^1 \cY_s \, \dd s \biggr)
 \end{align}
 as \ $T \to \infty$, \ where \ $Z_3$ \ is a standard normally distributed
 random variable independent of
 \ $\bigl(\bZ_2, \cY_1, \int_0^1 \cY_s \, \dd s\bigr)$, \ from \eqref{1/Y},
 \eqref{lim_int_1/Y}, \eqref{help10}, \eqref{tW2}, Slutsky's lemma, continuous
 mapping theorem, and
 \ $\PP\bigl(\int_0^1 \cY_s \, \dd s \in \RR_{++}\bigr) = 1$ \ (which has been
 shown in the proof of Theorem 3.1 in Barczy et al.~\cite{BarDorLiPap}).
Indeed,
 \[
   \begin{bmatrix}
    \sqrt{\log T} (\ha_T - a) \\
    \sqrt{\log T} (\halpha_T - \alpha) \\
    T \hb_T \\
    T(\hbeta_T-\beta)
   \end{bmatrix}
   \distr
   \frac{1}{1 - \frac{1}{\int_0^1 \cY_s \, \dd s} \cdot 0}
   \begin{bmatrix}
    \frac{1}{\bigl(a - \frac{\sigma_1^2}{2}\bigr)^{-1/2}} (\bS^{1/2} \bZ_2)_1
    - \frac{0}{\bigl(a - \frac{\sigma_1^2}{2}\bigr)^{-1}}
      \frac{\cY_1 - a}{\int_0^1 \cY_s \, \dd s} \\
    \frac{1}{\bigl(a - \frac{\sigma_1^2}{2}\bigr)^{-1/2}} (\bS^{1/2} \bZ_2)_2
    - \frac{0}{\bigl(a - \frac{\sigma_1^2}{2}\bigr)^{-1}}
      \frac{\cX_1 - \alpha}{\int_0^1 \cY_s \, \dd s} \\
    \frac{1}{\int_0^1 \cY_s \, \dd s} \cdot 0 \cdot (\bS^{1/2} \bZ_2)_1
    - \frac{\cY_1 - a}{\int_0^1 \cY_s \, \dd s} \\
    \frac{1}{\int_0^1 \cY_s \, \dd s} \cdot 0 \cdot (\bS^{1/2} \bZ_2)_2
    - \frac{\cX_1 - \alpha}{\int_0^1 \cY_s \, \dd s}
   \end{bmatrix}
 \]
 as \ $T \to \infty$, \ where
 \ $\bS^{1/2} \bZ_2
    =: \bigl( (\bS^{1/2} \bZ_2)_1 , (\bS^{1/2} \bZ_2)_2 \bigr)^\top$,
 \ since
 \begin{align}\label{D=}
   \Biggl( \bZ_2, \cY_1, \int_0^1 \cY_s \, \dd s,
           \frac{\sigma_2 \varrho}{\sigma_1}
           \frac{\cY_1 - a}{\int_0^1 \cY_s \, \dd s}
           + \frac{\sigma_2 \sqrt{1 - \varrho^2}}
                  {\bigl(\int_0^1 \cY_s \, \dd s\bigr)^{1/2}}
             Z_3 \Biggr)
   \distre
   \biggl( \bZ_2, \cY_1, \int_0^1 \cY_s \, \dd s,
           \frac{\cX_1 - \alpha}{\int_0^1 \cY_s \, \dd s} \biggr) .
 \end{align}
The statement \eqref{D=} is equivalent to
 \begin{align}\label{D==}
   \Biggl( \cY_1, \int_0^1 \cY_s \, \dd s,
           \frac{\sigma_2 \varrho}{\sigma_1}
           \frac{\cY_1 - a}{\int_0^1 \cY_s \, \dd s}
           + \frac{\sigma_2 \sqrt{1 - \varrho^2}}
                  {\bigl(\int_0^1 \cY_s \, \dd s\bigr)^{1/2}}
             Z_3 \Biggr)
   \distre
   \biggl( \cY_1, \int_0^1 \cY_s \, \dd s,
           \frac{\cX_1 - \alpha}{\int_0^1 \cY_s \, \dd s} \biggr) ,
 \end{align}
 since \ $\bZ_2$ \ is independent of \ $(Z_3,\cY_1,\int_0^1\cY_s\,\dd s)$ \ and
 of \ $(\cY_1,\int_0^1\cY_s\,\dd s,\cX_1)$.
\ The equality of the distributions in \eqref{D==} follows from the equality
 of their characteristic functions.
Namely, for all \ $(q_1, q_2, r) \in \RR^3$ \ and
 \ $T \in \RR_{++}$,
 \begin{align*}
  &\EE\Biggl(\exp\Biggl\{\ii q_1 \cY_1
                         + \ii q_2 \int_0^1 \cY_s \, \dd s
                         + \ii r \left(\frac{\sigma_2 \varrho}{\sigma_1}
                                 \frac{\cY_1 - a}{\int_0^1 \cY_s \, \dd s}
                                 + \frac{\sigma_2 \sqrt{1 - \varrho^2}}
                                    {\bigl(\int_0^1 \cY_s \, \dd s\bigr)^{1/2}}
                                    Z_3 \right) \Biggr\}
             \Bigg| \, \cY_1, \, \int_0^1 \cY_s \, \dd s \Biggr)
 \end{align*}
 \begin{align*}
  &= \exp\left\{\ii q_1 \cY_1
                 + \ii q_2 \int_0^1 \cY_s \, \dd s
                 + \ii r \frac{\sigma_2 \varrho}{\sigma_1}
                         \frac{\cY_1 - a}{\int_0^1 \cY_s \, \dd s} \right\}
     \EE\Biggl(\exp\Biggl\{ \ii r \frac{\sigma_2 \sqrt{1 - \varrho^2}}
                                    {\bigl(\int_0^1 \cY_s \, \dd s\bigr)^{1/2}}
                                    Z_3 \Biggr\}
             \Bigg| \, \cY_1, \, \int_0^1 \cY_s \, \dd s \Biggr) \\
  &= \exp\left\{\ii q_1 \cY_1
                 + \ii q_2 \int_0^1 \cY_s \, \dd s
                 + \ii r \frac{\sigma_2 \varrho}{\sigma_1}
                         \frac{\cY_1 - a}{\int_0^1 \cY_s \, \dd s} \right\}
     \exp\left\{ - \frac{1}{2} r^2 \frac{\sigma_2^2 (1 - \varrho^2)}
                                    {\int_0^1 \cY_s \, \dd s} \right\} ,
 \end{align*}
 thus
 \begin{align*}
  &\EE\Biggl(\exp\Biggl\{\ii q_1 \cY_1
                         + \ii q_2 \int_0^1 \cY_s \, \dd s
                         + \ii r \Biggl(\frac{\sigma_2 \varrho}{\sigma_1}
                                 \frac{\cY_1 - a}{\int_0^1 \cY_s \, \dd s}
                                 + \frac{\sigma_2 \sqrt{1 - \varrho^2}}
                                    {\bigl(\int_0^1 \cY_s \, \dd s\bigr)^{1/2}}
                                    Z_3 \Biggr) \Biggr\} \Biggr) \\
  &\qquad
   = \EE\left(\exp\Biggl\{\ii q_1 \cY_1
                          + \ii q_2 \int_0^1 \cY_s \, \dd s
                          + \ii r \frac{\sigma_2 \varrho}{\sigma_1}
                            \frac{\cY_1 - a}{\int_0^1 \cY_s \, \dd s}
                          - \frac{1}{2} r^2
                            \frac{\sigma_2^2 (1 - \varrho^2)}
                                 {\int_0^1 \cY_s \, \dd s} \Biggr\} \right).
 \end{align*}
Further, by \eqref{help_limit_YX},
 \[
   \cX_1 - \alpha
   = \sigma_2
     \int_0^1
      \sqrt{\cY_s} \,
      (\varrho \, \dd \cW_s + \sqrt{1 - \varrho^2} \, \dd \cB_s)
   = \frac{\sigma_2 \varrho}{\sigma_1} (\cY_1 - a)
     + \sigma_2 \sqrt{1 - \varrho^2} \int_0^1 \sqrt{\cY_s} \, \dd \cB_s ,
 \]
 hence for all \ $(q_1, q_2, r) \in \RR^3$ \ and
 \ $T \in \RR_{++}$, \ we have
 \begin{align*}
  &\EE\Biggl(\exp\Biggl\{\ii q_1 \cY_1
                         + \ii q_2 \int_0^1 \cY_s \, \dd s
                         + \ii r \frac{\cX_1 - \alpha}
                                      {\int_0^1 \cY_s \, \dd s} \Biggr\}
             \Bigg| \, \cY_s , s \in [0, 1] \Biggr) \\
  &= \EE\Biggl(\exp\Biggl\{\ii q_1 \cY_1
                         + \ii q_2 \int_0^1 \cY_s \, \dd s
                         + \ii r \left( \frac{\sigma_2 \varrho}{\sigma_1}
          \frac{\cY_1 - a}{\int_0^1 \cY_s \, \dd s}
          + \frac{\sigma_2 \sqrt{1 - \varrho^2}}
                 {\int_0^1 \cY_s \, \dd s}
            \int_0^1 \sqrt{\cY_s} \, \dd \cB_s \right) \Biggr\}
             \Bigg| \, \cY_s , s \in [0, 1] \Biggr) \\
  &= \exp\left\{\ii q_1 \cY_1
                 + \ii q_2 \int_0^1 \cY_s \, \dd s
                 + \ii r \frac{\sigma_2 \varrho}{\sigma_1}
                         \frac{\cY_1 - a}{\int_0^1 \cY_s \, \dd s} \right\} \\
  &\quad
     \times
     \EE\Biggl(\exp\Biggl\{\ii r \frac{\sigma_2 \sqrt{1 - \varrho^2}}
                                      {\int_0^1 \cY_s \, \dd s}
                                 \int_0^1 \sqrt{\cY_s} \, \dd \cB_s
                                     \Biggr\}
             \Bigg| \, \cY_s , s \in [0, 1] \Biggr) \\
  &= \exp\left\{\ii q_1 \cY_1
                 + \ii q_2 \int_0^1 \cY_s \, \dd s
                 + \ii r \frac{\sigma_2 \varrho}{\sigma_1}
                         \frac{\cY_1 - a}{\int_0^1 \cY_s \, \dd s} \right\}
     \exp\left\{ - \frac{1}{2} r^2 \frac{\sigma_2^2 (1 - \varrho^2)}
                                        {\int_0^1 \cY_s \, \dd s} \right\} ,
 \end{align*}
 where the last equality follows from the independence of \ $(\cY_t)_{t\in\RR_+}$
 \ and \ $(\cB_t)_{t\in\RR_+}$ \ yielding that the conditional distribution of
 \ $\int_0^1 \sqrt{\cY_s} \, \dd \cB_s$ \ given \ $(\cY_s)_{s\in[0,1]}$ \ is
 normal.
Thus
 \begin{align*}
  &\EE\Biggl(\exp\Biggl\{\ii q_1 \cY_1
                         + \ii q_2 \int_0^1 \cY_s \, \dd s
                         + \ii r
                           \frac{\cX_1 - \alpha}
                                {\int_0^1 \cY_s \, \dd s}\Biggr\}\Biggr) \\
  &\qquad
   = \EE\left(\exp\Biggl\{\ii q_1 \cY_1
                          + \ii q_2 \int_0^1 \cY_s \, \dd s
                          + \ii r \frac{\sigma_2 \varrho}{\sigma_1}
                            \frac{\cY_1 - a}{\int_0^1 \cY_s \, \dd s}
                          - \frac{1}{2} r^2
                            \frac{\sigma_2^2 (1 - \varrho^2)}
                                 {\int_0^1 \cY_s \, \dd s} \Biggr\} \right) ,
 \end{align*}
 and hence we obtain \eqref{D==}.

Now we turn to prove \eqref{aabbt}.
Using that
 \begin{align}
  \sigma_2 \int_0^T \frac{\dd \tW_s}{\sqrt{Y_s}}
  = \sigma_2 \varrho \int_0^T \frac{\dd W_s}{\sqrt{Y_s}}
    + \sigma_2 \sqrt{1 - \varrho^2} \int_0^T \frac{\dd B_s}{\sqrt{Y_s}} ,
  \qquad T \in \RR_{++} , \label{tW1}
 \end{align}
 and
 \begin{align}\label{S_1/2_helyett}
      \begin{bmatrix}
       \sigma_1 & 0 \\
       \sigma_2\rho & \sigma_2\sqrt{1-\rho^2} \\
      \end{bmatrix}
      \begin{bmatrix}
       \sigma_1 & 0 \\
       \sigma_2\rho & \sigma_2\sqrt{1-\rho^2} \\
      \end{bmatrix}^\top
      = \bS,
 \end{align}
 by continuous mapping theorem, to prove \eqref{aabbt}, it is sufficient to
 verify
 \begin{align}\label{aabb1}
  \Biggl( \frac{\int_0^T \frac{\dd W_s}{\sqrt{Y_s}}}
               {\bigl(\int_0^T \frac{\dd s}{Y_s}\bigr)^{1/2}} ,
          \frac{\int_0^T \frac{\dd B_s}{\sqrt{Y_s}}}
               {\bigl(\int_0^T \frac{\dd s}{Y_s}\bigr)^{1/2}} ,
          \frac{\int_0^T \sqrt{Y_s} \, \dd B_s}
               {\bigl(\int_0^T Y_s \, \dd s\bigr)^{1/2}} ,
          \frac{1}{T} Y_T,
          \frac{1}{T^2} \int_0^T Y_s \, \dd s \Biggr)
  \distr \biggl( \bZ_2, Z_3, \cY_1, \int_0^1 \cY_s \, \dd s \biggr)
 \end{align}
 as \ $T \to \infty$.
\ First we prove
 \begin{align}\label{YXintY}
  \biggl( \frac{1}{T} Y_T, \frac{1}{T^2} \int_0^T Y_s \, \dd s \biggr)
  \distr \biggl( \cY_1, \int_0^1 \cY_s \, \dd s \biggr) \qquad
  \text{as \ $T \to \infty$.}
 \end{align}
By part (ii) of Remark 2.7 in Barczy et al.\ \cite{BarDorLiPap}, we have
 \[
   \Bigl( \frac{1}{T} \cY_{Tt} \Bigr)_{t\in\RR_+}
   \distre \left( \cY_t \right)_{t\in\RR_+} \qquad
   \text{for all \ $T \in \RR_{++}$.}
 \]
Indeed, by Proposition \ref{Pro_Heston}, \ $(\cY_t)_{t\in\RR_+}$ \ is a regular
 affine process, and the so-called admissible set of parameters corresponding
 to \ $(\cY_t)_{t\in\RR_+}$ \ takes the form
 \ $(0, \frac{1}{2} \sigma_1^2, a, 0, 0, 0)$, \ and then part (ii) of Remark
 2.7 in Barczy et al.\ \cite{BarDorLiPap} can be applied.
Hence, by Lemma \ref{Lem_intf_cont}, we obtain
 \[
   \biggl( \cY_1, \int_0^1 \cY_s \, \dd s \biggr)
   \distre
   \biggl( \frac{1}{T} \cY_T, \frac{1}{T^2} \int_0^T \cY_s \, \dd s \biggr)
   \qquad \text{for all \ $T \in \RR_{++}$.}
 \]
Then, by Slutsky's lemma, in order to prove \eqref{YXintY}, it suffices to
 show convergences
 \begin{align}\label{beta2}
  \frac{1}{T} (Y_T - \cY_T) \stoch 0 , \qquad
  \frac{1}{T^2} \int_0^T (Y_s - \cY_s) \, \dd s \stoch 0, \qquad
  \text{as \ $T \to \infty$.}
 \end{align}
By (3.21) in Barczy et al.\ \cite{BarDorLiPap}, we have
 \begin{align}\label{E|Y-tY|}
  \EE(|Y_t - \cY_t|) \leq y_0 , \qquad t \in \RR_+ ,
 \end{align}
 hence
 \begin{align*}
  &\EE\left(\left|\frac{1}{T} (Y_T - \cY_T)\right|\right)
   \leq \frac{1}{T} y_0
   \to 0 , \\
  &\EE\left(\left|\frac{1}{T^2} \int_0^T (Y_s - \cY_s) \, \dd s\right|\right)
   \leq \frac{1}{T^2} \int_0^T \EE(|Y_s - \cY_s|) \, \dd s
   \leq \frac{1}{T} y_0
   \to 0 ,
 \end{align*}
 as \ $T \to \infty$ \ implying \eqref{beta2}.
Thus we conclude \eqref{YXintY}.

We will prove \eqref{aabb1} using continuity theorem.
Applying \eqref{logY}, one can write
 \begin{align}\label{help11}
   \sigma_1 \int_0^T \frac{\dd W_s}{\sqrt{Y_s}}
   = \log Y_T - \log y_0
     + \biggl( \frac{\sigma_1^2}{2} - a\biggr)
       \int_0^T \frac{\dd s}{Y_s} , \qquad T \in \RR_{++} ,
 \end{align}
 hence \ $\int_0^T \frac{\dd W_s}{\sqrt{Y_s}}$ \ is measurable with respect to
 the $\sigma$-algebra \ $\sigma(Y_s , s \in [0, T])$.
\ For all \ $(u_1, u_2, u_3, v_1, v_2) \in \RR^5$ \ and \ $T \in \RR_{++}$, \ we
 have
 \begin{align*}
  &\EE\Biggl( \exp\Biggl\{ \ii u_1
                           \frac{\int_0^T \frac{\dd W_s}{\sqrt{Y_s}}}
                                {\bigl(\int_0^T \frac{\dd s}{Y_s}\bigr)^{1/2}}
                           + \ii u_2
                             \frac{\int_0^T \frac{\dd B_s}{\sqrt{Y_s}}}
                                  {\bigl(\int_0^T
                                         \frac{\dd s}{Y_s}\bigr)^{1/2}}
                           + \ii u_3
                             \frac{\int_0^T \sqrt{Y_s} \, \dd B_s}
                                  {\bigl(\int_0^T Y_s \, \dd s\bigr)^{1/2}} \\
  &\phantom{\EE\Biggl( \exp\Biggl\{}
                           + \ii v_1 \frac{1}{T} Y_T
                           + \ii v_2 \frac{1}{T^2}
                             \int_0^T Y_s \, \dd s \biggr\}
              \Bigg| \, Y_s , s \in [0, T] \Biggr) \\
  &= \exp\Biggl\{ \ii u_1
                  \frac{\int_0^T \frac{\dd W_s}{\sqrt{Y_s}}}
                  {\bigl(\int_0^T \frac{\dd s}{Y_s}\bigr)^{1/2}}
                  + \ii v_1 \frac{1}{T} Y_T
                  + \ii v_2 \frac{1}{T^2} \int_0^T Y_s \, \dd s \Biggr\} \\
  &\quad
     \times
     \EE\Biggl( \exp\Biggl\{ \ii
                             \int_0^T
                              \Biggl( \frac{u_2}
                                 {\bigl(\int_0^T \frac{\dd t}{Y_t}\bigr)^{1/2}}
                                 \frac{1}{\sqrt{Y_s}}
                                 + \frac{u_3}
                                      {\bigl(\int_0^T Y_t \, \dd t\bigr)^{1/2}}
                                   \sqrt{Y_s} \Biggr) \dd B_s \Biggr\}
              \Bigg| \, Y_s , s \in [0, T] \Biggr) \\
  &= \exp\Biggl\{ \ii u_1
                  \frac{\int_0^T \frac{\dd W_s}{\sqrt{Y_s}}}
                  {\bigl(\int_0^T \frac{\dd s}{Y_s}\bigr)^{1/2}}
                  + \ii v_1 \frac{1}{T} Y_T
                  + \ii v_2 \frac{1}{T^2} \int_0^T Y_s \, \dd s \Biggr\} \\
  &\quad
     \times
     \exp\Biggl\{ - \frac{1}{2}
                    \int_0^T
                     \Biggl(\frac{u_2^2}
                                 {\int_0^T \frac{\dd t}{Y_t}}
                            \frac{1}{Y_s}
                            + \frac{u_3^2}
                                   {\int_0^T Y_t \, \dd t}
                              Y_s
                            + \frac{2 u_2 u_3}
                                   {\bigl(\int_0^T \frac{\dd t}{Y_t}
                                          \int_0^T Y_t \, \dd t\bigr)^{1/2}}
                     \Biggr) \dd s \Biggr\} \\
  &= \exp\Biggl\{ \ii u_1
                  \frac{\int_0^T \frac{\dd W_s}{\sqrt{Y_s}}}
                  {\bigl(\int_0^T \frac{\dd s}{Y_s}\bigr)^{1/2}}
                  + \ii v_1 \frac{1}{T} Y_T
                  + \ii v_2 \frac{1}{T^2} \int_0^T Y_s \, \dd s \Biggr\}
     \exp\Biggl\{ - \frac{1}{2} (u_2^2 + u_3^2)
                  - \frac{T u_2 u_3}
                         {\bigl(\int_0^T \frac{\dd t}{Y_t}
                                \int_0^T Y_t \, \dd t\bigr)^{1/2}} \Biggr\} ,
 \end{align*}
 where we used the independence of \ $Y$ \ and \ $B$.
\ Consequently, the joint characteristic function of the random vector on the
 left hand side of \eqref{aabb1} takes the form
 \begin{align*}
  &\EE\Biggl( \exp\Biggl\{ \ii u_1
                         \frac{\int_0^T \frac{\dd W_s}{\sqrt{Y_s}}}
                              {\bigl(\int_0^T \frac{\dd s}{Y_s}\bigr)^{1/2}}
                         + \ii u_2
                           \frac{\int_0^T \frac{\dd B_s}{\sqrt{Y_s}}}
                                {\bigl(\int_0^T \frac{\dd s}{Y_s}\bigr)^{1/2}}
                           + \ii u_3
                             \frac{\int_0^T \sqrt{Y_s} \, \dd B_s}
                                  {\bigl(\int_0^T Y_s \, \dd s\bigr)^{1/2}}
                         + \ii v_1 \frac{1}{T} Y_T
                         + \ii v_2 \frac{1}{T^2} \int_0^T Y_s \, \dd s \Biggr\}
      \Biggr) \\
  &\qquad\qquad
   = \ee^{-(u_2^2+ u_3^2)/2}
     \EE\Biggl(\exp\Biggl\{ \xi_T(u_1, v_1, v_2)
                           - \frac{T u_2 u_3}
                                  {\bigl(\int_0^T \frac{\dd t}{Y_t}
                                         \int_0^T Y_t \, \dd t\bigr)^{1/2}}
                   \Biggr\}\Biggr) ,
 \end{align*}
 where
 \begin{align*}
  \xi_T(u_1, v_1, v_2)
  := \ii u_1 \frac{\int_0^T \frac{\dd W_s}{\sqrt{Y_s}}}
                  {\bigl(\int_0^T \frac{\dd s}{Y_s}\bigr)^{1/2}}
     + \ii v_1 \frac{1}{T} Y_T
     + \ii v_2 \frac{1}{T^2} \int_0^T Y_s \, \dd s .
 \end{align*}
Ben Alaya and Kebaier \cite[proof of Theorem 6]{BenKeb2} proved
 \begin{equation}\label{BKconv}
   \Biggl( \frac{\log Y_T - \log y_0
                 + \bigl( \frac{\sigma_1^2}{2} - a \bigr)
                   \int_0^T \frac{\dd s}{Y_s}}
                {\sqrt{\log T}} ,
           \frac{Y_T}{T}, \frac{1}{T^2} \int_0^T Y_s \, \dd s \Biggr)
   \distr \left( \frac{\sigma_1}{\sqrt{a-\frac{\sigma_1^2}{2}}} Z_1 ,
                 \cY_1 , \int_0^1 \cY_s \, \dd s \right)
 \end{equation}
 as \ $T \to \infty$, \ where \ $Z_1$ \ is a $1$-dimensional standard normally
 distributed random variable independent of
 \ $\left(\cY_1, \int_0^1 \cY_t \, \dd t\right)$.
\ Using \eqref{help11} we have
 \[
   \frac{\int_0^T \frac{\dd W_s}{\sqrt{Y_s}}}
        {\bigl(\int_0^T \frac{\dd s}{Y_s}\bigr)^{1/2}}
   = \frac{\frac{1}{\sqrt{\log T}}\frac{1}{\sigma_1}
     \left( \log Y_T - \log y_0
            + \left(\frac{\sigma_1^2}{2} - a\right)
              \int_0^T\frac{\dd s}{Y_s} \right)}
          {\left(\frac{1}{\log T} \int_0^T\frac{\dd s}{Y_s}\right)^{1/2}} ,
  \qquad T \in \RR_{++} ,
 \]
 and, by \eqref{1/Y} and \eqref{BKconv}, we conclude
 \begin{equation}\label{BK}
  \Biggl( \frac{\int_0^T \frac{\dd W_s}{\sqrt{Y_s}}}
               {\bigl(\int_0^T \frac{\dd s}{Y_s}\bigr)^{1/2}} ,
          \frac{Y_T}{T}, \frac{1}{T^2} \int_0^T Y_s \, \dd s \Biggr)
  \distr \biggl( Z_1 , \cY_1 , \int_0^1 \cY_s \, \dd s \biggr) \qquad
  \text{as \ $T \to \infty$,}
 \end{equation}
 thus we derived joint convergence of three coordinates of the left hand side
 of \eqref{aabb1}.
Hence
 \begin{equation}\label{critxiconv}
  \EE(\exp\{\xi_T(u_1, v_1, v_2)\})
  \to
  \EE\biggl(\exp\biggl\{ \ii u_1 Z_1
                         + \ii v_1 \cY_1
                         + \ii v_2 \int_0^1 \cY_s \, \dd s \biggr\}\biggr)
  \qquad \text{as \ $T \to \infty$}
 \end{equation}
 for all \ $(u_1,v_1,v_2)\in\RR^3$.
\ Using \ $|\exp\{\xi_T(u_1, v_1, v_2)\}| = 1$, \ we have
 \begin{align*}
  &\Biggl|\EE\Biggl(\exp\Biggl\{ \xi_T(u_1, v_1, v_2)
                                 - \frac{T u_2 u_3}
                                        {\bigl(\int_0^T \frac{\dd t}{Y_t}
                                              \int_0^T Y_t \, \dd t\bigr)^{1/2}}
                        \Biggr\}\Biggr)
          -\EE(\exp\{\xi_T(u_1, v_1, v_2)\})\Biggr| \\
  &\qquad
   \leq \EE\Biggl(|\exp\{\xi_T(u_1, v_1, v_2)\}|
                  \Biggl|\exp\Biggl\{-\frac{T u_2 u_3}
                                           {\bigl(\int_0^T \frac{\dd t}{Y_t}
                                              \int_0^T Y_t \, \dd t\bigr)^{1/2}}
                             \Biggr\}
                         -1\Biggr|\Biggr) \\
  &\qquad
   = \EE\Biggl(\Biggl|\exp\Biggl\{-\frac{T u_2 u_3}
                                        {\bigl(\int_0^T \frac{\dd t}{Y_t}
                                              \int_0^T Y_t \, \dd t\bigr)^{1/2}}
                          \Biggr\}
                      -1\Biggr|\Biggr)
   \to 0 \qquad \text{as \ $T \to \infty$,}
 \end{align*}
 by the moment convergence theorem (see, e.g., Stroock \cite[Lemma 2.2.1]{Str}).
Indeed, by \eqref{lim_int_1/Y}, \eqref{BK}, continuous mapping theorem and
 Slutsky's lemma,
 \[
   \left| \exp\Biggl\{-\frac{T u_2 u_3}
                            {\bigl(\int_0^T \frac{\dd t}{Y_t}
                                   \int_0^T Y_t \, \dd t\bigr)^{1/2}}\Biggr\}
          - 1 \right|
   = \left| \exp\Biggl\{-\frac{u_2 u_3}
                              {\bigl(\int_0^T \frac{\dd t}{Y_t}
                                     \cdot \frac{1}{T^2}
                                     \int_0^T Y_t \, \dd t\bigr)^{1/2}}\Biggr\}
            - 1 \right|
   \stoch 0 \qquad \text{as \ $T \to \infty$,}
 \]
 and the family
 \[
   \left\{ \left| \exp\Biggl\{-\frac{T u_2 u_3}
                                    {\bigl(\int_0^T \frac{\dd t}{Y_t}
                                           \int_0^T
                                            Y_t \, \dd t\bigr)^{1/2}}\Biggr\}
                  - 1 \right| , \, \, T \in \RR_{++} \right\}
 \]
 is uniformly integrable, since, by Cauchy--Schwarz inequality,
 \[
   \Biggl| \exp\Biggl\{-\frac{T u_2 u_3}
                             {\bigl(\int_0^T \frac{\dd t}{Y_t}
                                    \int_0^T Y_t \, \dd t\bigr)^{1/2}}\Biggr\}
          - 1\Biggr|^2
   \leq \left( \exp\Biggl\{\frac{T |u_2 u_3|}
                                {\bigl(\int_0^T \frac{\dd t}{Y_t}
                                      \int_0^T Y_t \, \dd t\bigr)^{1/2}}\Biggr\}
               + 1 \right)^2
   \leq (\exp\{|u_2 u_3|\} + 1)^2
 \]
 for all \ $T \in \RR_{++}$.
\ Using \eqref{critxiconv}, we conclude
 \begin{align*}
  &\EE\Biggl(\exp\Biggl\{\ii u_1
                         \frac{\int_0^T \frac{\dd W_s}{\sqrt{Y_s}}}
                              {\bigl(\int_0^T \frac{\dd s}{Y_s}\bigr)^{1/2}}
                         + \ii u_2
                           \frac{\int_0^T \frac{\dd B_s}{\sqrt{Y_s}}}
                                {\bigl(\int_0^T \frac{\dd s}{Y_s}\bigr)^{1/2}}
                         + \ii u_3
                           \frac{\int_0^T \sqrt{Y_s} \, \dd B_s}
                                {\bigl(\int_0^T Y_s \, \dd s\bigr)^{1/2}}
                         + \ii v_1 \frac{1}{T} Y_T
                         + \ii v_2 \frac{1}{T^2}
                                   \int_0^T Y_s \, \dd s\Biggr\}\Biggr) \\
  &\qquad
   \to \ee^{-(u_2^2+u_3^2)/2}
       \EE\biggl(\exp\Biggl\{ \ii u_1 Z_1
                              + \ii v_1 \cY_1
                              + \ii v_2 \int_0^1 \cY_s \, \dd s \Biggr\}
          \biggr) \qquad \text{as \ $T \to \infty$.}
 \end{align*}
Note that, since \ $Z_1$ \ is independent of
 \ $\bigl( \cY_1, \int_0^1 \cY_s \, \dd s \bigr)$, \ we have
 \begin{align*}
  &\ee^{-(u_2^2+u_3^2)/2}
   \EE\biggl(\exp\biggl\{ \ii u_1 Z_1
                          + \ii v_1 \cY_1
                          + \ii v_2 \int_0^1 \cY_s \, \dd s \biggr\} \biggr) \\
  &\qquad
   = \EE(\ee^{\ii u_1 Z_1})\EE(\ee^{\ii u_2 Z_2}) \EE(\ee^{\ii u_3 Z_3})
     \EE\biggl(\exp\biggl\{ \ii v_1 \cY_1
                            + \ii v_2 \int_0^1 \cY_s \, \dd s \biggr\}\biggr) ,
 \end{align*}
 where \ $(Z_2, Z_3)$ \ is a
 2-dimensional standard normally distributed random vector, independent of
 \ $\bigl( Z_1, \cY_1, \int_0^1 \cY_s \, \dd s \bigr)$, \ hence we obtain
 \eqref{aabb1} with \ $\bZ_2 := (Z_1, Z_2)$.
\proofend

\begin{Rem}\label{Rem_comparison6}
(i) As a consequence of Theorem \ref{Thm_MLE} we get back the description of
 the asymptotic behavior of the MLE of \ $(a, b)$ \ for the CIR process
 \ $(Y_t)_{t\in\RR_+}$ \ in the critical case whenever
 \ $a \in \bigl(\frac{\sigma_1^2}{2}, \infty\bigr)$ \ proved by Ben Alaya and
 Kebaier \cite[Theorem 6, part 2]{BenKeb2}.
We note that Ben Alaya and Kebaier \cite[Theorem 6, part 1]{BenKeb2} described
 the asymptotic behavior of the MLE of \ $(a, b)$ \ in the critical case
 for the CIR process \ $(Y_t)_{t\in\RR_+}$ \ with \ $a = \frac{\sigma_1^2}{2}$
 \ as well.

\noindent
(ii) Theorem \ref{Thm_MLE} does not cover the case
 \ $a = \frac{\sigma_1^2}{2}$, \ we renounce to consider it.

\noindent
(iii) Ben Alaya and Kebaier's proof of part 2 of their Theorem 6 relies on an
 explicit form of the moment generating-Laplace transform of the quadruplet
 \[
   \biggl( \log Y_t, Y_t, \int_0^t Y_s \, \dd s,
           \int_0^t \frac{\dd s}{Y_s} \biggr) ,
   \qquad t \in \RR_+ .
 \]
Using this explicit form, they derived convergence \eqref{BKconv}, which is a
 corner stone of the proof of our Theorem \ref{Thm_MLE}.
\proofend
\end{Rem}

The next theorem can be considered as a counterpart of Theorem \ref{Thm_MLE}
 by incorporating random scaling.

\begin{Thm}\label{Thm_MLE_random}
If \ $a \in \bigl( \frac{\sigma_1^2}{2}, \infty \bigr)$, \ $b = 0$,
 \ $\alpha, \beta \in \RR$, \ $\sigma_1, \sigma_2 \in \RR_{++}$,
 \ $\varrho \in (-1, 1)$ \ and
 \ $(Y_0, X_0) = (y_0, x_0) \in \RR_{++} \times \RR$, \ then
 \begin{align}\label{aabbr}
  \begin{bmatrix}
   \bigl(\int_0^T \frac{\dd s}{Y_s}\bigr)^{1/2} \, (\ha_T - a) \\[1mm]
   \bigl(\int_0^T \frac{\dd s}{Y_s}\bigr)^{1/2} \,
   (\halpha_T - \alpha) \\[1mm]
   \bigl(\int_0^T Y_s \, \dd s\bigr)^{1/2} \, \hb_T \\[1mm]
   \bigl(\int_0^T Y_s \, \dd s\bigr)^{1/2} \, (\hbeta_T - \beta)
  \end{bmatrix}
  \distr
  \begin{bmatrix}
   \bS^{1/2} \bZ_2 \\
   \frac{a - \cY_1}{\bigl(\int_0^1 \cY_s \, \dd s\bigr)^{1/2}} \\[1mm]
   \frac{\alpha - \cX_1}{\bigl(\int_0^1 \cY_s \, \dd s\bigr)^{1/2}}
  \end{bmatrix}
  \qquad \text{as \ $T \to \infty$,}
 \end{align}
 where \ $(\cY_t, \cX_t)_{t\in\RR_+}$ \ is the unique strong solution of the SDE
 \eqref{help_limit_YX} with initial value \ $(\cY_0, \cX_0) = (0, 0)$,
 \ $\bZ_2$ \ is a
 $2$-dimensional standard normally distributed random vector independent of
 \ $\bigl(\cY_1, \int_0^1 \cY_t \, \dd t, \cX_1\bigr)$, \ and
 \ $\bS$ \ is defined in \eqref{bSigma}.
\end{Thm}

\noindent{\bf Proof.}
By Lemma \ref{LEMMA_MLE_exist}, there exists a unique MLE
 \ $\bigl(\ha_T, \hb_T, \halpha_T, \hbeta_T\bigr)$ \ of
 \ $(a, b, \alpha, \beta)$ \ for all \ $T\in\RR_{++}$, \ which has the form
 given in \eqref{MLE}.
By \eqref{MLE-}, we have
 \[
   \Bigl(\int_0^T\frac{\dd s}{Y_s}\Bigr)^{1/2} (\ha_T - a)
   = \frac{\frac{\sigma_1 \int_0^T \frac{\dd W_s}{\sqrt{Y_s}}}
                {\left(\int_0^T \frac{\dd s}{Y_s}\right)^{1/2}}
           - \frac{1}{\left(\int_0^T \frac{\dd s}{Y_s}\right)^{1/2}}
             \frac{T \sigma_1 \int_0^T \sqrt{Y_s} \, \dd W_s}
                  {\int_0^T Y_s \, \dd s}}
          {1 - \frac{1}{\frac{1}{T^2} \int_0^T Y_s \, \dd s}
               \frac{1}{\int_0^T \frac{\dd s}{Y_s}}} ,
 \]
 \[
   \Bigl(\int_0^T Y_s\,\dd s\Bigr)^{1/2} \hb_T
   = \frac{\frac{1}
                {\left(\frac{1}{T^2} \int_0^T Y_s \, \dd s\right)^{1/2}}
           \frac{1}{\left(\int_0^T \frac{\dd s}{Y_s}\right)^{1/2}}
           \frac{\sigma_1 \int_0^T \frac{\dd W_s}{\sqrt{Y_s}}}
                {\left(\int_0^T \frac{\dd s}{Y_s}\right)^{1/2}}
           - \frac{\sigma_1 \int_0^T \sqrt{Y_s} \, \dd W_s}
                  {\left(\int_0^T Y_s \, \dd s\right)^{1/2}}}
          {1 - \frac{1}{\frac{1}{T^2} \int_0^T Y_s \, \dd s}
               \frac{1}{\int_0^T \frac{\dd s}{Y_s}}} ,
 \]
 \[
   \Bigl(\int_0^T\frac{\dd s}{Y_s}\Bigr)^{1/2} (\halpha_T - \alpha)
   = \frac{\frac{\sigma_2 \int_0^T \frac{\dd \tW_s}{\sqrt{Y_s}}}
                {\left(\int_0^T \frac{\dd s}{Y_s}\right)^{1/2}}
           - \frac{1}{\left(\int_0^T \frac{\dd s}{Y_s}\right)^{1/2}}
             \frac{T \sigma_2 \int_0^T \sqrt{Y_s} \, \dd \tW_s}
                  {\int_0^T Y_s \, \dd s}}
          {1 - \frac{1}{\frac{1}{T^2} \int_0^T Y_s \, \dd s}
               \frac{1}{\int_0^T \frac{\dd s}{Y_s}}} ,
 \]
 and
 \[
   \Bigl(\int_0^T Y_s\,\dd s\Bigr)^{1/2} (\hbeta_T - \beta)
   = \frac{\frac{1}
                {\left(\frac{1}{T^2} \int_0^T Y_s \, \dd s\right)^{1/2}}
           \frac{1}{\left(\int_0^T \frac{\dd s}{Y_s}\right)^{1/2}}
           \frac{\sigma_2 \int_0^T \frac{\dd \tW_s}{\sqrt{Y_s}}}
                {\left(\int_0^T \frac{\dd s}{Y_s}\right)^{1/2}}
           - \frac{\sigma_2 \int_0^T \sqrt{Y_s} \, \dd \tW_s}
                  {\left(\int_0^T Y_s \, \dd s\right)^{1/2}}}
          {1 - \frac{1}{\frac{1}{T^2} \int_0^T Y_s \, \dd s}
               \frac{1}{\int_0^T \frac{\dd s}{Y_s}}} ,
 \]
 provided that \ $\int_0^T Y_s \, \dd s \int_0^T \frac{1}{Y_s} \, \dd s > T^2$
 \ which holds a.s.
\ We have
 \begin{align}\label{root1}
  \frac{\sigma_1 \int_0^T \sqrt{Y_s} \, \dd W_s}
       {\bigl(\int_0^T Y_s \, \dd s\bigr)^{1/2}}
  &= \frac{Y_T - y_0 -aT}{\bigl(\int_0^T Y_s \, \dd s\bigr)^{1/2}}
   = \frac{\frac{1}{T}(Y_T - y_0) - a}
          {\bigl(\frac{1}{T^2}\int_0^T Y_s \, \dd s\bigr)^{1/2}} ,
   \qquad T \in \RR_{++} , \\
  \frac{\sigma_2 \int_0^T \sqrt{Y_s} \, \dd \tW_s}
       {\bigl(\int_0^T Y_s \, \dd s\bigr)^{1/2}}
  &= \frac{\sigma_2 \varrho}{\sigma_1}
     \frac{\sigma_1 \int_0^T \sqrt{Y_s} \, \dd W_s}
          {\bigl(\int_0^T Y_s \, \dd s\bigr)^{1/2}}
     + \sigma_2 \sqrt{1 - \varrho^2}
       \frac{\int_0^T \sqrt{Y_s} \, \dd B_s}
           {\bigl(\int_0^T Y_s \, \dd s\bigr)^{1/2}} ,
    \qquad T \in \RR_{++} , \label{root2}
 \end{align}
 hence \eqref{aabbr} follows from \eqref{1/Y}, \eqref{lim_int_1/Y},
 \eqref{help10}, \eqref{tW2}, \eqref{aabbt}, \eqref{D=}, \eqref{tW1},
 Slutsky's lemma, continuous
 mapping theorem, and \ $\PP(\int_0^1 \cY_s \, \dd s \in \RR_{++}) = 1$ \ (which
 has been shown in the proof of Theorem 3.1 in Barczy et
 al.~\cite{BarDorLiPap}).
Indeed,
 \[
  \begin{bmatrix}
   \bigl(\int_0^T \frac{\dd s}{Y_s}\bigr)^{1/2} (\ha_T - a) \\
   \bigl(\int_0^T \frac{\dd s}{Y_s}\bigr)^{1/2} (\halpha_T - \alpha) \\
   \bigl(\int_0^T Y_s \, \dd s\bigr)^{1/2} \hb_T \\
   \bigl(\int_0^T Y_s \, \dd s\bigr)^{1/2} (\hbeta_T - \beta)
  \end{bmatrix}
   \distr
   \frac{1}{1 - \frac{1}{\int_0^1 \cY_s \, \dd s} \cdot 0}
   \begin{bmatrix}
    (\bS^{1/2} \bZ_2)_1
    - 0 \cdot \frac{\cY_1 - a}{\int_0^1 \cY_s \, \dd s} \\
    (\bS^{1/2} \bZ_2)_2
    - 0 \cdot \frac{\cX_1 - \alpha}{\int_0^1 \cY_s \, \dd s} \\
    \frac{1}{\left(\int_0^1 \cY_s \, \dd s\right)^{1/2}}
    \cdot 0 \cdot (\bS^{1/2} \bZ_2)_1
    - \frac{\cY_1 - a}{\left(\int_0^1 \cY_s \, \dd s\right)^{1/2}} \\
    \frac{1}{\left(\int_0^1 \cY_s \, \dd s\right)^{1/2}}
    \cdot 0 \cdot (\bS^{1/2} \bZ_2)_2
    - \frac{\cX_1 - \alpha}{\left(\int_0^1 \cY_s \, \dd s\right)^{1/2}}
   \end{bmatrix}
 \]
 as \ $T \to \infty$, \ where
 \ $\bS^{1/2} \bZ_2 = \bigl( (\bS^{1/2} \bZ_2)_1 , (\bS^{1/2} \bZ_2)_2 \bigr)^\top$,
 since
 \begin{align*}
   \Biggl( \bZ_2, \cY_1, \int_0^1 \cY_s \, \dd s,
           \frac{\sigma_2 \varrho}{\sigma_1}
           \frac{\cY_1 - a}{\bigl(\int_0^1 \cY_s \, \dd s\bigr)^{1/2}}
           + \sigma_2 \sqrt{1 - \varrho^2} Z_3 \Biggr)
   \distre
   \Biggl( \bZ_2, \cY_1, \int_0^1 \cY_s \, \dd s,
           \frac{\cX_1 - \alpha}
                {\bigl(\int_0^1 \cY_s \, \dd s\bigr)^{1/2}} \Biggr) ,
 \end{align*}
 which can be shown in the same way as \eqref{D=}.
\proofend

\begin{Rem}\label{Rem_comparison3}
For a critical (i.e., \ $b=0$) CIR models with
 \ $a \in\bigl(\frac{\sigma_1^2}{2}, \infty\bigr)$, \ using random
 scaling, Overbeck \cite[Theorem 3, part (ii)]{Ove} has already
 described the asymptotic behaviour
 of \ $\ha_T$ \  and \ $\hb_T$ \ separately, but he did not
 consider their joint asymptotic behaviour.
\proofend
\end{Rem}

\section{Asymptotic behaviour of MLE: supercritical case}
\label{section_AMLE_supercritical}

We consider supercritical Heston models, i.e., when
 \ $b \in \RR_{--}$.

\begin{Thm}\label{Thm_MLE_super}
If \ $a \in \left[ \frac{\sigma_1^2}{2}, \infty \right)$, \ $b \in \RR_{--}$,
 \ $\alpha, \beta \in \RR$, \ $\sigma_1, \sigma_2 \in \RR_{++}$,
 \ $\varrho \in (-1, 1)$, \ and
 \ $(Y_0, X_0) = (y_0, x_0) \in \RR_{++} \times \RR$, \ then
 \begin{align}\label{aabb_super}
  \begin{bmatrix}
   \ha_T - a \\
   \halpha_T - \alpha \\
   \ee^{-bT/2} (\hb_T - b) \\
   \ee^{-bT/2} (\hbeta_T - \beta)
  \end{bmatrix}
  \distr
  \begin{bmatrix}
   \tcV \\
   \varrho \frac{\sigma_2}{\sigma_1} \tcV
   + \sigma_2 \sqrt{1 - \varrho^2}
     \left(\int_0^{-1/b} \tcY_u \, \dd u\right)^{-1/2} Z_1 \\[2mm]
     \left(-\frac{\tcY_{-1/b}}{b}\right)^{-1/2} \bS^{1/2}\bZ_2
  \end{bmatrix}
 \end{align}
 as \ $T \to \infty$, \ where \ $(\tcY_t)_{t\in\RR_+}$ \ is a CIR process given
 by the SDE
 \begin{align*}
  \dd \tcY_t = a \dd t + \sigma_1 \sqrt{\tcY_t} \, \dd \cW_t ,
  \qquad t \in \RR_+ ,
 \end{align*}
 with initial value \ $\tcY_0 = y_0$, \ where \ $(\cW_t)_{t\in\RR_+}$ \ is a
 standard Wiener process,
 \[
  \tcV
   :=\frac{\log \tcY_{-1/b} - \log y_0}{\int_0^{-1/b}\tcY_u \, \dd u}
      + \frac{\sigma_1^2}{2} - a,
 \]
 \ $Z_1$ \ is a $1$-dimensional standard normally distributed random variable,
 \ $\bZ_2$ \ is a $2$-dimensional standard normally distributed random vector
 such that \ $(\tcY_{-1/b},\int_0^{-1/b}\tcY_u \, \dd u)$, \ $Z_1$ \ and
 \ $\bZ_2$ \ are independent, and \ $\bS$ \ is defined in \eqref{bSigma}.

With a random scaling, we have
 \begin{align}\label{MLE_super_random}
   \begin{bmatrix}
    \ha_T - a \\
    \halpha_T - \alpha \\
    \left(\int_0^T Y_s \, \dd s\right)^{1/2} (\hb_T - b) \\
    \left(\int_0^T Y_s \, \dd s\right)^{1/2} (\hbeta_T - \beta)
   \end{bmatrix}
   \distr
   \begin{bmatrix}
    \tcV \\
    \varrho \frac{\sigma_2}{\sigma_1} \tcV
    + \sigma_2 \sqrt{1 - \varrho^2}
      \left(\int_0^{-1/b} \tcY_u \, \dd u\right)^{-1/2}  Z_1 \\[2mm]
    \bS^{1/2} \bZ_2
  \end{bmatrix}
 \end{align}
 as \ $T \to \infty$.
\end{Thm}

\noindent{\bf Proof.}
By Lemma \ref{LEMMA_MLE_exist}, there exists a unique MLE
 \ $\bigl(\ha_T, \hb_T, \halpha_T, \hbeta_T\bigr)$ \ of
 \ $(a, b, \alpha, \beta)$ \ for all \ $T\in\RR_{++}$, \ which has the form
 given in \eqref{MLE}.
By \eqref{MLE-} and
 \[
  \sigma_2 \int_0^T \frac{\dd \tW_s}{\sqrt{Y_s}}
  = \sigma_2 \varrho \int_0^T \frac{\dd W_s}{\sqrt{Y_s}}
    + \sigma_2 \sqrt{1 - \varrho^2} \int_0^T \frac{\dd B_s}{\sqrt{Y_s}} ,
 \]
 we obtain
 \[
   \ha_T - a
   = \frac{\frac{\sigma_1 \int_0^T \frac{\dd W_s}{\sqrt{Y_s}}}
                {\int_0^T \frac{\dd s}{Y_s}}
           - \frac{T \ee^{bT/2}}{\int_0^T \frac{\dd s}{Y_s}}
             \frac{1}
                  {\left(\ee^{bT} \int_0^T Y_s \, \dd s\right)^{1/2}}
             \frac{\sigma_1 \int_0^T \sqrt{Y_s} \, \dd W_s}
                  {\left(\int_0^T Y_s \, \dd s\right)^{1/2}}}
          {1 - \frac{T^2 \ee^{bT}}
                    {\ee^{bT} \int_0^T Y_s \, \dd s
                     \int_0^T \frac{\dd s}{Y_s}}} ,
 \]
 \[
   \halpha_T - \alpha
   = \frac{\frac{\sigma_2 \varrho \int_0^T \frac{\dd W_s}{\sqrt{Y_s}}}
                {\int_0^T \frac{\dd s}{Y_s}}
           + \frac{\sigma_2 \sqrt{1 - \varrho^2}}
                  {\left(\int_0^T \frac{\dd s}{Y_s}\right)^{1/2}}
             \frac{\int_0^T \frac{\dd B_s}{\sqrt{Y_s}}}
                  {\left(\int_0^T \frac{\dd s}{Y_s}\right)^{1/2}}
           - \frac{T \ee^{bT/2}}{\int_0^T \frac{\dd s}{Y_s}}
             \frac{1}
                  {\left(\ee^{bT} \int_0^T Y_s \, \dd s\right)^{1/2}}
             \frac{\sigma_2 \int_0^T \sqrt{Y_s} \, \dd \tW_s}
                  {\left(\int_0^T Y_s \, \dd s\right)^{1/2}}}
          {1 - \frac{T^2 \ee^{bT}}
                    {\ee^{bT} \int_0^T Y_s \, \dd s
                     \int_0^T \frac{\dd s}{Y_s}}} ,
 \]
 \[
   \ee^{-bT/2} (\hb_T - b)
   = \frac{\frac{T \ee^{bT/2}}{\ee^{bT} \int_0^T Y_s \, \dd s}
           \frac{\sigma_1 \int_0^T \frac{\dd W_s}{\sqrt{Y_s}}}
                {\int_0^T \frac{\dd s}{Y_s}}
           - \frac{1}
                  {\left(\ee^{bT} \int_0^T Y_s \, \dd s\right)^{1/2}}
             \frac{\sigma_1 \int_0^T \sqrt{Y_s} \, \dd W_s}
                  {\left(\int_0^T Y_s \, \dd s\right)^{1/2}}}
          {1 - \frac{T^2 \ee^{bT}}
                    {\ee^{bT} \int_0^T Y_s \, \dd s
                     \int_0^T \frac{\dd s}{Y_s}}} ,
 \]
 and
 \[
   \ee^{-bT/2} (\hbeta_T - \beta)
   = \frac{\frac{T \ee^{bT/2}}{\ee^{bT} \int_0^T Y_s \, \dd s}
           \left( \frac{\sigma_2 \varrho \int_0^T \frac{\dd W_s}{\sqrt{Y_s}}}
                       {\int_0^T \frac{\dd s}{Y_s}}
                  + \frac{\sigma_2 \sqrt{1 - \varrho^2}}
                         {\left(\int_0^T \frac{\dd s}{Y_s}\right)^{1/2}}
                    \frac{\int_0^T \frac{\dd B_s}{\sqrt{Y_s}}}
                         {\left(\int_0^T \frac{\dd s}{Y_s}\right)^{1/2}} \right)
           - \frac{1}
                  {\left(\ee^{bT} \int_0^T Y_s \, \dd s\right)^{1/2}}
             \frac{\sigma_2 \int_0^T \sqrt{Y_s} \, \dd \tW_s}
                  {\left(\int_0^T Y_s \, \dd s\right)^{1/2}}}
          {1 - \frac{T^2 \ee^{bT}}
                    {\ee^{bT} \int_0^T Y_s \, \dd s
                     \int_0^T \frac{\dd s}{Y_s}}} ,
 \]
 provided that \ $\int_0^T Y_s \, \dd s \int_0^T \frac{1}{Y_s} \, \dd s > T^2$
 \ which holds a.s.
\ Applying \eqref{logY}, one can write
 \[
   \sigma_1 \int_0^T \frac{\dd W_s}{\sqrt{Y_s}}
   = \log Y_T - \log y_0
     + \left( \frac{\sigma_1^2}{2} - a\right)
       \int_0^T \frac{\dd s}{Y_s} + b T , \qquad T \in \RR_{++} ,
 \]
 thus, by \eqref{lim_Y} and \eqref{lim_1/Y},
 \begin{equation}\label{a_super}
   \frac{\sigma_1 \int_0^T \frac{\dd W_s}{\sqrt{Y_s}}}
         {\int_0^T \frac{\dd s}{Y_s}}
   = \frac{\log(\ee^{bT} Y_T) - \log y_0}{\int_0^T \frac{\dd s}{Y_s}}
      + \frac{\sigma_1^2}{2} - a
   \as \frac{\log V - \log y_0}{\int_0^\infty \frac{\dd s}{Y_s}}
        + \frac{\sigma_1^2}{2} - a
 \end{equation}
 as \ $T \to \infty$.
\ By Theorem 4 in Ben Alaya and Kebaier \cite{BenKeb2},
 \[
   \frac{\log V - \log y_0}{\int_0^\infty \frac{\dd s}{Y_s}}
   + \frac{\sigma_1^2}{2} - a
   \distre
   \frac{\log \tcY_{-1/b} - \log y_0}{\int_0^{-1/b}\tcY_u \, \dd u}
   + \frac{\sigma_1^2}{2} - a
   =: \tcV .
 \]
Moreover, \eqref{lim_intY} and \eqref{lim_1/Y} yield
 \begin{align}\label{help8}
  \frac{T^2 \ee^{bT}}{\ee^{bT} \int_0^T Y_s \, \dd s \int_0^T \frac{\dd s}{Y_s}}
  &\as \frac{0}{\left(-\frac{V}{b}\right) \int_0^\infty \frac{\dd s}{Y_s}}
   = 0 \qquad \text{as \ $T \to \infty$,} \\
  \frac{T\ee^{bT/2}}{\ee^{bT} \int_0^T Y_s \, \dd s}
  &\as \frac{0}{-\frac{V}{b}}
   =0 \qquad \text{as \ $T \to \infty$.} \label{help9}
 \end{align}
Consequently, \eqref{aabb_super} will follow from
 \begin{equation}\label{aabbt_super}
  \begin{aligned}
   &\left( \frac{\int_0^T \frac{\dd B_s}{\sqrt{Y_s}}}
                {\left(\int_0^T \frac{\dd s}{Y_s}\right)^{1/2}} ,
           \frac{\sigma_1 \int_0^T \sqrt{Y_s} \, \dd W_s}
                {\left(\int_0^T Y_s \, \dd s\right)^{1/2}} ,
           \frac{\sigma_2 \int_0^T \sqrt{Y_s} \, \dd \tW_s}
                {\left(\int_0^T Y_s \, \dd s\right)^{1/2}} ,
           \ee^{bT} Y_T , \ee^{bT} \int_0^T Y_s \, \dd s ,
           \int_0^T \frac{\dd s}{Y_s} \right) \\
    &\distr \left( Z_1, \bS^{1/2} \bZ_2, \tcY_{-1/b}, -\frac{\tcY_{-1/b}}{b} ,
                   \int_0^{-1/b} \tcY_u \, \dd u \right)
     \qquad \text{as \ $T \to \infty$,}
  \end{aligned}
 \end{equation}
 from \eqref{lim_intY}, \eqref{a_super}, \eqref{help8}, \eqref{help9},
 Slutsky's lemma, continuous mapping theorem and \
 $\PP(\tcY_{-1/b}\in \RR_{++}) = 1$,
 \ $\PP(\int_0^{-1/b} \tcY_u \, \dd u \in\RR_{++}) = 1$ \ (due to
 \ $\PP(\tcY_t\in\RR_{++}, \, \forall \, t \in \RR_+) = 1$).
\ Indeed,
 \begin{align*}
  &\begin{bmatrix}
    \ha_T - a \\
    \halpha_T - \alpha \\
    \ee^{-bT/2} (\hb_T - b) \\
    \ee^{-bT/2} (\hbeta_T - \beta)
   \end{bmatrix}
   \distr
 \end{align*}
 \begin{align*}
  &\qquad\distr
   \frac{1}{1 - \frac{0}{-\frac{\tcY_{-1/b}}{b} \int_0^{-1/b} \tcY_u \, \dd u}}
   \begin{bmatrix}
    \tcV
    - \frac{0}{\int_0^{-1/b} \tcY_u \, \dd u}
      \frac{1}{\left(-\frac{\tcY_{-1/b}}{b}\right)^{1/2}}
      (\bS^{1/2} \bZ_2)_1 \\[3mm]
    \varrho \frac{\sigma_2}{\sigma_1} \tcV
    + \frac{\sigma_2 \sqrt{1 - \varrho^2}}
           {\left(\int_0^{-1/b} \tcY_u \, \dd u\right)^{1/2}} Z_1
    - \frac{0}{\int_0^{-1/b} \tcY_u \, \dd u}
      \frac{1}{\left(-\frac{\tcY_{-1/b}}{b}\right)^{1/2}}
      (\bS^{1/2} \bZ_2)_2 \\[3mm]
    \frac{0}{-\frac{\tcY_{-1/b}}{b}} \tcV
    - \frac{1}{\left(-\frac{\tcY_{-1/b}}{b}\right)^{1/2}} (\bS^{1/2} \bZ_2)_1 \\[3mm]
    \frac{0}{-\frac{\tcY_{-1/b}}{b}}
    \left( \varrho \frac{\sigma_2}{\sigma_1} \tcV
           + \frac{\sigma_2 \sqrt{1 - \varrho^2}}
                  {\left(\int_0^{-1/b} \tcY_u \, \dd u\right)^{1/2}} Z_1 \right)
    - \frac{1}{\left(-\frac{\tcY_{-1/b}}{b}\right)^{1/2}} (\bS^{1/2} \bZ_2)_2
   \end{bmatrix}
 \end{align*}
 as \ $T \to \infty$, \ where
 \ $\bS^{1/2} \bZ_2 = \bigl( (\bS^{1/2} \bZ_2)_1 , (\bS^{1/2} \bZ_2)_2 \bigr)^\top$.

Using that
 \[
  \sigma_2\int_0^T \sqrt{Y_s} \, \dd \tW_s
    = \sigma_2\varrho \int_0^T \sqrt{Y_s} \, \dd W_s
       +\sigma_2\sqrt{1-\varrho^2} \int_0^T \sqrt{Y_s} \, \dd B_s,
       \qquad T\in\RR_{+},
 \]
 and \eqref{S_1/2_helyett}, by continuous mapping theorem, to prove
 \eqref{aabbt_super}, it is sufficient to verify
 \begin{equation}\label{aabb_super_1}
  \begin{aligned}
   &\left( \frac{\int_0^T \frac{\dd B_s}{\sqrt{Y_s}}}
                {\left(\int_0^T \frac{\dd s}{Y_s}\right)^{1/2}} ,
           \frac{\int_0^T \sqrt{Y_s} \, \dd W_s}
                {\left(\int_0^T Y_s \, \dd s\right)^{1/2}} ,
           \frac{\int_0^T \sqrt{Y_s} \, \dd B_s}
                {\left(\int_0^T Y_s \, \dd s\right)^{1/2}} ,
           \ee^{bT} Y_T , \ee^{bT} \int_0^T Y_s \, \dd s ,
           \int_0^T \frac{\dd s}{Y_s} \right) \\
    &\distr \left( Z_1, \bZ_2, \tcY_{-1/b}, -\frac{\tcY_{-1/b}}{b} ,
                   \int_0^{-1/b} \tcY_u \, \dd u \right)
     \qquad \text{as \ $T \to \infty$,}
  \end{aligned}
 \end{equation}
Applying Theorem \ref{THM_Zanten} for the continuous local martingale
 \ $M_t := \int_0^t \sqrt{Y_s} \,\dd W_s$, \ $t \in \RR_+$, \ with quadratic
 variation process \ $\langle M \rangle_t = \int_0^t Y_s \, \dd s$,
 \ $t \in \RR_+$, \ for \ $Q(t) := \ee^{bt/2}$, \ $t \in \RR_{++}$, \ and for
 \ $\bv := \left( V, -\frac{V}{b}, \int_0^\infty\frac{\dd s}{Y_s} \right)$
 \ (defined also on \ $(\Omega, \cF, \PP)$), \ we obtain
 \[
   \left( \ee^{bt/2} \int_0^t \sqrt{Y_s} \, \dd W_s ,
          V, -\frac{V}{b}, \int_0^\infty\frac{\dd s}{Y_s} \right)
   \distr
   \left( \left(-\frac{V}{b}\right)^{1/2} \xi_2,
         V, -\frac{V}{b}, \int_0^\infty\frac{\dd s}{Y_s} \right)
 \]
 as \ $t \to \infty$, \ where \ $\xi_2$ \ is a standard normally distributed
 random variable independent of \ $V$ \ and \ $\int_0^\infty\frac{\dd s}{Y_s}$.
\ Indeed, by \eqref{lim_intY}, we have
 \ $\ee^{bt} \langle M \rangle_t = \ee^{bt} \int_0^t Y_s \, \dd s
    \as - \frac{V}{b}$ \ as \ $t \to \infty$.
Here, by Ben Alaya and Kebaier \cite[Theorem 4]{BenKeb2},
 \[
 \left( \left(-\frac{V}{b}\right)^{1/2} \xi_2,
         V, -\frac{V}{b}, \int_0^\infty\frac{\dd s}{Y_s} \right)
 \distre
   \left( \left(-\frac{\tcY_{-1/b}}{b}\right)^{1/2} Z_2 ,
          \tcY_{-1/b}, -\frac{\tcY_{-1/b}}{b} ,
                   \int_0^{-1/b} \tcY_u \, \dd u \right),
 \]
 where \ $Z_2$ \ is a standard normally distributed random variable
 independent of \ $\tcY_{-1/b}$ \ and \ $\int_0^{-1/b} \tcY_u \, \dd u$.
\ By \eqref{lim_Y}, \eqref{lim_intY}, \eqref{lim_1/Y} and
 Lemma \ref{LEMMA_KATAI_MOGYORODI}, we obtain
 \[
   \left( \ee^{bt/2} \int_0^t \sqrt{Y_s} \, \dd W_s ,
          \ee^{bt} Y_t , \ee^{bt} \int_0^t Y_s \, \dd s ,
          \int_0^t \frac{\dd s}{Y_s} \right)
   - \left( \ee^{bt/2} \int_0^t \sqrt{Y_s} \, \dd W_s ,
            V, -\frac{V}{b} , \int_0^\infty\frac{\dd s}{Y_s} \right)
   \stoch 0
 \]
 as \ $t \to \infty$, \ hence
 \begin{align*}
  &\left( \ee^{bT/2} \int_0^T \sqrt{Y_s} \, \dd W_s ,
         \ee^{bT} Y_T , \ee^{bT} \int_0^T Y_s \, \dd s ,
         \int_0^T \frac{\dd s}{Y_s} \right) \distr
 \end{align*}
 \begin{align*}
  &\distr
   \left( \left(-\frac{\tcY_{-1/b}}{b}\right)^{1/2} Z_2 ,
          \tcY_{-1/b}, -\frac{\tcY_{-1/b}}{b} ,
                   \int_0^{-1/b} \tcY_u \, \dd u \right)
                 \qquad \text{as \ $T \to \infty$.}
 \end{align*}
Applying continuous mapping theorem, since \ $\PP(\tcY_{-1/b}\in\RR_{++})=1$,
 \ we obtain
 \begin{equation}\label{aabb_super_2}
   \left( \frac{\int_0^T \sqrt{Y_s} \, \dd W_s}
               {\left(\int_0^T Y_s \, \dd s\right)^{1/2}} ,
          \ee^{bT} Y_T , \ee^{bT} \int_0^T Y_s \, \dd s ,
          \int_0^T \frac{\dd s}{Y_s} \right)
   \distr
   \left( Z_2 , \tcY_{-1/b}, -\frac{\tcY_{-1/b}}{b} ,
                   \int_0^{-1/b} \tcY_u \, \dd u \right)
 \end{equation}
 as \ $T \to \infty$, \ hence we derived joint convergence of four coordinates
 of the left hand side of \eqref{aabb_super_1}.

We will prove \eqref{aabb_super_1} using continuity theorem.
Applying \eqref{Heston_SDE}, one can write
 \[
   \sigma_1 \int_0^T \sqrt{Y_s} \, \dd W_s
   = Y_T - y_0 - \int_0^T (a - b Y_s) \, \dd s , \qquad T \in \RR_{++} ,
 \]
 hence \ $\int_0^T \sqrt{Y_s} \, \dd W_s $ \ is measurable with respect to
 the $\sigma$-algebra \ $\sigma(Y_s , s \in [0, T])$.
\ For all \ $(u_1, u_2, u_3, v_1, v_2, v_3) \in \RR^6$ \ and \ $T \in \RR_{++}$,
 \ we have
 \begin{align*}
  &\EE\Biggl(\exp\Biggl\{\ii u_1
                         \frac{\int_0^T \frac{\dd B_s}{\sqrt{Y_s}}}
                              {\left(\int_0^T \frac{\dd s}{Y_s}\right)^{1/2}}
                         + \ii u_2
                         \frac{\int_0^T \sqrt{Y_s} \, \dd W_s}
                              {\left(\int_0^T Y_s \, \dd s\right)^{1/2}}
                         + \ii u_3
                           \frac{\int_0^T \sqrt{Y_s} \, \dd B_s}
                                {\left(\int_0^T Y_s \, \dd s\right)^{1/2}} \\
  &\phantom{\EE\Biggl(\exp\Biggl\{}
                         + \ii v_1 \ee^{bT} Y_T
                         + \ii v_2 \ee^{bT} \int_0^T Y_s \, \dd s
                         + \ii v_3 \int_0^T \frac{\dd s}{Y_s}\Biggr\}
             \Bigg| \, Y_s , s \in [0, T] \Biggr) \\
  &= \exp\Biggl\{\ii u_2
                 \frac{\int_0^T \sqrt{Y_s} \, \dd W_s}
                      {\left(\int_0^T Y_s \, \dd s\right)^{1/2}}
                 + \ii v_1 \ee^{bT} Y_T
                 + \ii v_2 \ee^{bT} \int_0^T Y_s \, \dd s
                 + \ii v_3 \int_0^T \frac{\dd s}{Y_s}\Biggr\} \\
  &\quad
     \times
     \EE\Biggl(\exp\Biggl\{\ii
                           \int_0^T
                            \Biggl( \frac{u_1}
                                  {\left(\int_0^T \frac{\dd t}{Y_t}\right)^{1/2}}
                                    \frac{1}{\sqrt{Y_s}}
                                    + \frac{u_3}
                                 {\left(\int_0^T Y_t \, \dd t\right)^{1/2}}
                                 \sqrt{Y_s} \Biggr)
                            \dd B_s\Biggr\}
               \Bigg| \, Y_s , s \in [0, T] \Biggr) \\
  &= \exp\Biggl\{\ii u_2
                 \frac{\int_0^T \sqrt{Y_s} \, \dd W_s}
                      {\left(\int_0^T Y_s \, \dd s\right)^{1/2}}
                 + \ii v_1 \ee^{bT} Y_T
                 + \ii v_2 \ee^{bT} \int_0^T Y_s \, \dd s
                 + \ii v_3 \int_0^T \frac{\dd s}{Y_s}\Biggr\} \\
  &\quad
     \times
     \exp\Biggl\{ - \frac{1}{2}
                    \int_0^T
                     \Biggl( \frac{u_1}
                               {\left(\int_0^T \frac{\dd t}{Y_t}\right)^{1/2}}
                               \cdot \frac{1}{\sqrt{Y_s}}
                             + \frac{u_3}
                                  {\left(\int_0^T Y_t \, \dd t\right)^{1/2}}
                             \sqrt{Y_s} \Biggr)^2 \dd s
         \Biggr\} \\
  &= \exp\Biggl\{\ii u_2
                 \frac{\int_0^T \sqrt{Y_s} \, \dd W_s}
                      {\left(\int_0^T Y_s \, \dd s\right)^{1/2}}
                 + \ii v_1 \ee^{bT} Y_T
                 + \ii v_2 \ee^{bT} \int_0^T Y_s \, \dd s
                 + \ii v_3 \int_0^T \frac{\dd s}{Y_s}\Biggr\} \\
  &\quad
     \times
     \exp\Biggl\{ - \frac{1}{2} (u_1^2 + u_3^2)
                  - \frac{T u_1 u_3}
                         {\left(\int_0^T Y_t \, \dd t
                                \int_0^T \frac{\dd t}{Y_t}\right)^{1/2}}
         \Biggr\} ,
 \end{align*}
 where we used the independence of \ $Y$ \ and \ $B$.
\ Consequently, the characteristic function of the random vector on the left
 hand side of \eqref{aabb_super_1} takes the form
 \begin{align*}
  &\EE\Biggl(\exp\Biggl\{\ii u_1
                         \frac{\int_0^T \frac{\dd B_s}{\sqrt{Y_s}}}
                              {\left(\int_0^T \frac{\dd s}{Y_s}\right)^{1/2}}
                         + \ii u_2
                           \frac{\int_0^T \sqrt{Y_s} \, \dd W_s}
                                {\left(\int_0^T Y_s \, \dd s\right)^{1/2}}
                         + \ii u_3
                           \frac{\int_0^T \sqrt{Y_s} \, \dd B_s}
                                {\left(\int_0^T Y_s \, \dd s\right)^{1/2}} \\
  &\phantom{\EE\Biggl(\exp\Biggl\{}
                         + \ii v_1 \ee^{bT} Y_T
                         + \ii v_2 \ee^{bT} \int_0^T Y_s \, \dd s
                         + \ii v_3 \int_0^T \frac{\dd s}{Y_s}\Biggr\}\Biggr) \\
  &= \ee^{-(u_1^2+u_3^2)/2}
     \EE\Biggl(\exp\Biggl\{\xi_T(u_2, v_1, v_2, v_3)
                           - \frac{T u_1 u_3}
                                  {\left(\int_0^T Y_t \, \dd t
                                         \int_0^T \frac{\dd t}{Y_t}\right)^{1/2}}
                   \Biggr\}\Biggr) ,
 \end{align*}
 where
 \[
   \xi_T(u_2, v_1, v_2, v_3)
   := \ii u_2
      \frac{\int_0^T \sqrt{Y_s} \, \dd W_s}
           {\left(\int_0^T Y_s \, \dd s\right)^{1/2}}
      + \ii v_1 \ee^{bT} Y_T
      + \ii v_2 \ee^{bT} \int_0^T Y_s \, \dd s
      + \ii v_3 \int_0^T \frac{\dd s}{Y_s} .
 \]
By \eqref{aabb_super_2}, for all \ $(u_2, v_1, v_2, v_3)\in\RR^4$,
 \begin{equation}\label{xiconv}
  \begin{split}
  &\EE(\exp\{\xi_T(u_2, v_1, v_2, v_3)\})\\
  &\qquad \to \EE\Biggl(\exp\Biggl\{ \ii u_2 Z_2
                         + \ii v_1 \tcY_{-1/b}
                         + \ii v_2 \left(-\frac{\tcY_{-1/b}}{b}\right)
                         + \ii v_3 \int_0^{-1/b} \tcY_u \, \dd u\Biggr\}\Biggr)
  \end{split}
 \end{equation}
 as \ $T \to \infty$.
\ Using \ $|\exp\{\xi_T(u_2, v_1, v_2, v_3)\}| = 1$, \ we have
 \begin{align*}
  &\Biggl| \EE\Biggl(\exp\Biggl\{ \xi_T(u_2, v_1, v_2, v_3)
                           - \frac{T u_1 u_3}
                                  {\left(\int_0^T Y_t \, \dd t
                                        \int_0^T \frac{\dd t}{Y_t}\right)^{1/2}}
                         \Biggr\}\Biggr)
         - \EE(\exp\{\xi_T(u_2, v_1, v_2, v_3)\}) \Biggr| \\
  &\leq \EE\Biggl( |\exp\{\xi_T(u_2, v_1, v_2, v_3)\}|
                 \Biggl|\exp\Biggl\{-\frac{T u_1 u_3}
                                          {\left(\int_0^T Y_t \, \dd t
                                        \int_0^T \frac{\dd t}{Y_t}\right)^{1/2}}
                            \Biggr\}
                        -1\Biggr|\Biggr) \\
  &= \EE\Biggl(\Biggl|\exp\Biggl\{-\frac{T u_1 u_3}
                                           {\left(\int_0^T Y_t \, \dd t
                                        \int_0^T \frac{\dd t}{Y_t}\right)^{1/2}}
                             \Biggr\}
                         -1\Biggr|\Biggr)
   \to 0 \qquad \text{as \ $T \to \infty$,}
 \end{align*}
 by dominated convergence theorem, since, by \eqref{lim_intY} and
 \eqref{lim_1/Y},
 \[
   \exp\Biggl\{-\frac{T u_1 u_3}
                {\left(\int_0^T Y_t \, \dd t
                       \int_0^T \frac{\dd t}{Y_t}\right)^{1/2}}\Biggr\}
   - 1
   = \exp\Biggl\{-\frac{T \ee^{bT/2} u_1 u_3}
                {\left(\ee^{bT} \int_0^T Y_t \, \dd t
                       \int_0^T \frac{\dd t}{Y_t}\right)^{1/2}}\Biggr\}
     - 1
   \as 0 \qquad \text{as \ $T \to \infty$,}
 \]
 and, by Cauchy--Schwarz inequality,
 \[
   \Biggl|\exp\Biggl\{-\frac{T u_1 u_3}
                            {\left(\int_0^T Y_t \, \dd t
                                   \int_0^T \frac{\dd t}{Y_t}\right)^{1/2}}
              \Biggr\}
          - 1\Biggr|
   \leq \exp\Biggl\{\frac{T |u_1 u_3|}
                         {\left(\int_0^T Y_t \, \dd t
                                \int_0^T \frac{\dd t}{Y_t}\right)^{1/2}}\Biggr\}
        + 1
   \leq \exp\{|u_1 u_3|\} + 1
 \]
 for all \ $T \in \RR_{++}$.
\ Using \eqref{xiconv}, we conclude
 \begin{align*}
  &\EE\Biggl(\exp\Biggl\{\ii u_1
                         \frac{\int_0^T \frac{\dd B_s}{\sqrt{Y_s}}}
                              {\left(\int_0^T \frac{\dd s}{Y_s}\right)^{1/2}}
                         + \ii u_2
                           \frac{\int_0^T \sqrt{Y_s} \, \dd W_s}
                                {\left(\int_0^T Y_s \, \dd s\right)^{1/2}}
                         + \ii u_3
                           \frac{\int_0^T \sqrt{Y_s} \, \dd B_s}
                                {\left(\int_0^T Y_s \, \dd s\right)^{1/2}} \\
  &\phantom{\EE\Biggl(\exp\Biggl\{}
                         + \ii v_1 \ee^{bT} Y_T
                         + \ii v_2 \ee^{bT} \int_0^T Y_s \, \dd s
                         + \ii v_3 \int_0^T \frac{\dd s}{Y_s}\Biggr\}\Biggr) \\
  &\to \ee^{-(u_1^2+u_3^2)/2}
       \EE\Biggl(\exp\Biggl\{ \ii u_2 Z_2
                              + \ii v_1 \tcY_{-1/b}
                              + \ii v_2 \left(- \frac{\tcY_{-1/b}}{b}\right)
                              + \ii v_3 \int_0^{-1/b} \tcY_u \, \dd u \Biggr\}
          \Biggr)
 \end{align*}
 as \ $T \to \infty$.
\ Note that, since \ $Z_2$ \ is independent of \ $\tcY_{-1/b}$ \ and
 \ $\int_0^{-1/b} \tcY_u \, \dd u$, \ we have
 \begin{align*}
  &\ee^{-(u_1^2+u_3^2)/2}
   \EE\Biggl(\exp\Biggl\{ \ii u_2 Z_2
                          + \ii v_1 \tcY_{-1/b}
                              + \ii v_2 \left(- \frac{\tcY_{-1/b}}{b}\right)
                              + \ii v_3 \int_0^{-1/b} \tcY_u \, \dd u \Biggr\}
      \Biggr) \\
  &= \EE(\ee^{\ii u_1 Z_1})\EE(\ee^{\ii u_2 Z_2}) \EE(\ee^{\ii u_3 Z_3})
     \EE\Biggl(\exp\Biggl\{ \ii v_1 \tcY_{-1/b}
                              + \ii v_2 \left(- \frac{\tcY_{-1/b}}{b}\right)
                              + \ii v_3 \int_0^{-1/b} \tcY_u \, \dd u \Biggr\} ,
 \end{align*}
 where \ $(Z_1, Z_3)$ \ is a 2-dimensional
 standard normally distributed random vector, independent of
 \ $(Z_2, \tcY_{-1/b}, \int_0^{-1/b} \tcY_u \, \dd u)$, \ hence we obtain
 \eqref{aabb_super_1} with \ $\bZ_2:=(Z_2, Z_3)$.

Finally, we prove \eqref{MLE_super_random}.
In a similar way, by \eqref{MLE-}, we have
 \[
   \left(\int_0^T Y_s \, \dd s\right)^{1/2} (\hb_T - b)
   = \frac{\frac{T \ee^{bT/2}}{\left(\ee^{bT} \int_0^T Y_s \, \dd s\right)^{1/2}}
           \frac{\sigma_1 \int_0^T \frac{\dd W_s}{\sqrt{Y_s}}}
                {\int_0^T \frac{\dd s}{Y_s}}
           - \sigma_1
             \frac{\int_0^T \sqrt{Y_s} \, \dd W_s}
                  {\left(\int_0^T Y_s \, \dd s\right)^{1/2}}}
          {1 - \frac{T^2 \ee^{bT}}
                    {\ee^{bT} \int_0^T Y_s \, \dd s
                     \int_0^T \frac{\dd s}{Y_s}}} ,
 \]
 and
 \[
   \left(\int_0^T Y_s \, \dd s\right)^{1/2} (\hbeta_T - \beta)
   = \frac{\frac{T \ee^{bT/2}}{\left(\ee^{bT} \int_0^T Y_s \, \dd s\right)^{1/2}}
           \left( \frac{\sigma_2 \varrho \int_0^T \frac{\dd W_s}{\sqrt{Y_s}}}
                       {\int_0^T \frac{\dd s}{Y_s}}
                  + \frac{\sigma_2 \sqrt{1 - \varrho^2}}
                         {\left(\int_0^T \frac{\dd s}{Y_s}\right)^{1/2}}
                    \frac{\int_0^T \frac{\dd B_s}{\sqrt{Y_s}}}
                         {\left(\int_0^T \frac{\dd s}{Y_s}\right)^{1/2}} \right)
           - \sigma_2
             \frac{\int_0^T \sqrt{Y_s} \, \dd \tW_s}
                  {\left(\int_0^T Y_s \, \dd s\right)^{1/2}}}
          {1 - \frac{T^2 \ee^{bT}}
                    {\ee^{bT} \int_0^T Y_s \, \dd s
                     \int_0^T \frac{\dd s}{Y_s}}} ,
 \]
 provided that \ $\int_0^T Y_s \, \dd s \int_0^T \frac{1}{Y_s} \, \dd s > T^2$
 \ which holds a.s.
By \eqref{lim_intY}, we get
 \[
   \frac{T\ee^{bT/2}}{ \left(\ee^{bT}\int_0^T Y_s \, \dd s\right)^{1/2} }
   \as \frac{0}{\left(-\frac{V}{b}\right)^{1/2}}
   = 0 \qquad \text{as \ $T \to \infty$,}
 \]
 hence \eqref{lim_intY}, \eqref{a_super}, \eqref{help8}, \eqref{help9},
 \eqref{aabbt_super}, \eqref{aabb_super_1},
 Slutsky's lemma, continuous mapping theorem and
 \ $\PP(\tcY_{-1/b}\in \RR_{++}) = 1$,
 \ $\PP(\int_0^{-1/b} \tcY_u \, \dd u \in\RR_{++}) = 1$ \ (due to
 \ $\PP(\tcY_t \in \RR_{++}, \, \forall \, t \in \RR_+) = 1$) \ yield the
 second statement.
Indeed,
 \begin{align*}
  &\begin{bmatrix}
    \ha_T - a \\
    \halpha_T - \alpha \\
    \left(\int_0^T Y_s \, \dd s\right)^{1/2} (\hb_T - b) \\
    \left(\int_0^T Y_s \, \dd s\right)^{1/2} (\hbeta_T - \beta)
   \end{bmatrix} \distr
 \end{align*}
 \begin{align*}
  &\distr
   \frac{1}{1 - \frac{0}{-\frac{\tcY_{-1/b}}{b}\int_0^{-1/b} \tcY_u \, \dd u}}
   \begin{bmatrix}
    \tcV
    - \frac{0}{\int_0^{-1/b} \tcY_u \, \dd u}
      \frac{1}{\left(-\frac{\tcY_{-1/b}}{b}\right)^{1/2}}
      (\bS^{1/2} \bZ_2)_1 \\[3mm]
    \varrho \frac{\sigma_2}{\sigma_1} \tcV
    + \frac{\sigma_2 \sqrt{1 - \varrho^2}}
           {\left( \int_0^{-1/b} \tcY_u \, \dd u \right)^{1/2}} Z_1
    - \frac{0}{\int_0^{-1/b} \tcY_u \, \dd u}
      \frac{1}{\left(-\frac{\tcY_{-1/b}}{b}\right)^{1/2}}
      (\bS^{1/2} \bZ_2)_2 \\[3mm]
    \frac{0}{\left(-\frac{\tcY_{-1/b}}{b}\right)^{1/2}} \tcV
    - (\bS^{1/2} \bZ_2)_1 \\[3mm]
    \frac{0}{\left(-\frac{\tcY_{-1/b}}{b}\right)^{1/2}}
    \left[ \varrho \frac{\sigma_2}{\sigma_1} \tcV
           + \frac{\sigma_2 \sqrt{1 - \varrho^2}}
                  {\left( \int_0^{-1/b} \tcY_u \, \dd u\right)^{1/2}} Z_1 \right]
    - (\bS^{1/2} \bZ_2)_2
   \end{bmatrix}
 \end{align*}
 as \ $T \to \infty$, \ where
 \ $\bS^{1/2} \bZ_2 = \bigl( (\bS^{1/2} \bZ_2)_1 , (\bS^{1/2} \bZ_2)_2 \bigr)^\top$.
\proofend

\begin{Rem}\label{Rem_comparison4}
Overbeck \cite[Theorem 3]{Ove} has already derived the asymptotic behaviour of
 \ $\hb_T$ \ with non-random and random scaling for supercritical CIR
 processes.
We also note that Ben Alaya and Kebaier \cite[Theorem 1, Case 3]{BenKeb1}
 described the asymptotic behavior of the MLE of \ $b$ \ for supercritical CIR
 processes supposing that \ $a \in \RR_{++}$ \ is known.
It turns out that in this case the limit distribution is different from that
 we have in \eqref{aabb_super}.
\proofend
\end{Rem}

\begin{Cor}\label{Cor_super}
Under the conditions of Theorem \ref{Thm_MLE_super}, the MLEs of \ $b$ \ and
 \ $\beta$ \ are weakly consistent, however, the MLEs of \ $a$ \ and \ $\alpha$
 \ are not weakly consistent.
(Recall also that earlier it turned out that the MLE of \ $b$ \ is in fact
 strongly consistent, see Theorem \ref{Thm_MLE_cons_super}.)
\end{Cor}

\noindent{\bf Proof.}
In order to show that the MLEs of \ $a$ \ and \ $\alpha$ \ are not weakly
 consistent, it suffices to show \ $\PP(\tcV \ne 0) > 0$, \ since \ $Z_1$ \ is
 independent of the random vector
 \ $(\tcY_{-1/b}, \int_0^{-1/b}\tcY_u \, \dd u)$, \ and
 \ $\PP\bigl(\int_0^{-1/b}\tcY_u \, \dd u > 0\bigr) = 1$ \ (see the end of
 Remark \ref{Thm_MLE_cons_sigma_rho}).
We have
 \[
   \PP(\tcV = 0)
   = \PP\biggl(\log \tcY_{-1/b} - \log y_0
               = \left(a - \frac{\sigma_1^2}{2}\right)
                 \int_0^{-1/b}\tcY_u \, \dd u\biggr)
   \leq \PP\bigl(\tcY_{-1/b} \geq y_0\bigr)
   < 1 ,
 \]
 where \ $y_0 \in \RR_{++}$.
\ Indeed, by Ikeda and Watanabe \cite[page 222]{IkeWat},
 \[
   \EE(\ee^{-\lambda\tcY_{-1/b}})
   = \left(1 + \frac{\sigma_1^2}{(-2b)}\lambda\right)^{-2a/\sigma_1^2} , \qquad
   \lambda \in \RR_+ ,
 \]
 hence \ $\tcY_{-1/b}$ \ has Gamma distribution with parameters
 \ $2a/\sigma_1^2$ \ and \ $-2b/\sigma_1^2$.
\proofend

\vspace*{3mm}

\appendix

\vspace*{5mm}

\noindent{\bf\Large Appendix}

\section{Limit theorems for continuous local martingales}

In what follows we recall some limit theorems for continuous local
 martingales.
We use these limit theorems for studying the asymptotic
 behaviour of the MLE of \ $(a, b, \alpha, \beta)$.
\ First we recall a strong law of large numbers for continuous local
 martingales.

\begin{Thm}{\bf (Liptser and Shiryaev \cite[Lemma 17.4]{LipShiII})}
\label{DDS_stoch_int}
Let \ $\bigl( \Omega, \cF, (\cF_t)_{t\in\RR_+}, \PP \bigr)$ \ be a filtered
 probability space satisfying the usual conditions.
Let \ $(M_t)_{t\in\RR_+}$ \ be a square-integrable continuous local martingale with respect to the
 filtration \ $(\cF_t)_{t\in\RR_+}$ \ such that \ $\PP(M_0 = 0) = 1$.
\ Let \ $(\xi_t)_{t\in\RR_+}$ \ be a progressively measurable process such that
 \[
   \PP\left( \int_0^t \xi_u^2 \, \dd \langle M \rangle_u < \infty \right) = 1 ,
   \qquad t \in \RR_+ ,
 \]
 and
 \begin{align}\label{SEGED_STRONG_CONSISTENCY2}
  \int_0^t \xi_u^2 \, \dd \langle M \rangle_u \as \infty \qquad
  \text{as \ $t \to \infty$,}
 \end{align}
 where \ $(\langle M \rangle_t)_{t\in\RR_+}$ \ denotes the quadratic variation
 process of \ $M$.
\ Then
 \begin{align}\label{SEGED_STOCH_INT_SLLN}
  \frac{\int_0^t \xi_u \, \dd M_u}
       {\int_0^t \xi_u^2 \, \dd \langle M \rangle_u} \as 0 \qquad
  \text{as \ $t \to \infty$.}
 \end{align}
If \ $(M_t)_{t\in\RR_+}$ \ is a standard Wiener process, the progressive
 measurability of \ $(\xi_t)_{t\in\RR_+}$ \ can be relaxed to measurability and
 adaptedness to the filtration \ $(\cF_t)_{t\in\RR_+}$.
\end{Thm}

The next theorem is about the asymptotic behaviour of continuous multivariate
 local martingales, see van Zanten \cite[Theorem 4.1]{Zan}.

\begin{Thm}{\bf (van Zanten \cite[Theorem 4.1]{Zan})}\label{THM_Zanten}
Let \ $\bigl( \Omega, \cF, (\cF_t)_{t\in\RR_+}, \PP \bigr)$ \ be a filtered
 probability space satisfying the usual conditions.
Let \ $(\bM_t)_{t\in\RR_+}$ \ be a $d$-dimensional square-integrable continuous local martingale
 with respect to the filtration \ $(\cF_t)_{t\in\RR_+}$ \ such that
 \ $\PP(\bM_0 = \bzero) = 1$.
\ Suppose that there exists a function \ $\bQ : \RR_+ \to \RR^{d \times d}$
 \ such that \ $\bQ(t)$ \ is an invertible (non-random) matrix for all
 \ $t \in \RR_+$, \ $\lim_{t\to\infty} \|\bQ(t)\| = 0$ \ and
 \[
   \bQ(t) \langle \bM \rangle_t \, \bQ(t)^\top \stoch \bfeta \bfeta^\top
   \qquad \text{as \ $t \to \infty$,}
 \]
 where \ $\bfeta$ \ is a \ $d \times d$ random matrix.
Then, for each $\RR^k$-valued random vector \ $\bv$ \ defined on
 \ $(\Omega, \cF, \PP)$, \ we have
 \[
   (\bQ(t) \bM_t, \bv) \distr (\bfeta \bZ, \bv) \qquad
   \text{as \ $t \to \infty$,}
 \]
 where \ $\bZ$ \ is a \ $d$-dimensional standard normally distributed random
 vector independent of \ $(\bfeta, \bv)$.
\end{Thm}

We note that Theorem \ref{THM_Zanten} remains true if the function \ $\bQ$
 \ is defined only on an interval \ $[t_0, \infty)$ \ with some
 \ $t_0 \in \RR_{++}$.

To derive consequences of Theorem \ref{THM_Zanten} one can use the following
 lemma which is a multidimensional version of Lemma 3 due to K\'atai and
 Mogyor\'odi \cite{KatMog}, see Barczy and Pap \cite[Lemma 3]{BarPap}.

\begin{Lem}\label{LEMMA_KATAI_MOGYORODI}
Let \ $(\bU_t)_{t\in\RR_+}$ \ be a $k$-dimensional stochastic process such that
 \ $\bU_t$ \ converges in distribution as \ $t \to \infty$.
\ Let \ $(\bV_t)_{t\in\RR_+}$ \ be an $\ell$-dimensional stochastic process such
 that \ $\bV_t \stoch \bV$ \ as \ $t \to \infty$, \ where \ $\bV$ \ is an
 \ $\ell$-dimensional random vector.
If \ $g: \RR^k \times \RR^\ell \to \RR^d$ \ is a continuous function, then
 \[
   g(\bU_t, \bV_t) - g(\bU_t, \bV) \stoch \bzero \qquad
   \text{as \ $t \to \infty$.}
 \]
\end{Lem}

\end{document}